\documentclass[11pt]{article}
\usepackage{longtable}

\usepackage{tocloft} 
\usepackage[caption = false]{subfig}\newcommand{\subfl}[1]{\subfloat[#1]}

\usepackage{amsmath}
\usepackage{amsfonts}
\usepackage[mathscr]{eucal}
\usepackage{amsthm}
\usepackage{amssymb}
\usepackage[noadjust]{cite}
\usepackage{enumerate}
\usepackage{tikz}\usetikzlibrary{matrix, calc, arrows}
\usepackage{hyperref}
\definecolor{myblue}{rgb}{0,0,0.6}
\hypersetup{pdftitle={Properties of IFS attractors with non-empty interiors},
     bookmarksopen=true, bookmarksopenlevel=4,
     colorlinks=true, linkcolor=myblue,  citecolor=myblue, filecolor=myblue,   urlcolor=myblue,  }

\usepackage{setspace}
\usepackage{color}
\usepackage{soul}

\graphicspath{{./}{figs/}}

\DeclareMathOperator{\supp}{supp}

\DeclareMathOperator{\diam}{diam}
\DeclareMathOperator{\dist}{dist}

\newcommand{\Cp}{\mathrm{Cap}}

\usepackage{bbm}
\usepackage{color}
\definecolor{amcol}{rgb}{0.8,0,0}
\definecolor{dhcol}{rgb}{0,0.5,0}

\begin{document}

\newcommand{\rf}[1]{(\ref{#1})}
\newcommand{\mmbox}[1]{\fbox{\ensuremath{\displaystyle{ #1 }}}}	

\newcommand{\hs}[1]{\hspace{#1mm}}
\newcommand{\vs}[1]{\vspace{#1mm}}

\newcommand{\ri}{{\mathrm{i}}}
\newcommand{\re}{{\mathrm{e}}}
\newcommand{\rd}{\mathrm{d}}

\newcommand{\R}{\mathbb{R}}
\newcommand{\Q}{\mathbb{Q}}
\newcommand{\N}{\mathbb{N}}
\newcommand{\Z}{\mathbb{Z}}
\newcommand{\C}{\mathbb{C}}
\newcommand{\K}{{\mathbb{K}}}

\newcommand{\cA}{\mathcal{A}}
\newcommand{\cB}{\mathcal{B}}
\newcommand{\cC}{\mathcal{C}}
\newcommand{\cS}{\mathcal{S}}
\newcommand{\cD}{\mathcal{D}}
\newcommand{\cH}{\mathcal{H}}
\newcommand{\cI}{\mathcal{I}}
\newcommand{\cItilde}{\tilde{\mathcal{I}}}
\newcommand{\cIhat}{\hat{\mathcal{I}}}
\newcommand{\cIcheck}{\check{\mathcal{I}}}
\newcommand{\cIstar}{{\mathcal{I}^*}}
\newcommand{\cJ}{\mathcal{J}}
\newcommand{\cM}{\mathcal{M}}
\newcommand{\cP}{\mathcal{P}}
\newcommand{\cV}{{\mathcal V}}
\newcommand{\cW}{{\mathcal W}}
\newcommand{\scrD}{\mathscr{D}}
\newcommand{\scrS}{\mathscr{S}}
\newcommand{\scrJ}{\mathscr{J}}
\newcommand{\sD}{\mathsf{D}}
\newcommand{\sN}{\mathsf{N}}
\newcommand{\sS}{\mathsf{S}}
 \newcommand{\sT}{\mathsf{T}}
 \newcommand{\sH}{\mathsf{H}}
 \newcommand{\sI}{\mathsf{I}}
 
\newcommand{\bs}[1]{\mathbf{#1}}
\newcommand{\bb}{\mathbf{b}}
\newcommand{\bd}{\mathbf{d}}
\newcommand{\bn}{\mathbf{n}}
\newcommand{\bp}{\mathbf{p}}
\newcommand{\bP}{\mathbf{P}}
\newcommand{\bv}{\mathbf{v}}
\newcommand{\bx}{\mathbf{x}}
\newcommand{\by}{\mathbf{y}}
\newcommand{\bz}{{\mathbf{z}}}
\newcommand{\bxi}{\boldsymbol{\xi}}
\newcommand{\boldeta}{\boldsymbol{\eta}}	

\newcommand{\ts}{\tilde{s}}
\newcommand{\tGamma}{{\tilde{\Gamma}}}
 \newcommand{\tbx}{\tilde{\bx}}
 \newcommand{\tbd}{\tilde{\bd}}
 \newcommand{\txi}{\xi}
 
\newcommand{\done}[2]{\dfrac{d {#1}}{d {#2}}}
\newcommand{\donet}[2]{\frac{d {#1}}{d {#2}}}
\newcommand{\pdone}[2]{\dfrac{\partial {#1}}{\partial {#2}}}
\newcommand{\pdonet}[2]{\frac{\partial {#1}}{\partial {#2}}}
\newcommand{\pdonetext}[2]{\partial {#1}/\partial {#2}}
\newcommand{\pdtwo}[2]{\dfrac{\partial^2 {#1}}{\partial {#2}^2}}
\newcommand{\pdtwot}[2]{\frac{\partial^2 {#1}}{\partial {#2}^2}}
\newcommand{\pdtwomix}[3]{\dfrac{\partial^2 {#1}}{\partial {#2}\partial {#3}}}
\newcommand{\pdtwomixt}[3]{\frac{\partial^2 {#1}}{\partial {#2}\partial {#3}}}
\newcommand{\bnabla}{\boldsymbol{\nabla}}
\newcommand{\dive}{\boldsymbol{\nabla}\cdot}
\newcommand{\curl}{\boldsymbol{\nabla}\times}
\newcommand{\Phixy}{\Phi(\bx,\by)}
\newcommand{\PhiOxy}{\Phi_0(\bx,\by)}
\newcommand{\dxPhixy}{\pdone{\Phi}{n(\bx)}(\bx,\by)}
\newcommand{\dyPhixy}{\pdone{\Phi}{n(\by)}(\bx,\by)}
\newcommand{\dxPhiOxy}{\pdone{\Phi_0}{n(\bx)}(\bx,\by)}
\newcommand{\dyPhiOxy}{\pdone{\Phi_0}{n(\by)}(\bx,\by)}

\newcommand{\eps}{\varepsilon}
\newcommand{\real}[1]{{\rm Re}\left[#1\right]} 
\newcommand{\im}[1]{{\rm Im}\left[#1\right]}
\newcommand{\ol}[1]{\overline{#1}}
\newcommand{\ord}[1]{\mathcal{O}\left(#1\right)}
\newcommand{\oord}[1]{o\left(#1\right)}
\newcommand{\Ord}[1]{\Theta\left(#1\right)}

\newcommand{\hsnorm}[1]{||#1||_{H^{s}(\bs{R})}}
\newcommand{\hnorm}[1]{||#1||_{\tilde{H}^{-1/2}((0,1))}}
\newcommand{\norm}[2]{\left\|#1\right\|_{#2}}
\newcommand{\normt}[2]{\|#1\|_{#2}}
\newcommand{\on}[1]{\Vert{#1} \Vert_{1}}
\newcommand{\tn}[1]{\Vert{#1} \Vert_{2}}

\newcommand{\xt}{\mathbf{x},t}
\newcommand{\PhiF}{\Phi_{\rm freq}}
\newcommand{\cone}{{c_{j}^\pm}}
\newcommand{\ctwo}{{c_{2,j}^\pm}}
\newcommand{\cthree}{{c_{3,j}^\pm}}

\newtheorem{thm}{Theorem}[section]
\newtheorem{lem}[thm]{Lemma}
\newtheorem{defn}[thm]{Definition}
\newtheorem{prop}[thm]{Proposition}
\newtheorem{cor}[thm]{Corollary}
\newtheorem{rem}[thm]{Remark}
\newtheorem{conj}[thm]{Conjecture}
\newtheorem{ass}[thm]{Assumption}
\newtheorem{example}[thm]{Example} 

\newcommand{\tH}{\widetilde{H}}
\newcommand{\Hze}{H_{\rm ze}} 	
\newcommand{\uze}{u_{\rm ze}}		
\newcommand{\dimH}{{\rm dim_H}}
\newcommand{\dimB}{{\rm dim_B}}
\newcommand{\IntClosOm}{\mathrm{int}(\overline{\Omega})}
\newcommand{\IntClosOmOne}{\mathrm{int}(\overline{\Omega_1})}
\newcommand{\IntClosOmTwo}{\mathrm{int}(\overline{\Omega_2})}
\newcommand{\Ccomp}{C^{\rm comp}}
\newcommand{\tCcomp}{\tilde{C}^{\rm comp}}
\newcommand{\uC}{\underline{C}}
\newcommand{\utC}{\underline{\tilde{C}}}
\newcommand{\oC}{\overline{C}}
\newcommand{\otC}{\overline{\tilde{C}}}
\newcommand{\capcomp}{{\rm cap}^{\rm comp}}
\newcommand{\Capcomp}{{\rm Cap}^{\rm comp}}
\newcommand{\tcapcomp}{\widetilde{{\rm cap}}^{\rm comp}}
\newcommand{\tCapcomp}{\widetilde{{\rm Cap}}^{\rm comp}}
\newcommand{\hcapcomp}{\widehat{{\rm cap}}^{\rm comp}}
\newcommand{\hCapcomp}{\widehat{{\rm Cap}}^{\rm comp}}
\newcommand{\tcap}{\widetilde{{\rm cap}}}
\newcommand{\tCap}{\widetilde{{\rm Cap}}}
\newcommand{\ccap}{{\rm cap}}
\newcommand{\ucap}{\underline{\rm cap}}
\newcommand{\uCap}{\underline{\rm Cap}}
\newcommand{\cCap}{{\rm Cap}}
\newcommand{\ocap}{\overline{\rm cap}}
\newcommand{\oCap}{\overline{\rm Cap}}
\DeclareRobustCommand
{\mathringbig}[1]{\accentset{\smash{\raisebox{-0.1ex}{$\scriptstyle\circ$}}}{#1}\rule{0pt}{2.3ex}}
\newcommand{\cirH}{\mathringbig{H}}
\newcommand{\cirHs}{\mathringbig{H}{}^s}
\newcommand{\cirHt}{\mathringbig{H}{}^t}
\newcommand{\cirHm}{\mathringbig{H}{}^m}
\newcommand{\cirHzero}{\mathringbig{H}{}^0}
\newcommand{\deO}{{\partial\Omega}}
\newcommand{\OO}{{(\Omega)}}
\newcommand{\Rn}{{(\R^n)}}
\newcommand{\Id}{{\mathrm{Id}}}
\newcommand{\gap}{\mathrm{Gap}}
\newcommand{\ggap}{\mathrm{gap}}
\newcommand{\isom}{{\xrightarrow{\sim}}}
\newcommand{\half}{{1/2}}
\newcommand{\mhalf}{{-1/2}}
\newcommand{\inter}{{\mathrm{int}}}

\newcommand{\Hsp}{H^{s,p}}
\newcommand{\Htq}{H^{t,q}}
\newcommand{\tHsp}{{{\widetilde H}^{s,p}}}
\newcommand{\SP}{\ensuremath{(s,p)}}
\newcommand{\Xsp}{X^{s,p}}

\newcommand{\dd}{{d}}\newcommand{\pp}{{p_*}}

\newcommand{\Rnn}{\R^{n_1+n_2}}
\newcommand{\Tr}{{\mathrm{Tr}}}

\renewcommand{\bs}[1]{\boldsymbol{#1}}
\newcommand{\be}{\bs{e}}
\renewcommand{\bn}{\bs{n}}
\renewcommand{\bx}{x}
\renewcommand{\by}{y}
\newcommand{\bbx}{\Psi}
\newcommand{\bby}{\widetilde{\Psi}}
\newcommand{\bg}{\bs{g}}
\newcommand{\bu}{\bs{u}}
\newcommand{\bw}{\bs{w}}
\newcommand{\bA}{\bs{A}}
\newcommand{\bC}{\bs{C}}
\newcommand{\bL}{\bs{L}}
\newcommand{\bS}{\bs{S}}
\newcommand{\bT}{\bs{T}}
\newcommand{\bU}{\bs{U}}
\newcommand{\bV}{\bs{V}}
\newcommand{\bX}{\bs{X}}
\newcommand{\bgamma}{{\bs{\gamma}}}
\newcommand{\bH}{\bs{H}}
\newcommand{\bnu}{\boldsymbol{\nu}}
\newcommand{\btau}{\boldsymbol{\tau}}
\newcommand{\bseta}{\boldsymbol{\eta}}
\newcommand{\bbeta}{{\boldsymbol{\beta}}}
\newcommand{\rD}{\mathrm{D}}
\newcommand{\rN}{\mathrm{N}}
\newcommand{\bm}{{\bs{m}}}
\newcommand{\bl}{\bs{l}}
\newcommand{\No}{{\mathbb{N}_0}}
\newcommand{\tOmega}{{\widetilde{\Omega}}}
\newcommand{\bj}{{\mathbf{j}}}
\newcommand{\tbX}{\tilde{\bX}}
\newcommand{\tA}{\widetilde{A}}
\newcommand{\tB}{\widetilde{B}}
\newcommand{\tF}{\widetilde{F}}
\renewcommand{\tH}{\widetilde{H}{}}
\newcommand{\tbH}{\widetilde{\bH}{}}
\newcommand{\sM}{\mathsf{M}}
\newcommand{\sE}{\mathsf{E}}
\newcommand{\cT}{\mathcal{T}}
\newcommand{\cU}{\mathcal{U}}
\newcommand{\cF}{\mathcal{F}}
\newcommand{\cL}{\mathcal{L}}
\newcommand{\cK}{\mathcal{K}}
\newcommand{\cN}{\mathcal{N}}
\newcommand{\cE}{\mathcal{E}}
\newcommand{\cR}{\mathcal{R}}
\newcommand{\tcA}{\tilde{\mathcal{A}}}
\newcommand{\tcL}{\tilde{\mathcal{L}}}
\newcommand{\tcK}{\tilde{\mathcal{K}}}
\newcommand{\bcD}{\boldsymbol{\mathcal{D}}}
\newcommand{\vbcU}{\vec{\boldsymbol{\cU}}}
\newcommand{\vbcE}{\vec{\boldsymbol{\cE}}}
\newcommand{\vbcS}{\vec{\boldsymbol{\cS}}}
\newcommand{\bscrS}{\boldsymbol{\scrS}}
\newcommand{\dudnjump}{\left[ \pdone{u}{n}\right]}
\newcommand{\dudnjumptext}{[ \pdonetext{u}{n}]}
\newcommand{\sumpm}[1]{\overbracket[0.5pt]{\underbracket[0.5pt]{\,#1\,}}}
\renewcommand{\dudnjump}{\left[ \pdone{u^s}{n}\right]}
\renewcommand{\dudnjumptext}{[ \pdonetext{u^s}{n}]}
\newcommand{\gradd}{\vec{\rm \bf grad}}	
\newcommand{\curld}{\vec{\rm \bf curl}}	
\newcommand{\dived}{{\rm div}}			
\newcommand{\gradp}{{\rm \bf grad}}	
\newcommand{\curlp}{{\rm \bf curl}}	
\newcommand{\scurlp}{{\rm curl}}		
\newcommand{\divep}{{\rm div}}			
\newcommand{\gradg}{{\rm \bf grad}_\Gamma}	
\newcommand{\curlg}{{\rm \bf curl}_\Gamma}	
\newcommand{\scurlg}{{\rm curl}_\Gamma}	
\newcommand{\diveg}{{\rm div}_\Gamma}		
\newcommand{\Deltag}{{\Delta_\Gamma}}		
\newcommand{\kap}{m}
\newcommand{\buk}{\bu_\kap}
\newcommand{\hbuk}{\hat{\bu}_\kap}
\newcommand{\bvk}{\bv_\kap}
\newcommand{\hbvk}{\hat{\bv}_\kap}
\newcommand{\bwk}{\bw_\kap}
\newcommand{\hbwk}{\hat{\bw}_\kap}
\newcommand{\hbu}{\hat{\bu}}
\newcommand{\hbv}{\hat{\bv}}
\newcommand{\hbw}{\hat{\bw}}
\newcommand{\hphi}{\hat{\phi}}
\newcommand{\hpsi}{\hat{\psi}}
\newcommand{\deG}{{\partial\Gamma}}
\newcommand{\GG}{(\Gamma)}
\newcommand{\tr}{\mathrm{tr}_\Gamma}
\newcommand{\trs}{\mathrm{tr}_{\Gamma,s}}
\newcommand{\trhalf}{\mathrm{tr}_{\Gamma,\frac{1}{2}}}
\newcommand{\IH}{\mathbb{H}}
\newcommand{\IS}{\mathbb{S}}
\newcommand{\IL}{\mathbb{L}}
\newcommand{\II}{\mathbb{I}}
\newcommand{\IV}{\mathbb{V}}
\newcommand{\IW}{\mathbb{W}}
\newcommand{\IC}{\mathbb{C}}
\newcommand{\IX}{\mathbb{X}}
\newcommand{\IP}{\mathbb{P}}
\newcommand{\IQ}{\mathbb{Q}}
\newcommand{\IY}{\mathbb{Y}}
\newcommand{\tf}{\tilde{f}}
\newcommand{\dn}{\partial_{\mathrm{n}}}
\newcommand{\ih}{\mathfrak{h}}
\newcommand{\Nn}{\mathbb{N}^{n}}%
\newcommand{\Non}{\mathbb{N}_{0}^{n}}
\newcommand{\Zn}{\mathbb{Z}^{n}}%
\newcommand{\Cn}{\mathbb{C}^{n}}%
\newcommand{\dual}[2]{\left\langle #1\,,\,#2\right\rangle} 
\newcommand{\wt}[1]{\widetilde{#1}}
\newcommand{\divepR}{{\rm div}_{\R^2}}
\newcommand{\scurlpR}{{\rm curl}_{\R^2}}

\definecolor{purple0}{rgb}{0.4,0,0.5}
\definecolor{orange}{rgb}{1,0.4,0}
\newcommand{\cu}[1]{{\color{purple0} #1 }}
\definecolor{orange0}{rgb}{1,0.3,0}
\newcommand{\ctodo}[1]{{\color{orange0} {\bf TODO:} #1 }}
\definecolor{green0}{rgb}{0.1,0.6,0}
\newcommand{\dnote}[1]{{\color{green0} {\bf DH:} #1 }}
\newcommand{\snote}[1]{{\color{red} {\bf SC:} #1 }}
\newcommand{\acnote}[1]{{\color{orange} {\bf AC:} #1 }}
\newcommand{\hBEM}{h_{\mathrm{BEM}}}
\newcommand{\hquad}{h_{\mathrm{quad}}}
\newcommand{\vb}{{\vec b}}
\newcommand{\vc}{{\vec c}}
\newcommand{\va}{{\vec a}}
\newcommand{\vphi}{{\vec\phi}}
\newcommand{\vpsi}{{\vec\psi}}
\newcommand{\vxi}{{\vec\xi}}
\newcommand{\dimA}{\dim_{\rm A}}

\allowdisplaybreaks[4]

\newtheorem{claim}[thm]{Claim}
\newtheorem{prob}[thm]{Problem}

\title{Properties of IFS attractors with non-empty interiors, related rough domains, and associated function spaces and scattering problems}
\author{A. M. Caetano$^{\text{a}}$,
S. N. Chandler-Wilde$^{\text{b}}$,
D. P. Hewett$^{\text{c}}$\\[8pt]
{\em 
This article is dedicated to Prof.~Hans Triebel on the Occasion of his 90th Birthday}
\footnote{{\bf Acknowledgements:} We would like to thank Prof.~Dr.~Hans Triebel for putting the first and third named authors in contact, inducing the fruitful collaboration which gave rise to this and other papers. DH acknowledges support from EPSRC grants EP/S01375X/1 and EP/V053868/1. SC-W and DH thank the Isaac Newton Institute for Mathematical Sciences for
support and hospitality during the programme “Mathematical Theory and Applications of Multiple Wave Scattering”, supported by EPSRC grant EP/R014604/1. AC is supported by CIDMA (https://ror.org/05pm2mw36)
under the Portuguese Foundation for Science and Technology 
(FCT, https://ror.org/00snfqn58), Grants 
UID/04106/2025 (https://doi.org/10.54499/UID/04106/2025)
and UID/PRR/04106/2025. We thank our collaborators Andrea Moiola (Pavia) and Andrew Gibbs (UCL) for useful discussions, and Lizaveta Ihnatsyeva (Kansas) for introducing us to the Aikawa and Assouad notions of fractal dimension.}\\[8pt]
$^{\text{a}}${\footnotesize
Center for R\&D in Mathematics and Applications, Departamento de Matem\'atica, Universidade de Aveiro, Aveiro, Portugal}\\ \footnotesize https://orcid.org/0000-0002-3655-4641\\
$^{\text{b}}${\footnotesize Department of Mathematics and Statistics, University of Reading, Reading, United Kingdom}\\ \footnotesize https://orcid.org/0000-0003-0578-1283\\
$^{\text{c}}${\footnotesize Department of Mathematics, University College London, London, United Kingdom}\\ \footnotesize https://orcid.org/0000-0003-3302-2567}

\maketitle
\renewcommand{\thefootnote}{\arabic{footnote}}

\begin{abstract}
We study properties of compact fractal sets $\Gamma\subset \R^n$ with 
non-empty interior $\Omega$, 
that are attractors of iterated function systems (IFSs) of contracting similarities satisfying the standard open set condition. 
Examples in the case $n=2$ are the closures of the Koch snowflake domain and the Gosper island domain. 
Our first main result is that 
$\Omega$ 
is {\em thick} in the sense of Triebel ({\em Function spaces and wavelets on domains}, EMS, 2008). 
An important consequence of this 
in the context of function spaces is that the set $C_0^\infty(\Omega)$ is dense in the Sobolev space $H^s_\Gamma:= \{\phi\in H^s(\R^n): \supp(\phi)\subset \Gamma\}$ for all $s\in\R$. 
Our second main result, accompanied by auxiliary results on pointwise multiplication by characteristic functions and uniform extension operators, is that the Sobolev spaces $\{H^s(\Omega)\}_{s\in \R}$, where $H^s(\Omega):=\{\phi|_\Omega: u\in H^s(\R^n)\}$, form an interpolation scale. 
This result is established as a special case of new extension and interpolation results for a range of Besov and Triebel-Lizorkin spaces, that apply to large classes of domains $\Omega$ that are thick and have boundary with Assouad dimension $<n$.
Our third main contribution is to prove best approximation error estimates in fractional negative-order Sobolev spaces for piecewise constant approximations on a ``fractal mesh'' of $\Gamma$, generated by the IFS, in which the mesh elements are self-similar copies of $\Gamma$. These estimates, similarly, are special cases of general new results, of independent interest, on piecewise constant approximation on rough domains with arbitrary meshes.
As an application of our results 
we study 
sound-soft acoustic scattering in $\R^{n+1}$ by the planar fractal screen $\Gamma\times \{0\}$. 
Using our density result 
we prove that 
the standard PDE formulation of this scattering problem is equivalent to the standard first kind boundary integral equation in which the boundary condition is imposed by restriction to the (relative) interior of the screen. To solve this equation we consider 
a piecewise-constant Galerkin boundary element method on a fractal mesh, 
and, using  
%
our best approximation error estimates, we prove convergence rates for the Galerkin approximation.





\end{abstract}
\tableofcontents

\section{Introduction}
\label{sec:Intro}


In this paper we address a number of problems involving function spaces defined on a class of subsets of $\R^n$ that we refer to as ``$n$-attractors''. 
Precisely (see Definition \ref{def:1} and \S\ref{sec:nsets} for details), an $n$-attractor $\Gamma\subset\R^n$ is the non-empty compact attractor of an iterated function system (IFS) satisfying the open set condition (OSC), such that $\Gamma$ has non-empty interior $\Gamma^\circ$ (and hence full Hausdorff dimension $\dimH(\Gamma)=n$). Such a $\Gamma$ possesses a certain ``self-similarity'', in the sense that $\Gamma$ can be written as a finite union of smaller copies of itself, the pairwise intersections of which have zero Lebesgue measure (see \S\ref{sec:nsets} for details). A collection of examples of $n$-attractors is given in Figure \ref{fig:2sets}. Typically, $n$-attractors have a fractal boundary, a classic example being the closure of the standard Koch snowflake domain, illustrated in Figure \ref{fig:2sets}(e), whose boundary is the union of three Koch curves. We note that $n$-attractors are related to the concept of ``fractal tilings'' 
--- see Remark \ref{rem:tiles} below for some references to the literature on this topic.  
We note also that $n$-attractors are examples of $n$-sets (see \eqref{eq:nset} for the definition) and their interiors are examples of open $n$-sets, often called interior regular domains (see, e.g., \cite[Definition~5.7]{caetano2019density}).

Our motivation for studying function spaces on $n$-attractors and related sets stems from their importance for the study of PDEs and integral equations in domains with fractal or other rough boundary (see, e.g., \cite{claret2024layer,hinz2023boundary,
lancia2002transmission,capitanelli2015weighted,FractalTransmission,WE}).
We shall explain this connection in more detail below, focussing in particular on the impact of our results for the mathematical and numerical analysis of acoustic scattering by fractal screens (studied recently in \cite{CoercScreen2,ScreenPaper,ScreenBEM,HausdorffBEM,HausdorffDomain} and related papers). However, we believe that our results will also be of independent interest to the function space and fractal analysis communities. In particular, as we make clear below, while our motivation in carrying out this work was to understand function spaces on $n$-attractors, many of our results apply to function spaces on rather general classes of domains with Lipschitz or rougher boundaries.


\subsection{Our main results}
We state and prove many of our function space results for general Besov and Triebel-Lizorkin spaces $B^s_{p,q}$ and $F^s_{p,q}$. However, in order to provide an accessible overview of our results, in this introduction we focus on their statements in the simpler setting of the 
{(fractional)} 
Sobolev spaces $H^s$ (i.e.\ Bessel potential spaces with integrability parameter $p=2$). As usual, for $s\in\R$ we let $H^s(\R^n)$ denote the set of tempered distributions $u$ for which $\|u\|_{H^s(\R^n)}^2:=\int_{\R^n} (1+|\xi|^2)^s |\hat{u}(\xi)|^2 \,\rd\xi<\infty$, where $\hat{u}$ denotes the Fourier transform of $u$. Given a non-empty closed set $E\subset \R^n$ we 
set $H^s_{E}:=\{u\in H^s(\R^n):\supp{u}\subset E\}$, and given a domain (i.e., a non-empty open set) $\Omega\subset \R^n$ we
define
$\tH^s(\Omega):=\overline{C_0^\infty(\Omega)}^{H^s(\R^n)}$ and we denote by $H^s(\Omega)$ the space of restrictions to $\Omega$ of elements of $H^s(\R^n)$, equipped with the quotient norm 
$\|u\|_{H^s(\Omega)} := \inf_{\varphi\in H^s\Rn,\,\varphi|_\Omega=u}\|\varphi\|_{H^s\Rn}$.
%
Restricted to this context, our main results are as follows for the $n$-attractor case: 

Let $\Gamma\subset\R^n$ be an $n$-attractor, and let $\Omega$ denote either the interior $\Gamma^\circ$ or the complement $\Gamma^c$. Then:
\begin{itemize}
\item[(R1)]
(Theorem \ref{thm:thick} and Corollary \ref{cor:Eporous}): 
$\Omega$ is thick, in the sense of Triebel (see Definition \ref{def:Thick} and \cite[Definition~3.1(ii)-(iv), Remark~3.2]{Tri08}), 
and $E$-porous, in the sense of Triebel (see Definition \ref{def:EPorous} and \cite[Definition~3.16(i)]{Tri08}).
\item[(R2)] 
{(Corollary \ref{cor:pointmult} and Proposition \ref{prop:extbyzero}):}
The characteristic function $\chi_\Omega$ of $\Omega$ is a pointwise multiplier on $H^s(\R^n)$ for 
{$|s|<\epsilon$,}  
{with $\epsilon=(n-\dimH(\partial\Omega))/2\in(0,1/2)$,  
where $\dimH(\partial\Omega)$ denotes the Hausdorff dimension of $\partial\Omega$}, 
which implies that ``extension by zero'' is well defined as a bounded operator from $H^s(\Omega)$ to $H^s(\R^n)$ 
{for such $s$.} 
{These statements fail for $s>\epsilon$.}
\item[(R3)]
(Corollary \ref{cor:tildesubscriptHs}): 
$\tH^s(\Omega)=H^s_{\overline\Omega}$ (i.e.\ $C^\infty_0(\Omega)$ 
is dense in 
$H^s_{\overline\Omega}$) for every $s\in \R$.
\item[(R4)] 
(Corollary
\ref{cor:Interpolationnset}):
$\{H^{s}(\Omega)\}_{s\in\R}$ and $\{\tH^{s}(\Omega)\}_{s\in\R}$ are interpolation scales ({in the sense of Definition \ref{d:IntScale}), with respect to both real and complex interpolation.} 
\end{itemize} 
Our results (R2)-(R4) apply to $n$-attractors, but also to a much larger class of domains. In relation to (R2) we prove as Proposition \ref{prop:pointmult} (for general Besov and Triebel-Lizorkin spaces) a version of Corollary \ref{cor:pointmult} which applies whenever $\Omega$ is a domain with $\dimA(\partial \Omega)<n$, where $\dimA(F)$ is the Aikawa/Assouad dimension of $F\subset \R^n$, as defined in Appendix \ref{sec:Appendix4}. Proposition \ref{prop:extbyzero} for (R2) has the same wider application. Corollary \ref{cor:tildesubscriptHs} in (R3) follows from a result of the same type for Besov/Triebel-Lizorkin spaces (Theorem \ref{cor:tildesubscriptGeneral}) that holds whenever $\Omega$ is thick, $\partial \Omega$ is bounded, and $\dimA(\partial \Omega)<n$. In relation to (R4), Corollary \ref{cor:Interpolationnset} is a consequence of Theorem \ref{prop:wolff}, a result on real and complex interpolation of Besov and Triebel-Lizorkin spaces that applies whenever $\Omega$ is thick, $\overline{\Omega}\neq \R^n$, and $\dimA(\partial \Omega)<n$. These more general results have in common that they all have application to rather general classes of domains with rough and fractal boundaries.

\subsection{Methodology}
Regarding methodology, our proof of the thickness result (R1) exploits the self-similar structure of $\Gamma$, which allows us to express $\Gamma$ as a union of scaled copies of $\Gamma$ of arbitrarily small diameter. 
To prove the pointwise multiplier result (R2) we appeal to a result from \cite{Sickel99a} (which is a special case of more general results proved in \cite{FrazierJawerth90}) which states that a sufficient condition for (R2) is that $\Omega$ belongs to the 
class $\cD^t$ (see \eqref{eq:D_t} for definition) for suitable $t$. We present two different proofs of the latter. The first, presented in 
{Theorem}
\ref{lem:newDt}, that applies in the $n$-attractor case but, more generally, whenever $\dimA(\partial \Omega)<n$, 
relates membership of the class $\cD^t$ to 
$\dimA(\partial \Omega)$. The second, presented in Lemma \ref{lem:Dt} in Appendix \ref{sec:Appendix3}, and tailored to the $n$-attractor case, is based on a more direct argument using the self-similar structure of $\Gamma$. 

In the Sobolev space setting the density result (R3) follows immediately from (R1) by \cite[Corollary 4.18]{caetano2019density}. In the general Besov/Triebel-Lizorkin setting, (R1) combined with \cite[Corollary 4.17]{caetano2019density} provides the analogous density result 
except in the case $s=0$, which we are able to handle separately (see Proposition \ref{prop:tilde-subscript}) using pointwise multiplier results generalising those in (R2) to the case where $\Omega$ is any domain with $\dimA(\partial \Omega)<n$ (see Proposition \ref{prop:pointmult}). 

For the proof of the interpolation scale result (R4) and its generalisation to the Besov/Triebel-Lizorkin setting we adapt the methodology of Triebel in  \cite[Theorem 4.17]{Tri08}, where similar results were proved under different assumptions on $\Omega$. 
The idea of Triebel's proof of \cite[Theorem 4.17]{Tri08} is to use restriction and extension operators to transfer standard interpolation properties of spaces on $\R^n$ to spaces on $\Omega$. A key component of this approach is to prove the existence of a common extension operator that acts in the same way for all $s$ in some open interval containing the range of $s$ for which one wishes to prove the interpolation result. While this may not be possible over the whole of $\R$, the result (R1) and results in \cite{Tri08} prove it is possible on any bounded subintervals of $(-\infty,0)$ and $(0,\infty)$, and our results in (R2) on extension by zero prove it 
for $s$ in {an} open interval $(-\epsilon,\epsilon)$ containing $0$. The fact that these three intervals overlap 
then permits the use of the so-called ``Wolff interpolation theorems'' from \cite{Wolff1982} to deduce interpolation results on the whole of $\R$. 

\subsection{Related literature}
Regarding comparisons with related literature, for (R1) we note that thickness of a family of snowflake domains, including the standard Koch snowflake shown in Figure \ref{fig:2sets}(e), was proved in \cite[Proposition~5.2]{caetano2019density} by a detailed geometrical analysis of a sequence of suitable polygonal prefractal approximations. Our result (R1) complements \cite[Proposition~5.2]{caetano2019density}, by proving thickness for general $n$-attractors, without the need for the construction of prefractal approximations. 
For (R2), our pointwise multiplier results complement and align with those presented for bounded Lipschitz domains in \cite[Remark 5.4]{Tri:02} and for certain non-Lipschitz domains in \cite{Triebel2003}, with a more detailed comparison being provided in Remark \ref{rem:multipliers}. 

{For (R3), 
a detailed exploration of conditions under which the equality $\tH^s(\Omega)=H^s_{\overline\Omega}$ does or does not hold can be found in \cite[\S3.5]{ChaHewMoi:13}.
In particular we note that $\tH^s(\Omega)=H^s_{\overline\Omega}$ for all $s\in\R$ if $\Omega$ is $C^0$ (see, e.g.,\ \cite[Theorem~3.20]{McLean}), and for a restricted range of $s$ 
if $\Omega$ is $C^0$ except at a closed, countable set of points of $\partial\Omega$ which has at most finitely many limit points in every bounded subset of $\partial\Omega$ \cite[Theorem~3.24]{ChaHewMoi:13}. 
As alluded to above, \cite[Corollary 4.18]{caetano2019density} shows that $\tH^s(\Omega)=H^s_{\overline\Omega}$ for all $s\in\R$ if $\Omega$ is thick {and $|\partial\Omega|=0$}. 

It was noted in \cite[Theorem 3.19]{ChaHewMoi:13} that 
$\tH^s(\Omega)\subsetneqq H^s_{\overline\Omega}$ for any $s\geq-n/2$ 
if 
$\Omega=B\setminus K$, where $B$ is an open ball and $K\subset B$ is any compact set with empty interior such that $H^{n/2}_{K}\neq \{0\}$ --- an example of such a $K$ was presented in \cite[Theorem~4]{Po:72a} (see also 
\cite[\S3.4]{HewMoi:15}). This example has the property that $\overline{\Omega}^\circ \neq \Omega$. Examples with $\overline{\Omega}^\circ = \Omega$ also exist. As noted in \cite[Lemma~3.20]{ChaHewMoi:13} (and the discussion before it),
results in \cite[\S11]{AdHe} imply that if $n\geq 2$ then for each $m\in \N$ there exists a bounded domain $\Omega\subset\R^n$ 
 with $\overline{\Omega}^\circ = \Omega$ such that $\tH^m(\Omega)\subsetneqq H^m_{\overline{\Omega}}$. 
Furthermore, by 
\cite[Lemma~3.26]{ChaHewMoi:13} 
the set $U=\overline{\Omega}^c$ then provides an example of a domain 
 $U\subset\R^n$ with $\overline{U}^\circ = U$ such that $\tH^{-m}(U)\subsetneqq {H^{-m}_{\overline{U}}}$. Moreover,
for any $R>\max_{x\in \partial\Omega}|x|$, 
$U_R:=\{x\in U:|x|<R\}$ is 
a bounded domain with $\overline{U_R}^\circ = U_R$ and $\tH^{-m}(U_R)\subsetneqq H^{-m}_{\overline{U_R}}$ (see Corollary \ref{cor:exnset} and its proof 
for more detail of these arguments).
This observation provides a negative answer to the speculative suggestion in \cite[\S3.5]{ChaHewMoi:13} (before Theorem 3.19) that $\tH^s(\Omega)= H^s_{\overline{\Omega}}$ might hold for all domains $\Omega\subset\R^n$ when $s<-n/2$. In Appendix \ref{sec:Appendix2} we extend the arguments of \cite[\S11]{AdHe} to Sobolev spaces of non-integral order, showing, for each $n\geq 2$ and $s\in [-n/2,n/2]\setminus\{0\}$, the existence of a bounded domain $\Omega\subset \R^n$ with $\overline{\Omega}^\circ=\Omega$ such that $\tH^s(\Omega)\subsetneqq H^s_{\overline{\Omega}}$ (Theorem \ref{thm:exnset} and Corollary \ref{cor:exnset}).
}

Finally, in relation to (R4), we note that function space interpolation is a well-studied topic with a large literature. However, as far as we are aware, our interpolation results {for Sobolev spaces on $n$-attractors} and our generalisations to Besov/Triebel-Lizorkin spaces and to the case that $\Omega$ is thick, $\overline{\Omega}\neq \R^n$, and $\dimA(\partial \Omega)<n$ are new. The closest similar results are those presented in Triebel \cite{Tri08}, in particular \cite[Theorem 4.17]{Tri08}, which provides similar interpolation results under the assumption that $\Omega$ is thick (as in our case), with $\partial\Omega$ being a $d$-set for some $n-1\leq d<n$ {(as is the case, e.g., for the Koch snowflake domain in Figure \ref{fig:2sets}(e)). Our results extend \cite[Theorem 4.17]{Tri08} to include cases where $\partial\Omega$ is fractal but not a $d$-set. Note that $\partial \Omega$ need not be a $d$-set even when $\Omega$ is the interior of an $n$-attractor, as shown in Remark \ref{rem:boundarydset}.} 
 

\subsection{Applications} \label{sec:app}
As key applications of our results within the current paper we highlight:
\begin{itemize}
\item[(A1)] Corollary \ref{cor:approx}, which provides new best approximation error {estimates} in fractional Sobolev spaces for approximation by piecewise-constant functions on a mesh of $\Omega$, the interior of the $n$-attractor $\Gamma$, comprising self-similar fractal elements, derived using the interpolation scale results in (R4).
\item[(A2)] Theorem \ref{thm:Closed}, which provides the first proof of well-posedness of the classical boundary value problem formulation (Problem \ref{prob:BVP}) for acoustic scattering in $\R^{n+1}$, $n=2,3$, by a sound-soft screen $\Gamma\subset\R^n\times\{0\}\cong \R^n$, for general $n$-attractors $\Gamma\subset\R^n$.
\item[(A3)] Theorem \ref{thm:ConvergenceRate}, which combines (A1) and (A2) to provide the first proof of  convergence rates for a piecewise-constant boundary element 
 method for the screen scattering problem for a general $n$-attractor screen under 
 solution regularity assumptions justified in Proposition \ref{rem:reg}.
\end{itemize}
We note that Corollary \ref{cor:approx} in (A1) is proved as a corollary of new general results on piecewise constant approximation on rough domains $\Omega$ with arbitrary meshes that will be of independent interest (Proposition \ref{prop:Approx} and Corollary \ref{cor:Interpolation}).

The results of the current paper are also crucial for the proofs of certain results in \cite{HausdorffDomain}, which studies acoustic scattering by general IFS attractors satisfying the OSC, including but not limited to the case considered in the current paper, where the scatterer is a subset of a hyperplane. 
However, the notation used in \cite{HausdorffDomain} differs slightly to that used here, and so for the benefit of readers of \cite{HausdorffDomain} we briefly comment {in Appendix \ref{sec:DomainComment}} on the connection between the two papers.

\subsection{Structure of the paper}
The structure of the paper is as follows. In \S\ref{sec:nsets} we recall some basic definitions relating to IFS attractors, and prove the thickness and $E$-porosity results mentioned in (R1) plus the fact that $\Omega$ is a member of the 
class $\cD^t$ for appropriate $t$. In \S\ref{sec:FunctionSpaces} we use the latter to prove the pointwise multiplier results mentioned in (R2), as well as the density and interpolation results mentioned in (R3) and (R4). At the end of this section we prove our best approximation error results for approximation by piecewise-constant functions on fractal meshes. In \S\ref{sec:BVPsBIEs} we apply the results of \S\S\ref{sec:nsets}-\ref{sec:FunctionSpaces} to the problem of acoustic scattering by a sound-soft planar screen, proving well-posedness of the scattering problem and convergence rates for a piecewise-constant boundary element method for its solution. 
{In Appendix \ref{sec:Notation} we provide a table of notation. 
In Appendix \ref{sec:DomainComment} we connect the results of the current paper with those of \cite{HausdorffDomain}.} 
In Appendices \ref{sec:Appendix}-\ref{sec:Appendix4} we collect technical definitions and auxiliary results that we use in the main paper and that may be of independent interest. Notably, in Appendix \ref{sec:Appendix4} we relate, in Proposition \ref{prop:upA}, Triebel's \cite{Tri:08,Tri08} notion of a uniformly porous compact set $F\subset \R^n$ to $\dimA(F)$.



\section{$d$-sets, $n$-sets, iterated function systems, and $n$-attractors}
\label{sec:nsets}


\subsection{Definitions, examples, and basic properties of $n$-attractors}
As in \cite[\S1.1]{JoWa84} and \cite[\S3]{Triebel97FracSpec}, given $0\leq d\leq n$ we say that a non-empty closed set $F\subset \R^n$ is a {\em$d$-set} (or an \textit{Ahlfors-David $d$-regular set}) if there exist constants $0<c_1\leq c_2$ such that
\begin{align}
\label{eq:dset}
c_{1}r^{d}\leq \cH^d(F\cap B(x,r))\leq c_{2}r^{d},\qquad x\in{F},\quad0<r\leq1,
\end{align}
where $B(x,r)\subset \R^n$ denotes the closed ball of radius $r$ centred on $x$ and $\cH^d$ denotes Hausdorff measure. 
Condition \rf{eq:dset} implies that $F$ is uniformly locally $d$-dimensional in the sense that $\dimH(F\cap B(x,r))=d$ for every $x\in {F}$ and $r>0$, where $\dimH$ denotes Hausdorff dimension.

Our particular focus in this paper is on the case $d=n$. We note that a non-empty closed set $F\subset\R^n$ is an {\em$n$-set} 
if and only if there exists $c_1>0$ such that
\begin{align}
\label{eq:nset}
c_{1}r^{n}\leq|F\cap B(x,r)|,\qquad x\in {F},\quad0<r\leq1,
\end{align}
where 
$|\cdot|$ 
denotes Lebesgue measure on $\R^n$.
If $\Omega\subset \R^n$ is a domain and there exists $c_1>0$ such that \eqref{eq:nset} holds with $F$ replaced by $\Omega$, then we call $\Omega\subset \R^n$ an {\em open $n$-set}.
If $\Omega$ is an open $n$-set then (see
	\cite[Proposition~1 on p.~205]{JoWa84}) $F := \overline{\Omega}$ is an $n$-set and $|\partial\Omega|=0$, so that also $|\partial F|=0$. Conversely, if $F\subset \R^{n}$ is an $n$-set and $|\partial F|=0$, it is easy to see that $F= \overline\Omega$, where $\Omega := F^\circ$, the interior of $F$, and that $\Omega$ is an open $n$-set. We note that open $n$-sets are also known as {\em interior regular domains},  see \cite[Definition 5.7]{caetano2019density}, or {\em $n$-thick domains}, e.g., \cite{Rychkov:00}, and that the condition \eqref{eq:nset}, applied when $F$ is replaced by an open set $\Omega$, is often termed the {\em measure density condition}, e.g., \cite{Zhou:15,HeIhTu:16}. 

By an iterated function system (IFS) of contracting similarities, we mean a collection $\{s_1,s_2,\ldots,s_M\}$, for some $M\geq 2$, where, for each $m=1,\ldots,M$, $s_m:\R^{n}\to\R^{n}$, with
$|s_m(x)-s_m(y)| = \rho_m|x-y|$, $x,y\in \R^{n}$,
for some $\rho_m\in (0,1)$.
It is standard that, given such an IFS, there exists a unique non-empty compact set $\Gamma\subset \R^{n}$, called the \emph{attractor} of the IFS, satisfying
\begin{equation} \label{eq:fixedfirst}
\Gamma = s(\Gamma), \quad \mbox{where} \quad s(E) := \bigcup_{m=1}^M  s_m(E), \quad  E\subset \R^{n}.
\end{equation}
{We define the {\em similarity dimension} of the IFS (following, e.g., \cite[\S9.2]{Fal}) to be the unique positive solution $d>0$ of the equation
\begin{align} \label{eq:dfirst}
\sum_{m=1}^M  (\rho_m)^d = 1,
\end{align}
which is no smaller than the Hausdorff dimension of $\Gamma$ \cite[Proposition 9.6]{Fal}.
If the \emph{open set condition (OSC)} \cite[(9.11)]{Fal} holds, i.e.\ there exists a non-empty bounded open set\footnote{We use here the version of the OSC from \cite{Fal}, where the open set $O$ is required to be bounded. The original OSC in Hutchinson \cite{hutchinson1981fractals} does not require boundedness of $O$, but  if $O$ satisfies the OSC condition in Hutchinson's sense and $\Gamma$ is the attractor of the IFS, then our OSC condition is satisfied by the bounded open set $O_1:=\{x\in O:\dist(x,\Gamma)<1\}$. Thus these two OSC versions are equivalent.} $O\subset \R^{n}$ such that
\begin{align} \label{oscfirst}
s(O) \subset O \quad \mbox{and} \quad s_m(O)\cap s_{m'}(O)=\emptyset, \quad m\neq m',
\end{align}
then $d\in (0,n]$, $\Gamma$ is a $d$-set, and $d$ coincides with the Hausdorff dimension of $\Gamma$ \cite[Theorem~4.7]{Triebel97FracSpec}.} 
We say that the IFS is \textit{homogeneous} if $\rho_1=\ldots = \rho_M = \rho$ for some $\rho\in(0,1)$. In that case, the similarity dimension satisfies the formula $d=\log{M}/\log(1/\rho)$.

As our focus is on the case where the OSC holds and $d=n$, for convenience we introduce some terminology for this case. 

\begin{defn}[``$n$-attractor''] 
\label{def:1}
If $\Gamma\subset\R^n$ is the non-empty compact attractor (in the sense of \rf{eq:fixedfirst}) of an IFS satisfying the OSC \eqref{oscfirst}, and \eqref{eq:dfirst} holds with $d=n$, then we say $\Gamma$ is an 
``$n$-attractor''. 
\end{defn}

\begin{rem}[Alternative characterisation of $n$-attractors] Suppose that $\Gamma$ is the attractor of an IFS of contracting similarities with similarity dimension $d\leq n$. Then $\Gamma$ is an $n$-attractor if and only if $\Gamma$ has interior points, alternatively, if and only if $|\Gamma|>0$. For if $\Gamma$ has interior points, which implies $|\Gamma|>0$, then $\dimH(\Gamma)=n$, so that $d=n$. And if $\Gamma$ is an attractor with similarity dimension $d=n$, then  $\Gamma$ is an $n$-attractor if and only if has interior points, alternatively if and only if $|\Gamma|>0$, by \cite[Theorem 2.2, Corollary 2.3]{Sc:94}.
\end{rem}

\begin{rem}
\label{rem:tiles}
The concept of an $n$-attractor is closely related to the concept of a ``fractal tiling'' --- see, e.g., \cite{bandt1991self,grochenig1994self,
strichartz1999geometry,Keesling:99,BaVi:14}. 
In particular\footnote{To see this claim
suppose that $\Gamma$ is a homogeneous $n$-attractor. Define $g:=s_1^{-1}$. Then $g(\Gamma)=\bigcup_{i=1}^M s_1^{-1}(s_i(\Gamma))$, and, since the IFS is homogeneous, the maps $s_1^{-1}\circ s_i$ are all congruences. Furthermore, the IFS satisfies the OSC with $O=\Gamma^\circ\neq \emptyset$ (see Proposition \ref{prop:nset}(\ref{viii})), so that the sets $\{s_1^{-1}(s_i(\Gamma))^\circ\}_{i=1}^M=\{s_1^{-1}(s_i(\Gamma^\circ))\}_{i=1}^M$ are mutually disjoint. Hence $\Gamma$ is a compact $M$-rep tile, since to say that a closed set $\Gamma_1$ is an $M$-rep tile, for some integer $M\geq 2$, as defined in \cite{bandt1991self}, means that $\Gamma_1^\circ\neq \emptyset$ and there are sets $\Gamma_2$, \ldots, $\Gamma_M$, congruent to $\Gamma_1$, such that the interiors of $\Gamma_1,\ldots,\Gamma_M$ are mutually disjoint, and a similarity mapping $g:\R^n\to \R^n$, such that $g(\Gamma_1)= \Gamma_1\cup\ldots \cup\Gamma_M$.
Conversely, suppose that $A_1$ is a compact $M$-rep tile, so that $A_1^\circ\neq \emptyset$ and $g(A_1) = \bigcup_{i=1}^M A_i$, for some similarity $g$ and sets $A_2,\ldots,A_M$ congruent to $A_1$, such that the interiors of $A_1,\ldots,A_M$ are mutually disjoint. Since $A_1^\circ\neq \emptyset$ and all the $A_i$ are congruent to $A_1$, it must be that $g$ is an expanding similarity. Hence $g^{-1}$ is a contracting similarity. Let $T_i:A_1\to A_i$ be congruences, for $i=1,\ldots, M$, with $T_1$ the identity map. For each $i$ define $s_i:=g^{-1}\circ T_i$. Then $\{s_1,\ldots,s_M\}$ is a homogeneous IFS of contracting similarities, of which $A_1$ is the unique attractor. Since $A_i^\circ\cap A_j^\circ=\emptyset$, for $i\neq j$, this IFS satisfies the OSC with $O=A_1^\circ$ and hence $A_1$ is a $d$-set. Finally, $d=n$ because $A_1^\circ\neq \emptyset$.}, the set of homogeneous $n$-attractors (i.e.,\ those for which all contraction factors are equal) coincides with the set of compact ``$M$-rep tiles'', as defined in \cite{bandt1991self}. 
\end{rem}

\begin{figure}[th!]
\centering
\subfl{A square}{\includegraphics[width = 32mm]{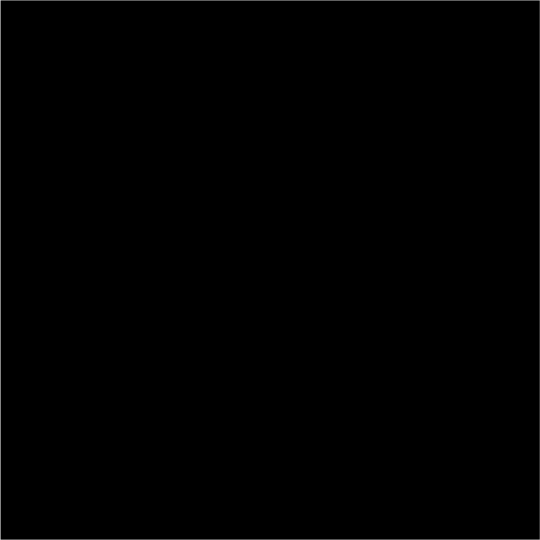}
\hspace{2mm}
\includegraphics[width = 32mm]{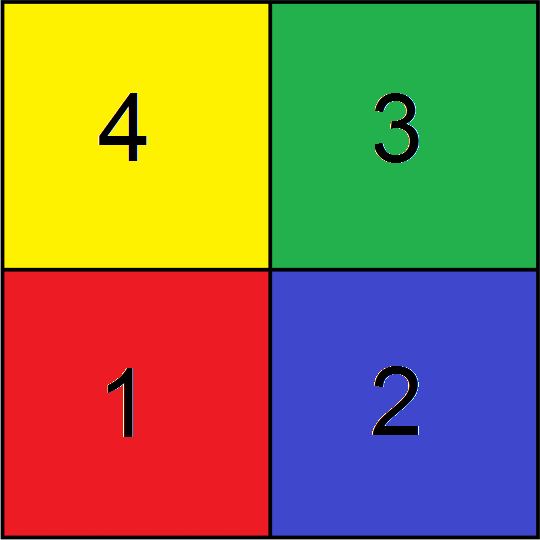}}
\hspace{6mm}
\subfl{Two squares touching at a point}{\includegraphics[width = 32mm]{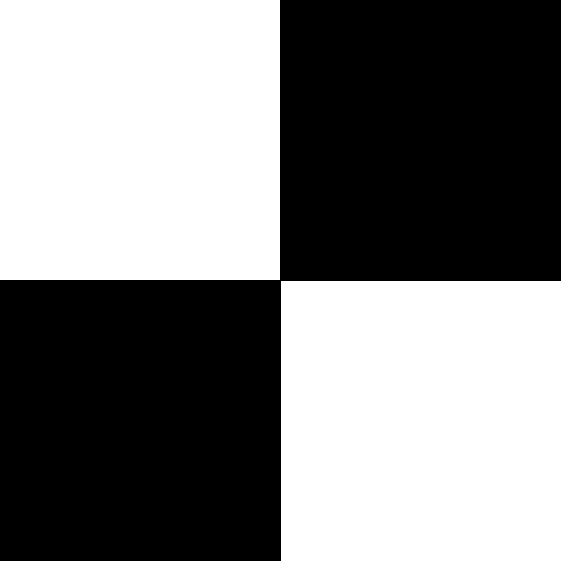}
\hspace{2mm}
\includegraphics[width = 32mm]{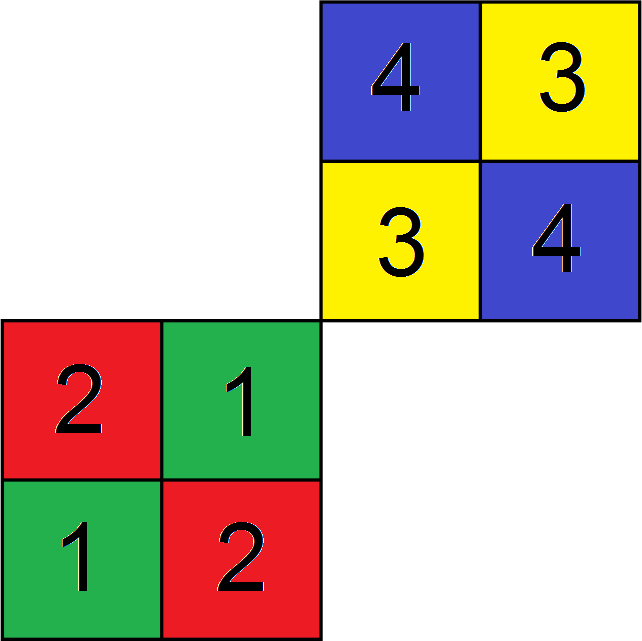}}
\hspace{2mm}
\subfl{Two separated squares}{\includegraphics[width = 42mm]{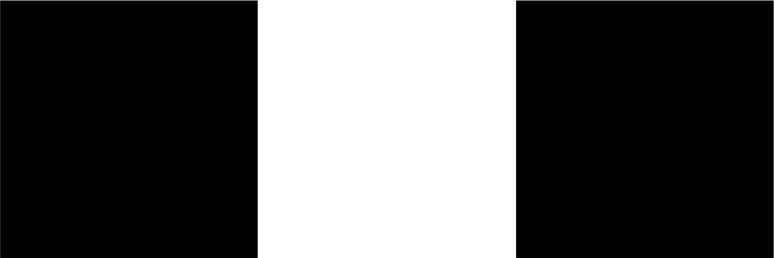}
\hspace{2mm}
\includegraphics[width = 42mm]{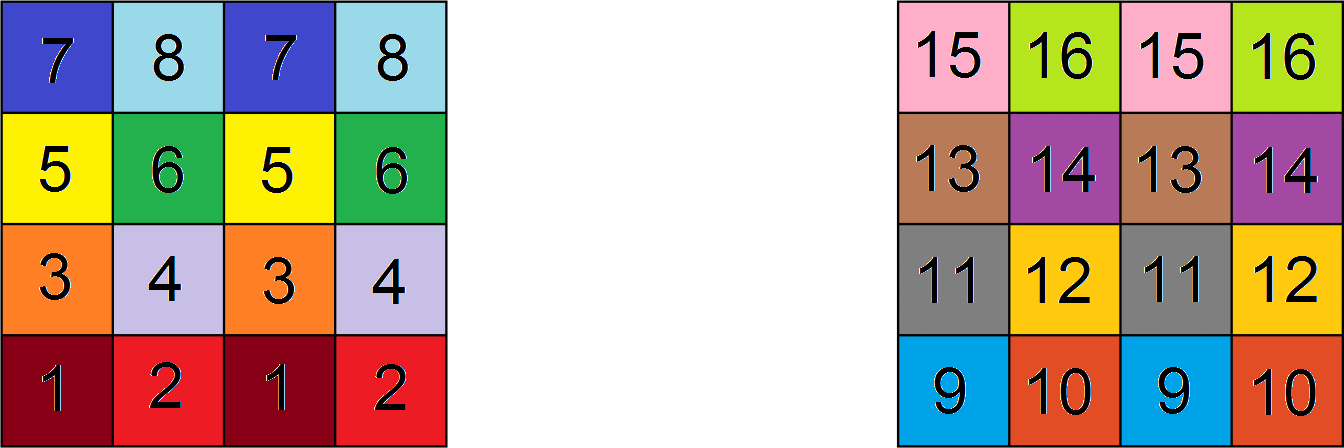}}
\hspace{4mm}
\subfl{An example with infinitely many components}{\includegraphics[width = 32mm]{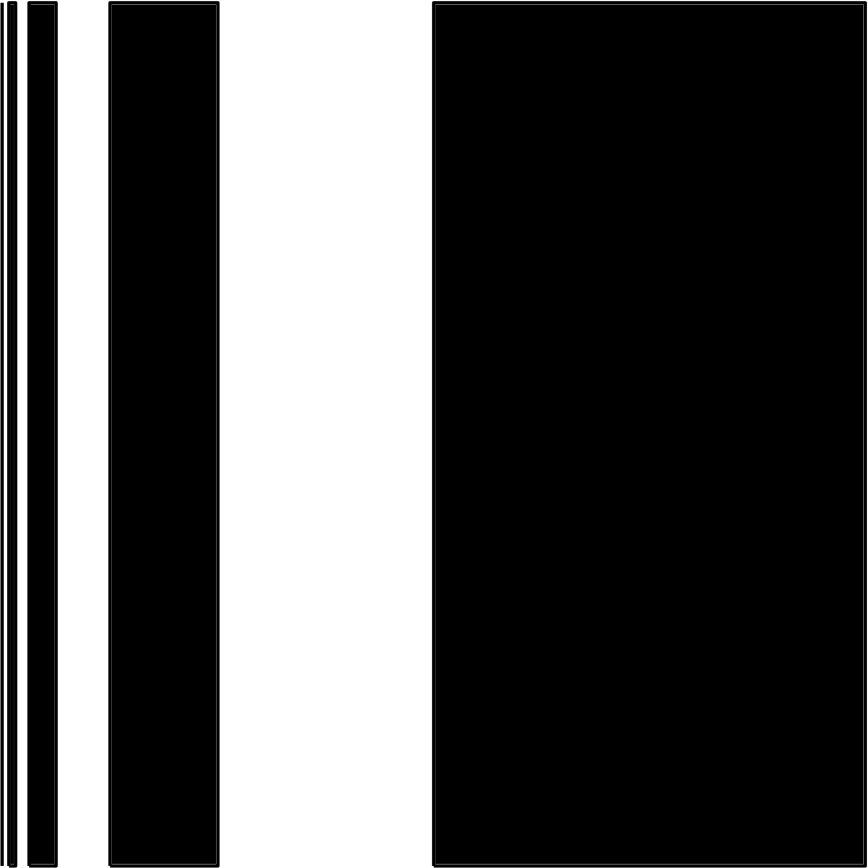}
\hspace{2mm}
\includegraphics[width = 32mm]{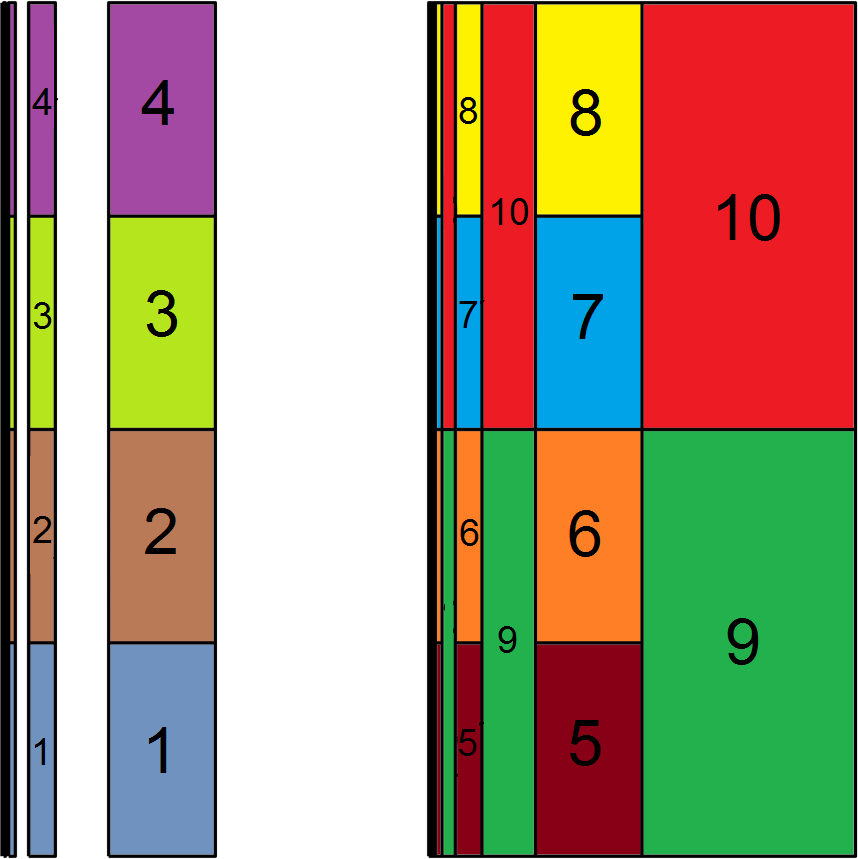}}
\hspace{4mm}
\subfl{Koch snowflake}{\includegraphics[width = 32mm]{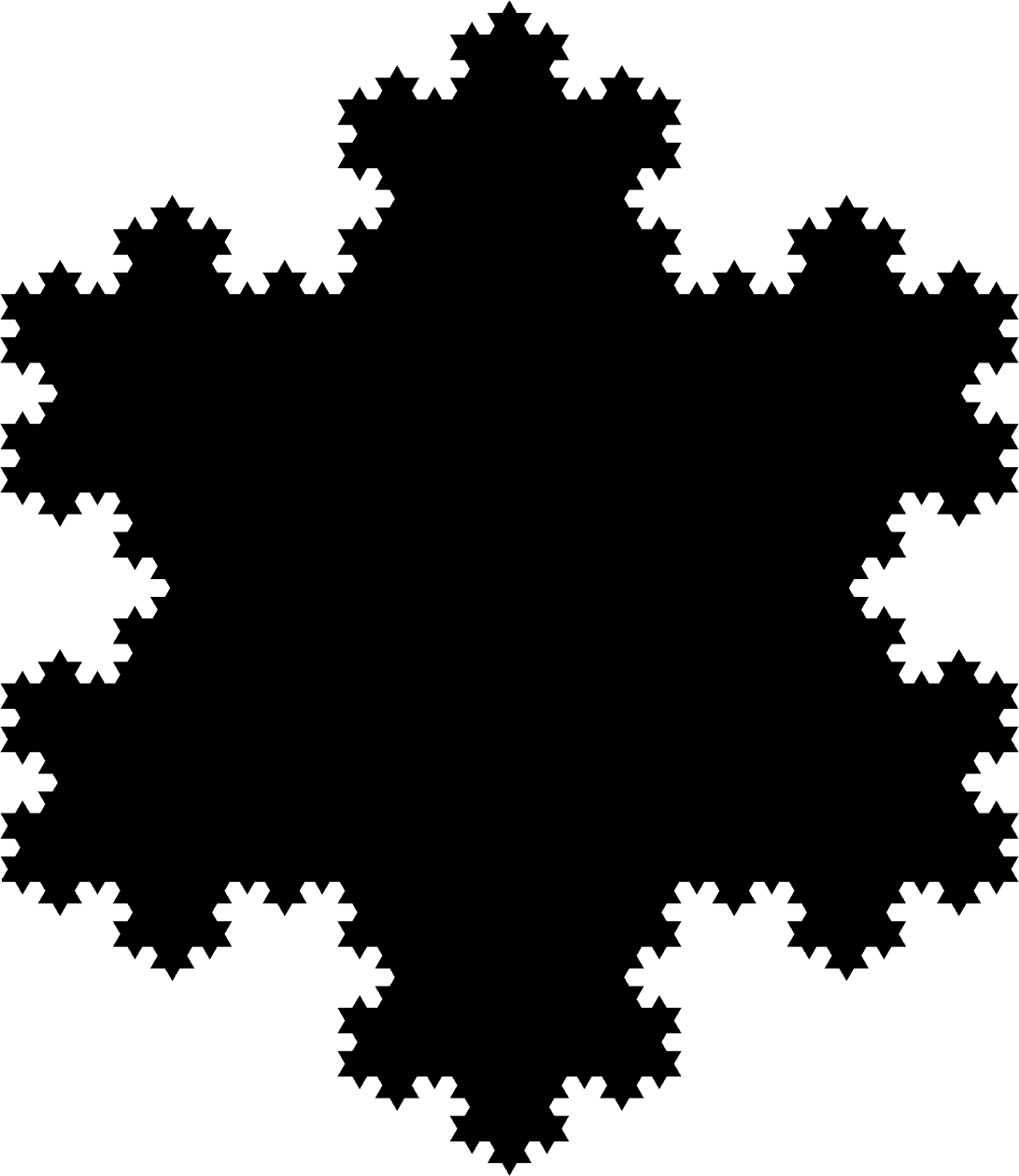}
\hspace{2mm}
\includegraphics[width = 32mm]{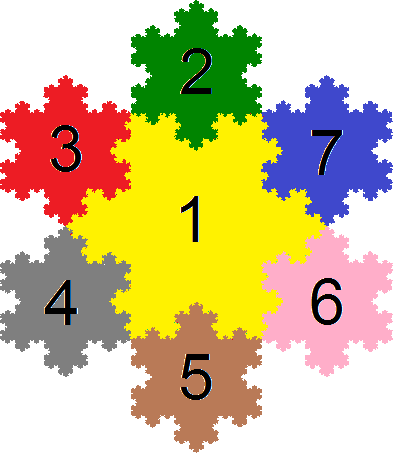}}
\hspace{4mm}
\subfl{Fudgeflake}{\includegraphics[width = 35mm]{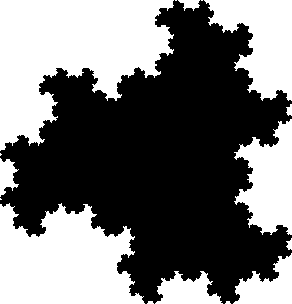}
\hspace{2mm}
\includegraphics[width = 35mm]{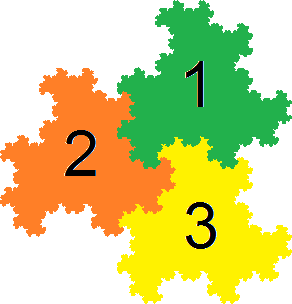}}
\hspace{4mm}
\subfl{Gosper island}{\includegraphics[width = 32mm]{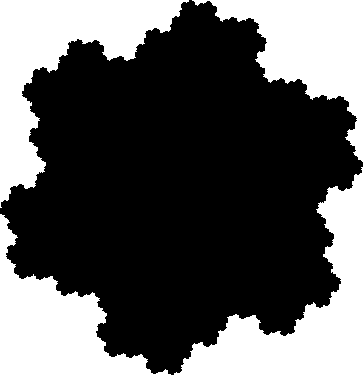}
\hspace{2mm}
\includegraphics[width = 32mm]{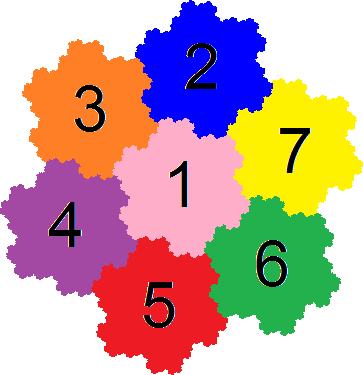}}
\hspace{4mm}
\subfl{Twindragon}{\includegraphics[width = 35mm]{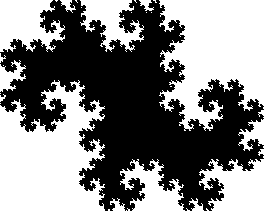}
\hspace{2mm}
\includegraphics[width = 35mm]{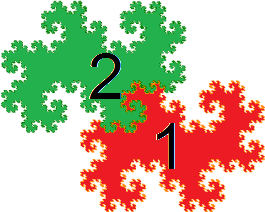}}
\hspace{4mm}
\subfl{Terdragon}{\includegraphics[width = 35mm]{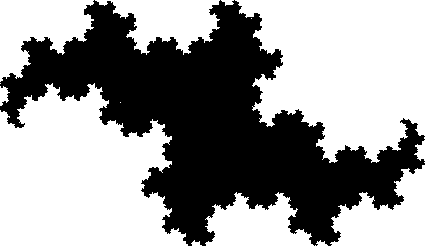}
\hspace{2mm}
\includegraphics[width = 35mm]{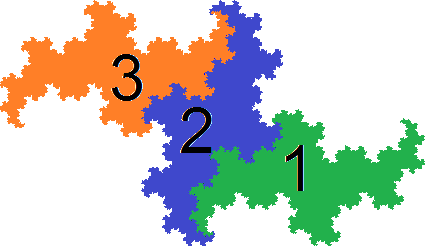}}
\hspace{4mm}
\subfl{Levy dragon}{\includegraphics[width = 35mm]{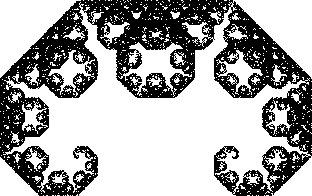}
\hspace{2mm}
\includegraphics[width = 35mm]{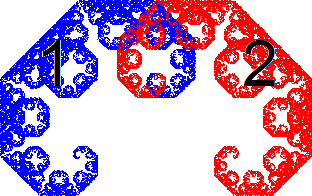}}
\caption{Examples of 
$n$-attractors, 
$\Gamma\subset\R^n$, 
in the case $n=2$. 
In each case the left-hand plot is $\Gamma$, the attractor of some IFS $\{s_1,\ldots,s_M\}$. The right-hand plot is a copy of $\Gamma$ with the subsets $s_1(\Gamma), \ldots s_M(\Gamma)$ numbered, and each shaded in a different colour. All except (d) and (e) are homogeneous, and hence examples of $M$-rep tiles in the sense of \cite{bandt1991self}.}
\label{fig:2sets}

\end{figure}

\begin{rem} [Examples of $n$-attractors]  \label{rem:exIFS} 

Given $0<\alpha<1$, a non-trivial example of an $n$-attractor in the case $n=1$ is
\begin{equation} \label{eq:Gammaeq}
\Gamma=\{0\}\cup \bigcup_{m=0}^\infty [\alpha^{2m+1},\alpha^{2m}],
\end{equation}
which is the attractor of $\{s_1,s_2,s_3\}$, where $s_1(x)=\alpha^2 x$, $s_2(x) = (1-\alpha)x + \alpha$, $s_3(x)=\alpha(1-\alpha)x + \alpha$. 
This example shows that $n$-attractors need not be connected. 

In Figure \ref{fig:2sets} we present a number of examples of $n$-attractors in the case $n=2$. 
Examples (a)-(c) are cases where $\Gamma$ is a finite union of polygons. Example (b) is a case where $\Gamma$ is connected but has disconnected interior, and example (c) is a case where $\Gamma$ is disconnected. 
Example (d) is an example where $\Gamma$ has infinitely many components. This example is the Cartesian product of \eqref{eq:Gammaeq}, for $\alpha=1/2$, with the interval $[0,1]$, and is the attractor of the IFS 
$\{s_1,\ldots,s_{10}\}$, where $s_m=\rho_mx+\delta_m$ for $m=1,\ldots,10$, with 
\[ \rho_1=\ldots=\rho_m=\frac{1}{4}, \qquad \rho_9=\rho_{10} = \frac{1}{2}, \qquad \delta_m=(0,\frac{m-1}{4}), \quad m=1,\ldots,4,\]
\[ \delta_m=(\frac{1}{2},\frac{m-5}{4}), \quad m=5,\ldots,8, \qquad
\delta_9 = (\frac{1}{2},0), \quad \delta_{10} = (\frac{1}{2},\frac{1}{2}).
\] 
Examples (e)-(j) have boundaries that are fractal, with Hausdorff dimensions in the interval $(1,2)$. 
In particular, example (e) is the closure of the standard Koch snowflake domain (e.g., \cite[\S2.1]{banjai2021poincar}). This is the attractor of an IFS with $M=7$; see \cite[Fig.~3]{HausdorffQuadrature}, 
\cite{RiddleWebSite}
and \cite{Burns94}, which argues using a standard inner-outer polygonal approximation, for which see, e.g., \cite[\S2.1]{banjai2021poincar} or \cite[\S5.1]{caetano2019density}. 
IFSs for examples (f)-(j) can be found at \cite{RiddleWebSite}. 
Example (j) is the Levy dragon, properties of which were studied, for example, in \cite{strichartz1999geometry,Keesling:99,DuKe:99,BaKiSt:02}. 
In particular, we note that,  while the Levy dragon is connected, its interior has countably many connected components, each of which has a polygonal boundary \cite{BaKiSt:02}.

\end{rem}

Proposition \ref{prop:nset} collects some fundamental facts about $n$-attractors. Before stating the proposition we 
{make some further preparations. 

First, we introduce some notation} that will allow us to conveniently describe certain decompositions of an $n$-attractor into similar copies of itself. 
Let $\Gamma\subset\R^n$ be an $n$-attractor, defined by some IFS $\{s_1,\ldots,s_M\}$. 
For $\ell\in \N$ we define the set of multi-indices $I_\ell:=\{1,\ldots,M\}^\ell\!=\{\bm = (m_1,m_2,\ldots,m_\ell):1\leq m_l\leq M , \text{ for } l=1,2,\ldots,\ell\}$, and for $E\subset\R^n$ and $\bm\in I_\ell$ we define
$E_{\bm}:=s_{\bm}(E)$, with $s_\bm:= s_{m_1}\circ s_{m_2}\circ \ldots \circ s_{m_\ell}$.
We also set $I_0:=\{0\}$ and adopt the convention that $E_0:=E$. 
We also use the notation $s^0$ for the identity map on the set of subsets of $\R^n$, and the notations $s^1:=s$ and $s^{\ell+1}:= s\circ s^{\ell}$, for $\ell\in \N$.
We will use these notations especially in the case $E=\Gamma$. Indeed, for any fixed $\ell\in \N_0$, $\{\Gamma_\bm\}_{\bm \in I_\ell}$ provides a decomposition of $\Gamma$ into $M^\ell$ similar copies of $\Gamma$, whose pairwise overlaps have zero Lebesgue measure (cf.\ Proposition \ref{prop:nset}\eqref{ii}). 
For $\bm\in I_\ell$ the diameter of $\Gamma_\bm$ satisfies
\begin{align}
\label{eq:diameterIell}
\rho_{\rm min}^\ell h_0 \leq \diam(\Gamma_\bm) \leq \rho_{\rm max}^\ell h_0,
\end{align}
where $\rho_{\rm min}:=\min_{m=1,\ldots,M}\rho_m$, $\rho_{\rm max}:=\max_{m=1,\ldots,M}\rho_m$ and $h_0:=\diam(\Gamma)$. 
We also define, for any $0<h\leq h_0$, the set $L_h$ for $h=h_0$ by $L_h:=\{0\}$, and for $0<h<h_0$ by 
\begin{align}
\label{eq:LhDef}
L_{h} := \big\{\bm \in \bigcup_{\ell\in \N} I_\ell
: \diam(\Gamma_{\bm})\leq h \text{ and } \diam(\Gamma_{\bm_-})>h\big\},
\end{align}
where, for $\bm = (m_1,\ldots,m_\ell)$, $\bm_-:= (m_1,...,m_{\ell-1})$ if $\ell\geq 2$, and $\bm_-:= 0$ if $\ell=1$. 
For a given $0<h\leq h_0$, $\{\Gamma_\bm\}_{\bm \in L_h}$ provides a decomposition of $\Gamma$ into similar copies of $\Gamma$, each of which has diameter approximately equal to $h$. 
Explicitly, for $\bm\in L_h$ the diameter of $\Gamma_\bm$ satisfies
\begin{align}
\label{eq:diameterLh}
\rho_{\rm min}h < \diam(\Gamma_\bm) \leq h.
\end{align}
When the IFS is homogeneous, with $\rho_1=\ldots = \rho_M = \rho$ for some $\rho\in(0,1)$, the families of decompositions corresponding to the index sets $I_\ell$ and $L_h$ coincide, under the appropriate correspondence between $\ell$ and $h$ (specifically, for $\ell\in \N$ we have $I_\ell=L_h$ for $\rho^{\ell}h_0\leq h < \rho^{\ell-1}h_0$), but in general they differ. 




{Second, we recall the definitions of some notions of fractal dimension other than those already introduced above. 
One} 
fractal dimension that will be of importance to us is the so-called Aikawa dimension, $\dimA(F)$, of a non-empty set $F\subset \mathbb{R}^n$ (see \cite{Aikawa:91}, \cite[Definition 3.2]{LeTu:13}, and Appendix \ref{sec:Appendix4}, where we make clear that the definitions of \cite{Aikawa:91} and \cite[Definition 3.2]{LeTu:13} are equivalent). Following \cite{Aikawa:91} this is defined as the infimum of the set 
\begin{align}
\label{eq:IADef}
I_{\rm A}(F):=\{s>0:\sup_{x\in F}\sup_{r>0}r^{-s}\int_{B(x, r)} \operatorname{dist}(y, F)^{s-n} \, \rd y <\infty\},
\end{align}
with the convention that, if $\dist(y,F)=0$, then $\dist(y,F)^{s-n} := +\infty$ if $s<n$, $:=1$ if $s=n$, so that $\dimA(F):=\inf I_{\rm A}(F) =n$ if $\{y\in\ \R^n:\dist(y,F)=0\}=\overline{F}$ has positive measure. As discussed in  Appendix \ref{sec:Appendix4}, $s\in I_{\rm A}(F)$  if $s>\dimA(F)$.

It is shown in \cite[Theorem 1.1]{LeTu:13} that $\dimA(F)$ coincides with the better known Assouad dimension \cite{assouad1979etude} of $F$ (see Appendix \ref{sec:Appendix4} for a discussion of the equivalence of different definitions of this dimension; our primary definition is that of \cite{Fraser:21}, stated as \eqref{eq:IASDef} below, {but we note that, as is pointed out in \cite{luukkainen1998assouad}, a related concept was studied in much earlier work by Bouligand \cite{bouligand1928}.}). Given this result, $\dimA(F)$ should be read, below, as meaning either the Aikawa or the Assouad dimension of $F$.
We note (e.g., \cite[Lemma 2.4.3]{Fraser:21}) that $\dim_{\rm A}(E)\geq \dimH(E)$ for any set $E\subset\R^n$ and, if $E$ is bounded,
\begin{equation} \label{eq:uMin}
\dimH(E) \leq \underline{\dim}_{\rm M}(E)\leq \overline{\dim}_{\rm M}(E) \leq \dimA(E),
\end{equation}
where $\underline{\dim}_{\rm M}(E)$ and $\overline{\dim}_{\rm M}(E)$ denote, respectively, the  lower and upper Minkowski/box-counting dimensions of $E$ (e.g., \cite[\S2.1]{Fal}). (If $\underline{\dim}_{\rm M}(E)=\overline{\dim}_{\rm M}(E)$ the common value is $\dim_{\rm M}(E)$, the Minkowski/box-counting dimension of $E$.)
 
If $F$ is a compact $d$-set, in particular if $F$ is the attractor of an IFS of contracting similarities satisfying the OSC, then $\dimA(F)={\dim_{\rm M}(F)}=\dimH({F})=d$: 
the second equality is well known, and the first follows from \cite[Equation (3.1)]{Fraser:21} and \cite[Proposition 2.8]{Fraser:14}, noting that, as discussed in Appendix \ref{sec:Appendix4}, the definitions of the Assouad dimension of a non-empty $F\subset \R^n$ in \cite{Fraser:21} and \cite{Fraser:14} are subtly different, but coincide if $F$ is bounded. 
{That the equality between the three dimensions also holds when $F$ is the boundary of an $n$-attractor (noting that in this case $F$ is not necessarily a $d$-set --- see Remark \ref{rem:boundarydset}) is part \eqref{v} of the following proposition.} 

\begin{prop} \label{prop:nset}
Let $\Gamma\subset\R^n$ be an $n$-attractor. Then
\begin{enumerate}[(i)]
\item \label{i} $\Gamma$ is an $n$-set;
\item \label{ii} $\Gamma$ is \textit{self-similar} 
in the sense that the sets $\Gamma_m:=s_m(\Gamma)$, $m=1,\ldots,M$,
satisfy
\begin{align}
\label{eq:SelfSim}
|\Gamma_{m}\cap \Gamma_{m'}|=0, \qquad m\neq m';
\end{align}
\item \label{iii} $\Gamma = \overline{O}$, for any open $O$ satisfying the OSC. In particular, $O\subset\Gamma$ and hence $\Gamma^\circ\neq\emptyset$;
\item \label{iv} $|\partial O| =0$, indeed $\dimH(\partial O)<n$, for any open $O$ satisfying the OSC; hence also $|\partial \Gamma|=0$, indeed $\dimH(\partial \Gamma)<n$;
\item \label{v} $\dimA(\partial \Gamma)  {=\dim_{\rm M}(\partial \Gamma)} = \dimH(\partial \Gamma)$, so that $\dimA(\partial \Gamma) < n$;
\item \label{vi} $\Gamma = \overline{\Gamma^\circ}$, and hence $\partial\Gamma = \partial(\Gamma^\circ)=\partial(\Gamma^c)$ {and $\dimH(\partial\Gamma)\geq n-1$};
\item \label{vii} $\Gamma^\circ$ is an open $n$-set;
\item \label{viii} the OSC holds with $O=\Gamma^\circ$.
\item \label{ix} $\Gamma^c$ is an open $n$-set.
\item \label{x} $\partial\Gamma$ is a subset of a compact $d$-set, for some $n-1\leq d<n$. 
\item \label{xi} $\partial\Gamma$ is porous (satisfies the ball condition) in the sense of Definition \ref{def:Porous} (and see \cite[Definition~3.4]{Tri08}). 
\end{enumerate}
\end{prop}


\begin{proof}
Parts \eqref{i} and \eqref{ii} are a special case of \cite[Theorem~4.7]{Triebel97FracSpec}. Part \eqref{iii} is shown in the proof of \cite[Corollary 2.3]{Sc:94} (or see \cite[Lemma 2.6]{HausdorffBEM}).

If $O$ is an open set satisfying the OSC, then $O\cap \Gamma = O$ by part \eqref{iii}, so that $O\cap \Gamma\neq \emptyset$, i.e., $O$ satisfies the so-called {\em strong} OSC (e.g., \cite{Keesling:99}). Further $\overline{O}=\Gamma$ by \eqref{iii}, so that $\partial \Gamma\subset \partial O= \Gamma\setminus O$, and part \eqref{iv} follows from \cite[Corollary 3.2]{Keesling:99}, which gives that $\dimH(\Gamma\setminus O)<\dimH(\Gamma)=n$.

{For part \eqref{v}, the equality $\dimA(\partial \Gamma) = \dimH(\partial \Gamma)$  holds by \cite[Corollary 2.11]{Fraser:14} and the equality $\dimH(\partial \Gamma)=\dim_{\rm M}(\partial \Gamma)$ holds by \cite[2.(f) and Theorem 3.5]{Falconer:95}. The fact that $\dimA(\partial \Gamma)<n$ then follows from part \eqref{iv}.}

Parts \eqref{vi} and \eqref{vii} follow from parts \eqref{i} and \eqref{ii}, 
using the relationship between $n$-sets and open $n$-sets mentioned after \eqref{eq:nset}. 
{The statement in \eqref{vi} that $\dimH(\partial\Gamma)\geq n-1$ holds because $\partial\Gamma=\partial(\Gamma^\circ)$ and $\dimH(\partial(\Gamma^\circ))\geq n-1$. The latter holds because $\dimH(\partial\Omega)\geq n-1$ for any domain $\Omega\subset\R^n$ such that $\overline\Omega\neq \R^n$, which follows from the fact that projections from $\R^n$ to $\R^{n-1}$ do not increase Hausdorff dimension, cf.\ \cite[Proof of Proposition~3.18(iii)]{Tri08}.}

Part \eqref{viii} is proved as part of the proof of \cite[Claim 1]{Keesling:99}. 
For the convenience of the reader we include the short argument.
Note that $s(\Gamma^\circ)\subset s(\Gamma)=\Gamma$, so that, since $s(\Gamma^\circ)$ is open, $s(\Gamma^\circ)\subset \Gamma^\circ$. Moreover, $s_m(\Gamma^\circ)$ is open, for $m=1,\ldots,M$, so that $s_m(\Gamma^\circ)\cap s_{m'}(\Gamma^\circ)=\emptyset$, for $m\neq m'$, for otherwise $|s_m(\Gamma)\cap s_{m'}(\Gamma)|\geq |s_m(\Gamma^\circ)\cap s_{m'}(\Gamma^\circ)| >0$, contradicting \eqref{eq:SelfSim}. Thus the OSC holds with $O=\Gamma^\circ$.

For part \eqref{ix}, to see that $\Gamma^c$ is an open $n$-set we will adapt well-known arguments which show that $\R^n$ can be tiled (in the sense of \cite[Definition.~3.1]{BaVi:14}) by similar copies of $\Gamma$ (see, e.g., \cite[Theorem 1.2]{Keesling:99}, \cite[Theorem 3.8]{BaVi:14}). To show that $\Gamma^c$ satisfies \eqref{eq:nset} (with $F$ replaced by $\Gamma^c$) it is clearly enough to show that \eqref{eq:nset} holds (with $F$ replaced by $\Gamma^c$) for all $x\in U^\prime\cap \Gamma^c$, where $U^\prime$ is some open neighbourhood of $\partial \Gamma$, and this follows if $\Omega\subset \Gamma^c$ is an open $n$-set and $U\cap \Gamma^c= U\cap \Omega$, with $U$ some open neighbourhood of $\Gamma$.

We have shown above that the OSC holds with $O:=\Gamma^\circ$. Making this choice for $O$, recall also that $\overline{O}=\Gamma$ and that, for all $\ell\in \N$, $\Gamma = s^\ell(\Gamma)= \cup_{\bm\in I_\ell}s_\bm(\Gamma)$,
so that every point in $O$ is contained in $s_\bm(\Gamma)$ for some $\bm\in I_\ell$. Noting also that $\max_{\bm\in I_\ell}\diam(s_\bm(\Gamma))\to 0$ as $\ell\to\infty$, we see that we can
choose $\ell \in \N$ large enough so that $s_\bm(\Gamma)\subset O$ for some $\bm\in I_\ell$. For
any such choice of $\ell$ and $\bm$, $s_\bm^{-1}$ is a similarity that maps $s_\bm(\Gamma)$ to $\Gamma$ and maps $O$ to the open set $U:=s_\bm^{-1}(O)$ which, since $O\supset s_\bm(\Gamma)$, is an open neighbourhood of $\Gamma$. Further, 
for $\bm'\in I_\ell$, $\Omega^\bm_{\bm'} := (s_{\bm}^{-1}(s_{\bm'}(\Gamma)))^\circ= s_{\bm}^{-1}(s_{\bm'}(O))$ is similar to $O$ and so is an open $n$-set (indeed $\Omega^\bm_\bm=O$), so that
\begin{equation} \label{eq:Omega*}
\Omega_*:=\cup_{\bm'\in I_\ell, \bm'\neq \bm} \Omega^\bm_{\bm'}
\end{equation}
is an open $n$-set. Moreover, by the OSC, for $\bm'\neq \bm$, $\Omega^\bm_{\bm'}\cap \Omega^\bm_\bm=\emptyset$, so that $\overline{\Omega_*}\cap O = \overline{\Omega_*}\cap \Omega^\bm_\bm = \emptyset$, but also $\overline{\Omega_*}\cup \Gamma = \cup_{\bm'\in I_\ell} \overline{\Omega^\bm_{\bm'}}=s_\bm^{-1}(\cup_{\bm'\in I_\ell}s_{\bm'}(\Gamma))=s_\bm^{-1}(\Gamma)$. Thus
$$
\tilde{\Omega}_* := \left(\overline{\Omega_*}\right)^\circ = \left(s_\bm^{-1}(\Gamma)\setminus O\right)^\circ = U\setminus \Gamma \subset \Gamma^c.
$$
Since $\tilde{\Omega}_*$ is an open $n$-set (using the relationship between $n$-sets and open $n$-sets mentioned after \eqref{eq:nset}, and since $\Omega_*$ is an open $n$-set),  $U$ is an open neighbourhood of $\Gamma$, and $U\cap \Gamma^c= U\cap \tilde{\Omega}_*$, we conclude that $\Gamma^c$ is an open $n$-set.

For part \eqref{x}, to prove that, 
for some $d<n$, $\partial \Gamma\subset \Gamma^\prime$, where $\Gamma^\prime$ is a compact $d$-set, we adapt the argument of \cite{Keesling:99} (which builds on \cite{Falconer:95}), where a similar approach was used to prove the weaker statement in \eqref{iv} that $\dimH(\partial \Gamma)<n$. Note first that $\Gamma$, the attractor of $\{s_1,\ldots, s_M\}$, satisfying $\Gamma=s(\Gamma)$, is also  the attractor of the IFS $S_\ell:=\{s_\bm:\bm\in I_\ell\}$, i.e.\ $\Gamma=s^\ell(\Gamma)$, for every $\ell\in \N$. Now consider the IFS $S_\ell^\prime:=\{s_\bm:\bm\in I^\prime_\ell\}$, where $I_\ell^\prime:=\{\bm\in I_\ell: s_\bm(\Gamma)\not\subset \Gamma^\circ\}$. For all $\ell$ sufficiently large (as argued in our proof of \eqref{ix}), $I^\prime_\ell\subsetneqq I_\ell$, so that $S_\ell^\prime$ is a strict subset of $S_\ell$. Further, let $O^\prime:=s^\prime(\Gamma^\circ)$, where $s^\prime(E):= \cup_{\bm\in I_\ell^\prime} s_\bm(E)$, for $E\subset \R^n$. Since our original IFS satisfies the OSC with $O=\Gamma^\circ$, we have that  $s_\bm(\Gamma^\circ)\subset \Gamma^\circ$ and $s_\bm(\Gamma^\circ)\cap s_{\bm'}(\Gamma^\circ)=\emptyset$, for all $\bm, \bm'\in I_\ell$, so that $O^\prime\subset \Gamma^\circ$. 
{Hence} 
$s_\bm(O^\prime)\subset s_\bm(\Gamma^\circ)\subset O^\prime$, for $\bm\in I_\ell^\prime$, so that $S_\ell^\prime$ satisfies the OSC with $O=O^\prime$. 
{Hence} its attractor $\Gamma^\prime$ is a compact $d$-set for some $d\in (0,n]$; further, if we choose $\ell$ sufficiently large, $S_\ell^\prime$ is a strict subset of $S_\ell$ whose attractor is the $n$-set $\Gamma$, so that $d<n$ (compare the equations \eqref{eq:dfirst} defining the similarity dimensions of $\Gamma$ and $\Gamma^\prime$). 
Further, $\Gamma=s(\Gamma)$ and $s(\Gamma^\circ)\subset \Gamma^\circ$, so that $\partial \Gamma \subset s(\partial \Gamma)$ and $\partial \Gamma \subset s^\ell(\partial \Gamma)$, for $\ell\in \N$. 
But this implies, by the definition of $I_\ell^\prime$, that $\partial \Gamma \subset s^\prime(\partial \Gamma)$ (because if $\bm\in I_\ell\setminus {I_\ell'}$ then $s_\bm(\partial\Gamma)\cap \partial\Gamma=\emptyset$), which in turn implies that $\partial \Gamma \subset \Gamma^\prime$. 
{The fact that $d\geq n-1$ follows from the fact that $\dimH(\partial\Gamma)\geq n-1$, as established in part \eqref{vi}.}

Finally, part \eqref{xi} follows from the fact that the $d$-set $\Gamma^\prime$ constructed above is porous, by \cite[Proposition~4.3]{Cae02}, 
and the fact that $\Gamma\subset \Gamma^\prime$. 
\end{proof}

\begin{rem}
\label{rem:boundarydset}
If $\Gamma\subset\R^n$ is an $n$-attractor it is not guaranteed that $\partial\Gamma$ is a $d$-set. A counterexample in the case $n=1$ is provided by the set $\Gamma$ of \eqref{eq:Gammaeq}, 
for which 
$\partial\Gamma$ 
is a countable set with an accumulation point at $0$, and hence $\partial\Gamma$ is neither a $0$-set (because the $0$-set condition fails for $x=0$) nor a $d$-set for any $0<d\leq 1$ (because $\partial\Gamma$ is countable and hence $\dimH(\partial\Gamma)=0$). 
A counterexample in the case $n=2$ is $\Gamma:=\Gamma_{1/2}\times[0,1]$, 
where $\Gamma_{1/2}$ is the set defined by \eqref{eq:Gammaeq} with $\alpha=1/2$, as shown in Figure \ref{fig:2sets}(d). This has $\dimH(\partial\Gamma)=1$ but $\partial\Gamma$ is not a $1$-set, because the $1$-set condition fails on $\Gamma_{1/2}\times\{0\}$. Similarly, $\Gamma:=\Gamma_{1/2}\times[0,1]^{n-1}$ is a counterexample for $n\geq 3$. Generalising this observation, it is possible to show that the boundary $\partial \Gamma$ of this set is also not an $h$-set as defined in Definition \ref{def:hset}, for any choice of $h:[0,1]\to[0,1]$ satisfying the conditions of the definition. The relationship of this observation  to our results in \S\ref{sec:Multipliers}, in particular to Corollary \ref{cor:pointmult}, is explained in Remark \ref{rem:multipliers}. 
\end{rem}

\begin{rem}
In relation to part \eqref{x} of Proposition \ref{prop:nset}, it is natural to ask whether every compact set $K\subset\R^n$ with $\dimH(K)<n$ is contained in a $d$-set for some $0\leq d<n$. This is not true. 
One way to produce a counterexample is to consider a compact set $K$ with empty interior that is not porous (i.e.\ does not satisfy the ball condition) in the sense of Definition \ref{def:Porous}, since by \cite[Proposition~4.3 and Remark~3.5(i)]{Cae02} any $d$-set with $0<d<n$ is porous, and it is immediate from Definition \ref{def:Porous} that any subset of a porous set is porous. 
An example in $\R$ is the set $K:=\{0\}\cup\{1/m:m\in\N\}$, which was shown not to be porous in \cite[\S2]{Cae02}. 
Note that while in this case $\dimH(K)=0$, $K$ cannot be contained in a $0$-set, because that would contradict the fact that ${\cal H}^{0}$ is just the counting measure. 
An example in $\R^n$ for $n\geq 2$ is given by $K\times [0,1]^{n-1}$, where $K\subset \R$ is defined as above. 
A further example in $\R^{2}$ is the anisotropic fractal considered in \cite[\S5, Example 6]{Cae02}. 
\end{rem}

\subsection{Thickness and the class $\cD^t$}
The following Theorem \ref{thm:thick} is the main result of this section.
The thickness of $O:=\Gamma^\circ$ established in this theorem has significance for proving a range of results involving function spaces on $O$, including results on the continuous extension to $\R^n$ of functions defined on $O$, and associated results on interpolation of function spaces, see \cite[Chapter 4]{Tri08}.
One key application we will make of Theorem \ref{thm:thick} will be to deduce the equality of the Sobolev spaces $\widetilde{H}^{s}(\Gamma^\circ)$ and $H_{\Gamma}^{s}$ (defined in \S\ref{sec:FunctionSpaces}), for $s\in \R$, 
which we do in Corollary \ref{cor:tildesubscriptHs}.

\begin{thm}
\label{thm:thick}
Let $\Gamma\subset\R^n$ be an $n$-attractor. 
Then both
$\Gamma^\circ$ and $\Gamma^c$ are thick in the sense of Triebel (see Definition \ref{def:Thick}, \cite[Definition~3.1(ii)-(iv), Remark~3.2]{Tri08} or \cite[Definition~4.5]{caetano2019density}).
\end{thm}

\begin{proof}
Setting $O=\Gamma^\circ$, let us show that $O$ is thick, which requires showing that $O$ is both $I$-thick and $E$-thick, in the sense of \cite[Definition~4.5]{caetano2019density} (see Definition \ref{def:Thick}). 

Our argument to show that $O$ is $I$-thick will rely on a decomposition of $\Gamma$ into similar copies of $\Gamma$ of approximately equal diameter. Following the notation introduced before Proposition \ref{prop:nset}, for every $\ell\in \N$ and every $\bm\in I_\ell$, let $O_\bm := s_\bm(O)$ and $\Gamma_\bm:= s_\bm(\Gamma)=\overline{O_\bm}$. 
Let $\nu_0$ be the unique integer such that
$2^{-\nu_0}\leq \diam(\Gamma)<2^{-\nu_0+1}$.
For $\nu\geq \nu_0$ let (cf.\ the notation $L_h$ in \eqref{eq:LhDef}) 
\begin{equation} \label{eq:Knudef}
K_{\nu}  := \big\{\bm \in
\cup_{\ell'\in\N}I_{\ell'}
: \diam(\Gamma_{\bm})< 2^{-\nu} \text{ and } \diam(\Gamma_{\bm_-})\geq 2^{-\nu}\big\}.
\end{equation}
Then, for $\nu\geq \nu_0$, by repeated application of \eqref{eq:fixedfirst} and \eqref{eq:SelfSim},
\begin{equation} \label{eq:KnuGamma}
\Gamma = \bigcup_{\bm \in K_\nu} \Gamma_\bm \quad \mbox{and} \quad |\Gamma_\bm\cap \Gamma_{\bm'}|=0, \quad \bm, \bm'\in K_\nu, \;\; \bm\neq \bm'.
\end{equation}
As a further preliminary, choose $r>0$ and $x\in O$ such that $B(x,r)\subset O$, and let $Q\subset B(x,r)$ denote an open cube contained in $B(x,r)$ such that the corner points of $\overline{Q}$ lie on $\partial B(x,r)$. For every $\ell\in \N$ and every $\bm=(m_1,\ldots,m_\ell)\in I_\ell$, 
  let $\tilde{Q}_\bm\subset O_\bm$ denote the open cube contained in $s_\bm(B(x,r))$ which has sides parallel to the coordinate axes and is such that $\tilde{Q}_\bm$ is congruent to $Q_\bm=s_\bm(Q)$. Then, where $l(Q^\prime)$ denotes the length of one of the sides of a cube $Q^\prime$, for some constants $C_1$, $C_2$, $C_3$, $C_4>0$, independent of $\bm=(m_1,\ldots,m_\ell)$,
\begin{equation}
C_1 \diam(\Gamma_\bm) \leq \ell(\tilde{Q}_\bm) \leq C_2 \diam(\Gamma_\bm) \; \mbox{and} \; C_3 \diam(\Gamma_\bm) \leq \dist(Q_\bm,\partial \Gamma_\bm) \leq  C_4 \diam(\Gamma_\bm).
\end{equation}

Now, to prove that $O$ is $I$-thick, suppose that $0< c_1\leq c_2$, $0<c_3\leq c_4$, and $j_0\in \N$. Suppose that $Q^e\subset O^c$ is an open cube such that
$$
c_1 2^{-j} \leq l(Q^e) \leq c_2 2^{-j} \quad \mbox{and} \quad c_3 2^{-j} \leq\dist(Q^e,\partial O) \leq c_4 2^{-j},
$$
for some $j\in \N$ with $j\geq j_0$. Choose $x^*\in \partial O = \partial \Gamma$ such that $D^e:= \dist(Q^e,\partial O)=\dist(Q^e,x^*)$ and let 
$J:=\max(1,\nu_0-1)\in \N$.
Then, where $\nu:= j+J$ and $K_\nu$ is defined by \eqref{eq:Knudef}, it holds that $\nu\geq \nu_0$ and, by the definition of $K_\nu$, that
$$
\min_m\rho_m 2^{-\nu}\leq \diam(\Gamma_\bm)<2^{-\nu}, \quad \bm\in K_\nu.
$$
Now choose $\bm\in K_\nu$ such that $x^*\in \Gamma_\bm$, which is possible by \eqref{eq:KnuGamma}. Then
$$
\dist(\tilde{Q}_\bm,Q^e)\leq D^e + \dist(\tilde{Q}_\bm,x^*)\leq D^e+\diam(\Gamma_\bm) < (c_4+2^{-J}) 2^{-j},
$$
$$
\dist(\tilde{Q}_\bm,\partial O) \geq \dist(\tilde{Q}_\bm,\partial \Gamma_\bm)\geq C_3 \diam(\Gamma_\bm) \geq C_32^{-J}\min_m\rho_m 2^{-j},
$$
and
$$
C_12^{-J}\min_m\rho_m2^{-j}\leq l(\tilde{Q}_\bm) < C_22^{-J}2^{-j}.
$$
Since also $\dist(\tilde{Q}_\bm,\partial O)\leq \dist(Q^e,\tilde{Q}_\bm)$ (since $Q_\bm\subset O$ and $Q^e\subset O^c$), we see, inspecting \cite[Definition~4.5]{caetano2019density}, that we have shown that $O$ is $I$-thick.

To see that $O$ is $E$-thick it is enough, by \cite[Proposition~3.6]{Tri08} (see also \cite[Theorem~4.7]{caetano2019density}), and since $O=\Gamma^\circ=(\overline O)^\circ$, to show that $\Gamma^c$ is $I$-thick. But this follows from the following observations, in which, to avoid having to repeatedly consider the properness of the subsets involved, we trivially extend the definition of $I$-thick sets
to include $\R^n$ itself.
Firstly, noting \cite[Remark~4.6]{caetano2019density}, it is clear that a domain $\Omega\subset \R^n$ is $I$-thick if and only if $U^\prime\cap \Omega$ is $I$-thick, for some open set $U^\prime$ such that $\partial \Omega\subset U^\prime$. Regarding the not so obvious implication here, we note that choosing $j_0$ in the definition of $I$-thick large enough we get for the $Q^e \subset \Omega^c$ we start with that $\dist(Q^e,\partial\Omega)=\dist(Q^e,\partial(U'\cap\Omega))$. Similarly, for the $Q^i \subset U'\cap\Omega$ we arrive at, we have that $\dist(Q^i,\partial\Omega)=\dist(Q^i,\partial(U'\cap\Omega))$. 
Secondly, if a domain $\Omega\subset \R^n$ is $I$-thick, then also $(\overline{\Omega})^\circ$ is $I$-thick,
since, for every open cube
$Q^e\subset ((\overline{\Omega})^\circ)^c$,
$\dist(Q^e,\partial \Omega)$ = $\dist(Q^e,\partial((\overline{\Omega})^\circ))$.
Thirdly, if $\Omega_1$, $\Omega_2\subset \R^n$ are $I$-thick domains, then $\Omega:=\Omega_1\cup \Omega_2$ is $I$-thick;
this is clear from the $I$-thick definition since, for every open cube
$Q^e\subset \Omega^c$,
$\dist(Q^e,\partial \Omega)=\dist(Q^e,\partial \Omega_j)$, for $j=1$ or $2$,
and that for every open cube $Q^i\subset \Omega_j$, $j=1$ or $2$, $\dist(Q^{i},\partial\Omega_j)\leq\dist(Q^i,\partial\Omega)$. 
To prove the first claim here, suppose that $Q^e\subset \Omega^c$ and $x\in \partial \Omega$ is such that $\dist(Q^e,\partial \Omega)=\dist(Q^e,x)$. Since $\partial \Omega\subset \partial \Omega_1\cup \partial \Omega_2$ it follows that $x\in \partial \Omega_j$, for $j=1$ or $2$; let us suppose without loss of generality that $x\in \partial \Omega_1$. Either $\dist(Q^e,x)=\dist(Q^e,\partial \Omega_1)$, in which case we are done, or $\dist(Q^e,x)>\dist(Q^e,\partial \Omega_1)$, in which case there exists $y\in \partial \Omega_1$ with $\dist(Q^e,\partial \Omega)=\dist(Q^e,x)>\dist(Q^e,y)$. But this implies that $y\not\in \partial \Omega$, so that $y\in \Omega$, which implies in turn that there exists $z\in \partial \Omega$ with $\dist(Q^e,z)<\dist(Q^e,y)$, a contradiction.
To prove the second claim, let us suppose without loss of generality that $Q^i\subset\Omega_1$. Let $x\in\partial\Omega$ be such that $\dist(Q^{i},x)=\dist(Q^{i},\partial\Omega)$.
If $x\in\partial\Omega_1$, then $\dist(Q^{i},\partial\Omega_1)\leq\dist(Q^{i},x)=\dist(Q^{i},\partial\Omega)$
and we are done. If $x\notin\partial\Omega_1$, then $x\not\in\Omega_1$ too. To see this, with $x\not\in\partial\Omega_1$ we have that there
exists $r>0$ such that $B(x,r)\cap\Omega_1=\emptyset$ or $B(x,r)\cap\Omega_1^{c}=\emptyset$;
however, from $x\in\partial\Omega$ it follows that $B(x,r)\cap\Omega_1^{c}\cap\Omega_{2}^{c}\not=\emptyset$,
therefore it must be $B(x,r)\cap\Omega_1=\emptyset$, hence $x\not\in\Omega_1$. Hence $\dist(Q^{i},x)\geq\dist(Q^{i},\partial\Omega_1)$,
therefore $\dist(Q^{i},\partial\Omega_1)\leq\dist(Q^{i},\partial\Omega)$
and we are also done.

To make use of these observations we argue as we did in the proof of Proposition \ref{prop:nset}\eqref{ix} to show that $\Gamma^c$ is an open $n$-set. Choose $\ell\in \N$ and $\bm\in I_\ell$ so that $s_\bm(\Gamma)\subset O$, and, as above, define $U:=s_\bm^{-1}(O)$ and define $\Omega_*$ by \eqref{eq:Omega*}. Then, as we have noted above, $U$ is an open neighbourhood of $\Gamma$, and, defining $\tilde{\Omega}_*:= \left(\overline{\Omega_*}\right)^\circ$, $\tilde{\Omega}_* = U\cap \Gamma^c$. But $\Omega_*$ is a finite union of similar copies of $O$. Thus, using the above observations, since $O$ is $I$-thick, $\Omega_*$ is $I$-thick. By the above observations this implies that $\tilde{\Omega}_*$ is $I$-thick, which implies in turn that $\Gamma^c$ is $I$-thick.

Having established the thickness of $\Gamma^\circ$, the thickness of $\Gamma^c$ follows by \cite[Proposition~3.6(iii)-(iv)]{Tri08}.
\end{proof}

The following corollary, which we include 
because it may be of independent interest, follows from Theorem \ref{thm:thick}, Proposition \ref{prop:nset}\eqref{xi} and Proposition \ref{prop:ThickPorous}.

\begin{cor}
\label{cor:Eporous}
Let $\Gamma\subset\R^n$ be an $n$-attractor. 
Then both
$\Gamma^\circ$ and $\Gamma^c$ are $E$-porous in the sense of Triebel (see Definition \ref{def:EPorous} and \cite[Definition~3.16(i)]{Tri08}).
\end{cor}

Next we prove, in 
{Theorem} 
\ref{lem:newDt}, that, for suitable $0<t<1$, 
both the interior $\Gamma^\circ$ and complement $\Gamma^c$ of an $n$-attractor $\Gamma$ belong to the class $\cD^t$ defined in \cite[Definition 4.2]{Sickel99a} as a variant on the Frazier-Jawerth class $\cD_t$ of \cite{FrazierJawerth90}. This class $\cD^t$ is defined, for $0<t<1$,
 to be
 the class of domains $\Omega\subsetneqq\R^n$ such that 
\begin{equation}
\sup_{x\in\partial\Omega}\sup_{0<r\leq1}r^{t-n}\int_{B(x,r)\setminus\partial\Omega}{\rm dist}(y,\partial\Omega)^{-t}\,\rd y<\infty.\label{eq:D_t}
\end{equation}
As will be noted in Proposition \ref{prop:pointmult} below, membership of $\cD^t$ implies \cite{FrazierJawerth90,Sickel99a} that the characteristic function of $\Omega$ is a pointwise multiplier on a certain range of function spaces. We will use this multiplier property in \S\ref{sec:FunctionSpaces} to construct extension operators and prove interpolation properties. 
{Theorem} \ref{lem:newDt}  
 connects the class $\cD^t$ with the Aikawa/Assouad dimension. 
In connection with part (i) of the 
{theorem}, 
a sufficient condition for $n-{\rm dim_{A}}(\partial\Omega)\leq 1$ is that $\overline{\Omega}\not=\R^n$, because in that case 
${\rm dim_{H}}(\partial\Omega)\geq n-1$, which 
follows by the argument used to prove the last part of Proposition \ref{prop:nset}\eqref{vi}, 
and then the fact that 
${\rm dim_{A}}(\partial\Omega)\geq {\rm dim_{H}}(\partial\Omega)$ (see \cite[Lemma 2.4.3]{Fraser:21}) 
gives that also 
${\rm dim_{A}}(\partial\Omega)\geq n-1$.
Note that $|\partial \Omega|=0$ if $\Omega\in \cD^t$ for some $t\in (0,1)$, by the arguments of \cite[Remark 4.3]{Sickel99a}.






\begin{thm}
\label{lem:newDt}
Let $\Omega\subsetneqq\R^n$ be a domain. 
\begin{itemize}
\item[(i)] If 
${\rm dim_{A}}(\partial\Omega)<n$ then 
$\Omega\in \cD^t$ 
for all $0< t<\min(n-{\rm dim_{A}}(\partial \Omega),1)$.
\item[(ii)] If $\partial\Omega$ is bounded 
and 
$\Omega\in \cD^t$, for some $0<t<1$, then $\dim_{\rm A}(\partial\Omega)\leq n-t$.
\end{itemize}



\end{thm}

\begin{proof}
(i)	Observe first that the range 
$0 < t<\min(n-{\dimA}{(\partial \Omega)},1)$ 
is not void (by the assumption that ${\dimA}(\partial\Omega)<n$).
Now take $s:=n-t$, so that $s>\dimA{(\partial\Omega)}$. Then, by the comment after \eqref{eq:IADef}, $s\in I_{\rm A}(F)$ so there exists 
a constant $ c_{s} > 0 $ such that	
		\[
		\int_{B(x, r)} \operatorname{dist}(y, \partial\Omega)^{s-n} \, \rd y \leq c_{s} \, r^{s}
		\]
		\noindent for every \( x \in \partial\Omega \) and all \( r > 0 \). From this it follows
that \eqref{eq:D_t} holds, i.e.\ that
$\Omega \in \cD^t$. 

(ii)
If $\partial\Omega$ is bounded,  $0<t<1$, and $\Omega\in \cD^t$, in which case $|\partial\Omega|=0$, then, defining $s:=n-t>0$, $s\in I_{\rm A}^{r_0}(\partial \Omega)$, for $0<r_0<1$, where $I_{\rm A}^{r_0}(\partial \Omega)$ is defined by \eqref{eq:IADefr2}. It follows that $s\in I_{\rm A}(\partial \Omega)$ by Lemma \ref{lem:equiv}, so that $\dimA(\partial \Omega)\leq s = n-t$.
\end{proof}

\begin{cor}
\label{cor:newDt}
Suppose that $\Gamma\subset\R^n$ is an $n$-attractor and $0<t<1$, and let $\Omega$ denote either $\Gamma^\circ$ or $\Gamma^c$.
Then $\Omega$ belongs to the
class $\cD^t$  
if $t<n-\dimH(\partial \Omega)$, while $\Omega$ does not belong to $\cD^t$ if $t>n-\dimH(\partial \Omega)$.
\end{cor}
\begin{proof}
By Definition \ref{def:1}, $\Gamma$ is bounded and, by Proposition \ref{prop:nset}\eqref{v},\eqref{vi}, we have $\partial\Omega=\partial\Gamma$ and $n-1\leq \dimH(\partial\Gamma)=\dim_{\rm A}(\partial\Omega)<n$. Thus the result follows by 
{Theorem} 
\ref{lem:newDt}.
\end{proof}

\begin{rem}
\label{rem:direct}
The proofs of two of our main results about $n$-attractors, namely Corollaries \ref{cor:tildesubscriptHs} and \ref{cor:Interpolationnset} below, do not need a result as strong as that provided by 
Corollary \ref{cor:newDt}.  
Instead they merely need that $\Omega$ belongs to ${\cal D}^t$ for some small positive $t$, which we can prove more directly using the self-similarity properties of $\Gamma$, without resorting to the comparison made between different notions of dimension as we did in the proof of 
Corollary \ref{cor:newDt}. 
More precisely, using such direct arguments
we can prove that $\Omega$ belongs to ${\cal D}^t$ for all $0 < t < n-d$, for the $d$ of Proposition \ref{prop:nset}\eqref{x} (which may be strictly larger than $\dimH(\partial\Gamma)$ but is guaranteed to be smaller than $n$). For the interested reader, details can be found in Appendix \ref{sec:Appendix3}.
\end{rem}

\section{Function spaces on $n$-sets, $n$-attractors, and related rough domains}
\label{sec:FunctionSpaces}

In this section we present results on properties of Sobolev, Triebel-Lizorkin and Besov spaces on $n$-sets and $n$-attractors and their interiors, and on related classes of rough domains. In \S\ref{sec:FunctionSpaceNotation} we introduce our notation and basic definitions. In \S\ref{sec:Multipliers} we prove results about pointwise multiplication by characteristic functions, which we use to define extension operators and hence prove density results. In \S\ref{sec:Interpolation} we prove results about function space interpolation. 
Finally, in \S\ref{sec:Approximation} we 
derive best approximation error estimates in fractional Sobolev spaces for approximation by piecewise constant functions on certain fractal meshes of $n$-attractors. 


\subsection{Definitions and notation}
\label{sec:FunctionSpaceNotation}

{For a domain $\Omega\subset \R^n$, let ${\cal D}(\Omega)=C_0^\infty(\Omega)$ and ${\cal D}'(\Omega)$ denote the usual spaces of test functions and distributions on $\Omega$. Let ${\cal S}(\R^n)\subset C^\infty(\R^n)$ denote the usual Schwartz space of rapidly decreasing functions,  and ${\cal S}'\Rn$ the space of tempered distributions.} 
For $s\in\R$, $0<p<\infty$ and $0<q\leq \infty$ let $A^s_{p,q}(\R^n)$ denote either the Besov space $B^s_{p,q}(\R^n)$ or the Triebel-Lizorkin space $F^s_{p,q}(\R^n)$ (for definitions see e.g.\ \cite[Definition~1.1]{Tri08}, and for some basic properties see e.g.\ \cite[\S2.3.2, \S2.3.3]{Triebel83ThFS}). We note that $F^s_{2,2}(\R^n)=B^s_{2,2}(\R^n)=H^s(\R^n)$, with equivalent norms, where $H^s\Rn$ denotes the standard 
{(fractional)} 
Sobolev space introduced in \S\ref{sec:Intro}, a Hilbert space. Indeed, $F_{p,2}(\R^n)=H^s_p(\R^n)$, for $s\in \R$ and $1<p<\infty$, where $H^s_p(\R^n)$ denotes the standard $L_p$-version {(fractional)} Sobolev space/Bessel potential space of order $s$, as defined, e.g., in \cite[p.~37]{Triebel83ThFS}. Further, $F^s_{p,p}(\R^n)=B^s_{p,p}(\R^n)$, with equivalent {quasi-}norms, for $0<p<\infty$ and $s\in \R$.

For a domain $\Omega\subset \R^n$ we
define
$\tA^s_{p,q}(\Omega):=\overline{C_0^\infty(\Omega)}^{A^s_{p,q}(\R^n)}$.
We denote by $A^s_{p,q}(\Omega)$ the space of restrictions to $\Omega$ of elements of $A^s_{p,q}(\R^n)$,
equipped with the quotient {quasi-}norm
\[\|u\|_{A^s_{p,q}(\Omega)} := \inf_{\substack{\varphi\in A^s_{p,q}\Rn\\ \varphi|_\Omega=u}}\|\varphi\|_{A^s_{p,q}\Rn}.\]
For a closed set $E\subset \R^n$ we denote by $A^s_{p,q,E}$ the set of all elements of $A^s_{p,q}\Rn$ whose distributional support is contained in $E$. 
Analogously, we define $\tH^s_p(\Omega)$, $H^s_p(\Omega)$, and $H^s_{p,E}$, for $s\in \R$ and $1<p<\infty$, so that $\tH^s_p(\Omega)=\tF^s_{p,2}(\Omega)$, $H^s_p(\Omega)=F^s_{p,2}(\Omega)$, and $H^s_{p,E}=F^s_{p,2,E}$.

Now suppose that $s\in \R$ and $1<p,q<\infty$. 
For a domain $\Omega\subset\R^n$ the space $A^s_{p,q}(\Omega)$ is isometrically isomorphic to the quotient space $A^s_{p,q}(\R^n)/A^s_{p,q,\Omega^c}$, since $A^s_{p,q,\Omega^c}$ is the kernel of the restriction operator $|_\Omega:A^s_{p,q}(\R^n)\to A^s_{p,q}(\Omega)$. 
Also, it is well known that 
$A^s_{p,q}\Rn$ is reflexive\footnote{Indeed, that the canonical embedding into the second dual is surjective is immediate from our discussion below of the duality between  $A^s_{p,q}(\R^n)$ and $A^{-s}_{p',q'}(\R^n)$.}, and hence so are  $\tA^s_{p,q}(\Omega)$, $A^s_{p,q,E}$ and $A^s_{p,q}(\Omega)$, since they are all either a closed subspace of $A^s_{p,q}\Rn$ or a quotient of $A^s_{p,q}\Rn$ by a closed subspace of $A^s_{p,q}\Rn$.
Furthermore, with $p'$ and $q'$ denoting the usual conjugate exponents of $p$ and $q$, $A^s_{p,q}(\R^n)$ is a realisation of the dual space of $A^{-s}_{p',q'}(\R^n)$, with the duality pairing $\langle \cdot,\cdot\rangle_{A^s_{p,q}(\R^n)\times A^{-s}_{p',q'}(\R^n)}$
extending the {restriction of the $L_2(\R^n)$ inner product to Schwartz
functions}\footnote{As in \cite{ChaHewMoi:13}, for convenience we adopt the convention that distributions are antilinear functionals.}.
With respect to this duality the spaces $\tA^{-s}_{p',q'}(\Omega)$ and $A^s_{p,q,\Omega^c}$ are mutual annihilators \cite[Proposition\ 3.5]{caetano2019density} (and see \cite[Lemma~3.2]{ChaHewMoi:13} for the $p=2$ Sobolev space case). Hence, recalling the observation above about $A^s_{p,q}(\Omega)$ being isometrically isomorphic to $A^s_{p,q}(\R^n)/A^s_{p,q,\Omega^c}$, 
it follows by \cite[Theorem~1.10.16]{megginson2012introduction} that $A^s_{p,q}(\Omega)$ is a realisation of the dual space of $\tA^{-s}_{p',q'}(\Omega)$, with the pairing $\langle u,v\rangle_{A^s_{p,q}(\Omega)\times \tA^{-s}_{p',q'}(\Omega)}:=\langle U,v\rangle_{A^s_{p,q}(\R^n)\times A^{-s}_{p',q'}(\R^n)}$ for $u\in A^s_{p,q}(\Omega)$ and $v\in \tA^{-s}_{p',q'}(\Omega)$, where $U$ is any element of $A^s_{p,q}(\R^n)$ such that $U|_\Omega=u$. 
Similarly, by \cite[Theorem~1.10.17]{megginson2012introduction}, $\tA^{-s}_{p',q'}(\Omega)$ is a realisation of the dual space of $A^s_{p,q}(\Omega)$, with the pairing $\langle u,v\rangle_{\tA^{-s}_{p',q'}(\Omega)\times A^s_{p,q}(\Omega)}:=\overline{\langle v,u\rangle_{A^s_{p,q}(\Omega)\times \tA^{-s}_{p',q'}(\Omega)}}$, where the overline denotes complex conjugation.

\subsection{Multipliers, extensions and density results}
\label{sec:Multipliers}

{Following \cite[\S4.2.1, Definition 1]{RunstSickel}, given $f, g \in \mathcal{S}'(\mathbb{R}^n)$, we define the product $f \cdot g$ by
\begin{equation}\label{eq:product-def}
f \cdot g := \lim_{j \to \infty} S^j f \cdot S^j g \quad \text{in } \mathcal{S}'(\mathbb{R}^n),
\end{equation}
whenever this limit exists, where $S^j f := \mathcal{F}^{-1}[\psi(2^{-j}\,\cdot\,)\,\mathcal{F}f]$, $j = 0, 1, \ldots$, are the de la Vall\'ee-Poussin means of $f$ (with $\cF$ and $\cF^{-1}$ the Fourier transform and its inverse), and $\psi \in \mathcal{S}(\mathbb{R}^n)$ is a fixed function satisfying $0 \leq \psi \leq 1$, $\psi(\xi) = 1$ for $|\xi| \leq 1$, and $\psi(\xi) = 0$ for $|\xi| \geq 3/2$. Each $S^j f$ is an entire analytic function of exponential type (by the Paley--Wiener theorem), so the pointwise products $S^j f \cdot S^j g$ are well-defined, and $S^j f \to f$ in $\mathcal{S}'(\mathbb{R}^n)$ as $j \to \infty$.}

The following result is a consequence of 
{Theorem} 
\ref{lem:newDt} and results in \cite{FrazierJawerth90,Sickel99a} (in particular, \cite[Theorem 4.4]{Sickel99a}). The statement that $\chi_\Omega$ is a pointwise multiplier for $A^s_{p,q}(\R^n)$ means 
that the product $\chi_\Omega f$ is 
well-defined
 for $f\in A^s_{p,q}(\R^n)$, 
{in the above sense,} 
 as a product of tempered distributions\footnote{
If $1<p<\infty$ and $f\in A^s_{p,q}(\R^n)$, then $\chi_\Omega f$ is well-defined in this sense, and 
this definition coincides with the ordinary pointwise product of functions, whenever $s$ and $q$ are such as to guarantee that $f$ {\em is} a function, i.e., that $f\in L_1^{\rm loc}(\R^n)$. By \cite[\S2.2.4, Theorem 2]{RunstSickel} this holds if and only if $s>0$, or $A=F$ with $s=0$ and $q\leq 2$, or $A=B$ with $s=0$ and $q\leq \min(p,2)$. That $\chi_\Omega f$ is well-defined in the sense of  \cite[\S4.2.1, Definition 1]{RunstSickel} in these cases can be seen by applying part (ii) of \cite[\S4.2.1, Proposition]{RunstSickel}, since, where $p'$ is the usual conjugate exponent of $p$,  $\chi_\Omega\in L_{p'}^{\rm loc}(\R^n)$ and satisfies the bound (16) of that proposition, and, in each of the indicated cases, $f\in A^s_{p,q}(\R^n)\subset L_p(\R^n)$ (see \cite[\S2.2.4, Corollary 1]{RunstSickel}). Further, in these cases, $\chi_\Omega f$, understood as an ordinary pointwise product of functions, is also well-defined as a tempered distribution, since $\chi_\Omega f\in L_p(\R^n)$. That these two notions of the product coincide, in the first instance as distributions acting on $C_0^\infty(\R^n)$, follows from part (i) of \cite[\S4.2.1, Proposition]{RunstSickel}.}, and that $f\mapsto \chi_\Omega f$ defines a bounded linear operator from $A^s_{p,q}(\R^n)$ to itself. 
We exclude the case $\overline{\Omega}=\R^n$ in the following proposition because if $\overline{\Omega}=\R^n$ and $\dimA(\partial \Omega)<n$ then  $0=|\partial \Omega| = |\R^n\setminus \Omega|$ (see the discussion below \eqref{eq:IADefr}, where we note that $\dim_{\rm A}(\partial\Omega)=n$ if
$|\partial\Omega|>0$), so that $\chi_\Omega = 1$ almost everywhere, and $\chi_\Omega f = f$ for every tempered distribution $f$. Note that if the domain $\Omega$ satisfies 
$\overline{\Omega}\neq \R^n$ then 
$\dim_{\rm A}(\partial\Omega)\geq\dimH(\partial\Omega)\geq n-1$ (see the proof of Proposition \ref{prop:nset}\eqref{vi}).

\begin{prop}
\label{prop:pointmult}
Suppose that $\Omega\subset\R^n$ is  a domain with $\overline{\Omega}\neq \R^n$ and $\dim_{\rm A}(\partial\Omega)<n$. 
Then $\chi_{\Omega}$ is a pointwise
multiplier for $F_{p,q}^{s}\Rn $ if
\begin{equation}
1<p<\infty,\quad1\leq q\leq\infty,\quad\delta\left(\frac{1}{p}-1\right)<s<\frac{\delta}{p},\label{eq:sufcondmultF}
\end{equation}
and for $B_{p,q}^{s}\Rn $ if
\begin{equation}
1<p<\infty,\quad0<q\leq\infty,\quad\delta\left(\frac{1}{p}-1\right)<s<\frac{\delta}{p},\label{eq:sufcondmultB}
\end{equation}
where
\begin{equation}\label{delta}
\delta:= 
 n-\dimA(\partial \Omega)\in (0,1].
\end{equation}
\end{prop}

\begin{proof}
By 
{Theorem} \ref{lem:newDt}, $\Omega$ belongs to the 
 class 
$\cD^t$ for every 
$0< t<\delta$.
It follows then by \cite[Theorem 4.4]{Sickel99a} (see also \cite{FrazierJawerth90}) that $\chi_{\Omega}$ is a pointwise multiplier for $F_{p,q}^{s}\Rn $
under the conditions (\ref{eq:sufcondmultF}). The corresponding
result for $B_{p,q}^{s}\Rn $, under the conditions (\ref{eq:sufcondmultB}),
follows by real interpolation (using \cite[Theorem 2.4.2.(ii)]{Triebel83ThFS}).
\end{proof}


In the case that $\Omega$ is a bounded Lipschitz domain, in which case $\dimA(\Omega)=n-1$ and $\delta=1$, it is known that the above result is sharp in terms of its $s$-dependence, for the specified ranges of $p$ and $q$.
This is captured by the following proposition, largely taken from \cite{Tri:02}.
\begin{prop} \label{prop:Lip}
If $\Omega$ is a bounded Lipschitz domain then, for the ranges of $p$ and $q$ indicated in \eqref{eq:sufcondmultF} and \eqref{eq:sufcondmultB} for the cases $A=F$ and $A=B$, respectively, $\chi_\Omega$ is a pointwise multiplier for $A^s_{p,q}(\R^n)$ if and only if
\begin{equation} \label{eq:Lcase}
\frac{1}{p}-1<s<\frac{1}{p}.
\end{equation}
\end{prop}
\begin{proof} 
The ``if'' part of this result follows by combining Triebel's \cite[Proposition 5.3]{Tri:02} with his \cite[Proposition 5.1]{Tri:02} (or see \cite[Corollary 13.6]{FrazierJawerth90} for the case $A=F$ or \cite[\S4.6.3, Theorem 2]{RunstSickel} for the case $A=B$ for the larger class of domains with locally finite perimeter). The same results from \cite{Tri:02} (note that \cite[Proposition 5.3]{Tri:02} restricts attention, in the case $1<p<\infty$, only to $s$ in the range $(-1,1)$) show that $\chi_\Omega$ is a pointwise multiplier only if \eqref{eq:Lcase} holds or $|s|\geq 1$. But if $\chi_\Omega$ were a pointwise multiplier for some $s$ with $|s|\geq 1$, with $p$ and $q$ in the respective ranges \eqref{eq:sufcondmultF} and \eqref{eq:sufcondmultB}, it would follow by real interpolation (e.g., \cite[Theorems 2.4.2.(i), (ii)]{Triebel83ThFS}) that $\chi_\Omega$ would be a pointwise multiplier for $B^s_{p,q}(\R^n)$ for some $s\in (-1,1)$ outside the range \eqref{eq:Lcase}, and some $p,q\in (1,\infty)$, contradicting \cite[Proposition 5.3]{Tri:02}.
\end{proof}
 

The next result is a partial converse of Proposition \ref{prop:pointmult}, valid for Lipschitz domains\footnote{If $\Omega$ is a bounded Lipschitz domain then $\Omega$ is an open $n$-set, $\Omega^c$ is an $n$-set,  and $\partial \Omega$ is a  $d$-set with $d=n-1$, so that $\dimH(\partial \Omega)=\dim_{\rm M}(\partial \Omega)=\dimA(\partial \Omega)=n-1$; recall \eqref{eq:uMin} and the subsequent discussion. Hence Proposition \ref{lem:sharp} applies in this case, with $\delta'=1$, but gives a weaker conclusion than that provided by the ``only if'' part of Proposition \ref{prop:Lip}, since Proposition \ref{lem:sharp} applies only for $p,q\in(1,\infty)$ and does not rule out the endpoint cases $s=\delta/p$ and $s=\delta(p^{-1}-1)$.} but also for  certain non-Lipschitz domains, including the case where $\Omega$ is the interior or the complement of an $n$-attractor (see Corollary \ref{cor:pointmult}).

\begin{prop} \label{lem:sharp}
Suppose that $\Omega \subsetneqq\R^n$ is an open $n$-set such that $\partial \Omega$ is bounded,  $\dimA(\partial \Omega)<n$, and $\Omega^c$ is an $n$-set. 
 Suppose that $\chi_\Omega$ is a pointwise multiplier for $A^s_{p,q}(\R^n)$, for some $s\in \R$ and some $p,q\in (1,\infty)$.
 Then
$$
\delta'\left(\frac{1}{p}-1\right)\leq s\leq \frac{\delta'}{p},
$$
where 
$\delta':= n-\overline{\dim}_{\rm M}(\partial \Omega)\in (0,1]$.
\end{prop}
\begin{proof}
Note that $\partial \Omega=\partial F$, where $F:= \Omega^c$. Since $\dimA(\partial \Omega)<n$, $|\partial \Omega|=0$ so that, as discussed below \eqref{eq:nset}, $F^\circ$ is an open $n$-set and $\overline{F^\circ}=F$, so that $\overline{\Omega}\neq \R^n$. As noted above, this implies that $\dimH(\partial \Omega)\geq n-1$ so that $\delta'\in (0,1]$ follows from \eqref{eq:uMin}. One of $F^\circ$ and $\Omega$ is bounded, since $\partial \Omega$ is bounded. We assume, without loss of generality, that $\Omega$ is bounded. The assumption that $\Omega$ is an open $n$-set and $F$ is an $n$-set imply that $\partial \Omega$ is weakly regular of order $(0,0)$ in the sense\footnote{There is an error in the typesetting of (4.1) and (4.2) in  \cite{Sickel99a}. Inspecting the proof of the proposition below (4.2) it is clear that the backslash in these equations should be a minus sign.} of \cite[\S4.1]{Sickel99a} (to see this use \eqref{eq:dset} with $d=n$ to see that (4.1) and (4.2) in  \cite[\S4.1]{Sickel99a} hold with $\alpha=\beta=0$ for some sufficiently large $M$). 

Suppose that $\chi_\Omega$ is a pointwise multiplier of $A^{s}_{p,q}(\R^n)$, for some $p,q\in (1,\infty)$ 
and some $s>\delta'/p$. Then, since $\chi_\Omega$ is a pointwise multiplier for  $A^0_{p,q}(\R^n)$ by Proposition \ref{prop:pointmult}, it follows by interpolation, using \cite[Theorems 2.4.2.(i), (ii)]{Triebel83ThFS}, that $\chi_\Omega$ is a pointwise multiplier for  $B^{s'}_{p,p}(\R^n)=F^{s'}_{p,p}(\R^n)$, for some $s'> \delta'/p$. Applying \cite[Proposition 4.1]{Sickel99a}, we deduce that $\max(\overline{D}(\partial F),\overline{D}(\partial \Omega))\leq n-s'p$, where, for a set $E\subset \R^n$, $\overline{D}(\partial E)$ denotes the upper Minkowski dimension of $\partial E$ relative to $E$, as defined in \cite[(3.7)]{Sickel99a}. But it is easy to see from this definition that $\overline{\dim}_{\rm M}(\partial \Omega)= \max(\overline{D}(\partial F),\overline{D}(\partial \Omega))$, for any domain $\Omega\subsetneqq \R^n$ for which $\partial \Omega$ is bounded, where $F:= \Omega^c$. Thus 
$$
\overline{\dim}_{\rm M}(\partial \Omega) \leq n-s'p < n-\delta' = \overline{\dim}_{\rm M}(\partial \Omega),
$$
contradicting our assumption that $s>\delta'/p$.

Again assuming $p,q\in (1,\infty)$, 
suppose instead that $\chi_\Omega$ is a pointwise multiplier for $A^{s}_{p,q}(\R^n)$, for some $s<\delta'(p^{-1}-1)=-\delta'/p'$. (We denote by $p'$ and $q'$ the usual conjugate exponents.) It follows by duality (see \cite[Theorem 2.11.2]{Triebel83ThFS}), using that $\cS(\R^n)$ is dense in both $A^{s}_{p,q}(\R^n)$ and $A^{-s}_{p',q'}(\R^n)$ if $p,q\in (1,\infty)$ (e.g.,  \cite[\S2.1.1, Convention 1, \S2.1.3, Proposition 1]{RunstSickel}), that $\chi_\Omega$ is a pointwise multiplier for $A^{-s}_{p',q'}(\R^n)$. 
But by the argument in the above paragraph this implies that $-s\leq \delta'/p'$, i.e., $s\geq -\delta/'p'$, contradicting our assumption that $s< \delta' (p^{-1}-1)$.
\end{proof}


\begin{cor}
\label{cor:pointmult}
Let $\Gamma\subset \R^n$ be an $n$-attractor and let $\Omega$ denote either $\Gamma^\circ$ or $\Gamma^c$.
Then
$$
\delta:= n- \dimA(\partial \Omega)=n-\dimH(\partial \Omega){\in (0,1]}.
$$
{If} 
\eqref{eq:sufcondmultF} or \eqref{eq:sufcondmultB} hold, 
{then} $\chi_{\Omega}$ is a pointwise multiplier for $F_{p,q}^{s}\Rn $, respectively for $B_{p,q}^{s}\Rn $. Conversely, if $s\in \R$, $p,q\in (1,\infty)$,  and $\chi_\Omega$ is a pointwise multiplier for $F^s_{p,q}(\R^n)$ or $B^s_{p,q}(\R^n)$, 
then 
$$
\delta\left(\frac{1}{p}-1\right)\leq s\leq \frac{\delta}{p}.
$$
\end{cor}
\begin{proof}
{By Proposition \ref{prop:nset}\eqref{v},\eqref{vi} we have that $\dimH(\partial \Gamma)=\dim_{\mathrm{M}}(\partial\Gamma)=\dimA(\partial \Gamma)\in[n-1,n)$ and $\partial \Omega=\partial \Gamma$, from which the first conclusion follows.}
The {second conclusion} follows from Proposition \ref{prop:pointmult}, and the {third conclusion follows} from Proposition \ref{lem:sharp}, noting that $\Omega$ is an open $n$-set and $\Omega^c$ an $n$-set by  Proposition \ref{prop:nset}\eqref{vi},\eqref{vii},\eqref{ix}, since the closure of an open $n$-set is an $n$-set.
\end{proof}

\begin{rem}[Comparison of Proposition \ref{prop:pointmult} to the related literature]
\label{rem:multipliers}
As the above proof of Proposition \ref{prop:pointmult} makes clear, the proposition, as it applies to the spaces $F_{p,q}^s\Rn$, is essentially an immediate corollary of \cite[Theorem 4.4]{Sickel99a}, which establishes that $\chi_\Omega$ is a  pointwise multiplier if $\Omega\in \cD^t$ and the parameters $p,q,s$ satisfy \eqref{eq:sufcondmultF} with $\delta$ replaced by $t$. Conversely, in the case that $\partial \Omega$ is bounded, our Proposition \ref{prop:pointmult} implies \cite[Theorem 4.4]{Sickel99a}, with the help of 
{Theorem}
\ref{lem:newDt}(ii),  except that Sickel also states a multiplier result for the case $p=1$. In turn the multiplier result stated as \cite[Theorem 4.4]{Sickel99a} is (as noted in \cite[Remark 4.5]{Sickel99a}) a special case of \cite[Theorem 13.3, Corollary 13.5] {FrazierJawerth90}, where the class $\cD^t$ of \cite[Theorem 4.4]{Sickel99a} is replaced by a class $\cD_t$ of domains $\Omega$ that satisfy \eqref{eq:D_t} but with the integral over $B(x,r)\setminus \partial \Omega$ replaced by one over the smaller set $B(x,r)\cap \Omega$.  

As noted already in the discussion around Proposition \ref{prop:Lip}, Proposition \ref{prop:pointmult} fits well with what is known classically in the case where
$\Omega$ is a bounded Lipschitz domain, 
in which case $\dim_{\rm A}(\partial\Omega)=\dimH(\partial\Omega)=n-1$ and $\delta=1$.
See Proposition \ref{prop:Lip} or the more general \cite[Remark 5.4]{Tri:02}, which gives the multiplier property for $\max(n(1/p-1),1/p-1)<s<\min(1/p,1)$, with $0<p<\infty$ and $\min(1,p)\leq q\leq\infty$ if $A=F$ and with $0<p\leq \infty$ and $0<q\leq\infty$ if $A=B$.
Compare also with \cite[Corollaries 3 and 4]{Triebel2003}, where related results are presented for non-Lipschitz cases. In particular,
part (ii) of \cite[Corollary 4]{Triebel2003} gives that if
we merely assume that a bounded domain $\Omega$ has a porous boundary
then $\chi_{\Omega}$ is a pointwise multiplier for $F_{p,q}^{0}\Rn $
when $1<p<\infty$, $1\leq q\leq\infty$. Recalling Proposition \ref{prop:nset}\eqref{xi}, this provides another route to deducing the result for $\Omega=\Gamma^\circ$ for the $F^0_{p,q}$ spaces in Proposition \ref{prop:pointmult}, {in the case where $\Gamma$ is an $n$-attractor}. 

Again assuming that $\Omega$ is bounded, part (i) of \cite[Corollary 4]{Triebel2003} states, among other things, that $\chi_\Omega$ is a pointwise multiplier for $F^s_{p,q}(\R^n)$ if \eqref{eq:sufcondmultF} holds with $\delta$ replaced by some $\varepsilon\in (0,n)$ (whose value is not specified in the corollary or its proof), under the assumption that $\partial \Omega$ is uniformly porous in the sense of Definition \ref{def:up}, an assumption that, by Proposition \ref{prop:upA}, is stronger than our assumption that $\dimA(\partial \Omega)<n$.  
Note that $\partial \Omega$ being uniformly porous means that $\partial \Omega$ is a porous set that is an $h$-set, for some $h:[0,1]\to [0,1]$ that satisfies the conditions of Definition \ref{def:hset}. By Proposition \ref{prop:hsetp} it follows, moreover, that \eqref{eq:po} holds for some $\lambda\in (0,n)$. Further, tracing the arguments that lead to \cite[Corollary 4]{Triebel2003} back to \cite[Remark 9]{Triebel2003}, it appears that the value of $\varepsilon$ in \cite[Corollary 4]{Triebel2003} is precisely the $\lambda\in (0,n)$ of \eqref{eq:po} (see also the discussion in the proof of \cite[Proposition 3.19]{Tri08}). Thus the result of part (i) of \cite[Corollary 4]{Triebel2003} is, among other things, that, under the assumption that $\partial \Omega$ is uniformly porous, $\chi_\Omega$ is a pointwise multiplier for $F^s_{p,q}(\R^n)$ if \eqref{eq:sufcondmultF} holds with $\delta$ replaced by the $\lambda$ of \eqref{eq:po}. But note that 
$\lambda \leq \delta = n-\dimA(\partial \Omega)$, by Proposition \ref{prop:upA}.
Thus the pointwise multiplier property of $\chi_\Omega$ for $F^s_{p,q}(\R^n)$ for the parameter ranges  \eqref{eq:sufcondmultF} is established in Proposition \ref{prop:pointmult}  under a weaker assumption on $\partial \Omega$ than the assumption of \cite{Triebel2003}, 
and 
for a range of $s$ that is no smaller than that in \cite{Triebel2003}. Of course this means, in view of the first paragraph above, that the pointwise multiplier results of \cite[Corollary 4(i)]{Triebel2003}, for $1<p<\infty$, are also a corollary of results in \cite{Sickel99a,FrazierJawerth90}.

\end{rem}

Proposition \ref{prop:pointmult} will be used in Proposition \ref{prop:extbyzero} to prove the existence of extension operators on function spaces on domains. 
First, however, we prove a preliminary result,
which might be of independent interest. In this lemma and the following propositions, for a domain $\Omega\subset \R^n$ and $\varphi\in{\cal D}(\Omega)$, $\widetilde{\varphi}$ denotes the extension of $\varphi$ 
by zero to $\R^n$.




\begin{lem}
\label{lem:restriction}
Let $f,g\in{\cal S}'\Rn $ and let $\Omega\subset \R^n$ be a domain. 
 Assume that the product $f\cdot g$ is well-defined in
the sense of 
\eqref{eq:product-def}
and
that the distributional restriction $f|_{\Omega}$ of $f$ to $\mathcal{D}'(\Omega)$
equals $0$. Then also $(f\cdot g)|_{\Omega}=0$.
\end{lem}

\begin{proof}
Using the notation in 
\eqref{eq:product-def}
we have the following,
for any $\varphi\in{\cal D}(\Omega)$: 
\[
\left\langle (f\cdot g)|_{\Omega},\varphi\right\rangle =\left\langle f\cdot g,\widetilde{\varphi}\right\rangle =\lim_{j\to\infty}\left\langle S^{j}f\cdot S^{j}g,\widetilde{\varphi}\right\rangle =\lim_{j\to\infty}\left\langle \left(S^{j}f\cdot S^{j}g\right)|_{\Omega},\varphi\right\rangle =0,
\]
where the last identity follows from \cite[Lemma 4.2.1/2 (p.~145)]{RunstSickel}.
\end{proof}

\begin{prop}
\label{prop:extbyzero}
Let $\Omega$, $p$,  and $s$
be as in Proposition \ref{prop:pointmult}. Let $1\leq q<\infty$ in
the $F$ case and $0<q<\infty$ in the $B$ case. 
{Given $f\in A_{p,q}^{s}(\Omega)$, and any $g\in A_{p,q}^{s}\Rn \subset\mathcal{S}'\Rn $
such that $f$ equals the distributional restriction $g|_{\Omega}$
of $g$ to $\mathcal{D}'(\Omega)$, the map}
\begin{equation}
{{\rm ext}: f\mapsto \chi_{\Omega}\cdot g}
\label{eq:extbyzero}
\end{equation}
defines a linear and bounded
extension operator from $A_{p,q}^{s}(\Omega)$ into $A_{p,q}^{s}\Rn $. 
Furthermore,
${\rm supp}({\rm ext}f)\subset\overline{\Omega}$,
for any $f\in A_{p,q}^{s}(\Omega)$, indeed ${\rm ext}$ is an isomorphism from $A_{p,q}^{s}(\Omega)$ to $A_{p,q,\overline{\Omega}}^{s}$. Moreover, the restriction operator $|_\Omega$ is an isomorphism from $A_{p,q,\overline{\Omega}}^{s}$ to $A_{p,q}^{s}(\Omega)$, with $|_\Omega^{-1}={\rm ext}$.
\end{prop}

\begin{proof}

Clearly, ${\rm ext}$ {given by \eqref{eq:extbyzero}} takes elements of $A_{p,q}^{s}(\Omega)$ to
elements of $A_{p,q}^{s}\Rn $, due to the pointwise multiplier property
in Proposition \ref{prop:pointmult}, but in order that this is a well-defined map
we need to show that the image of any given $f$ does not depend
on the $g$ chosen. Or, equivalently, that $\chi_{\Omega}\cdot h=0$
whenever $h\in A_{p,q}^{s}\Rn $ with $h|_{\Omega}=0$.

By assumption, $\dim_{\rm A}(\partial\Omega)<n$, so that
$|\partial\Omega|=0$ (see the discussion below \eqref{eq:IADefr}), 
and also $\dimH(\partial \Omega)\leq \dimA(\partial \Omega)$,   as noted above \eqref{eq:uMin}.  
We thus have, from \cite[Proposition 3.7]{caetano2019density} applied
to $\partial\Omega$ and under our hypotheses, that $\partial\Omega$
is $A_{p,q}^{s}$-null, i.e., $A_{p,q,\partial\Omega}^{s}=\{0\}$.
Given any $h\in A_{p,q}^{s}\Rn $ with $h|_{\Omega}=0$, in order
to show that $\chi_{\Omega}\cdot h=0$ it is then enough to guarantee
that ${\rm supp}(\chi_{\Omega}\cdot h)\subset\partial\Omega$, which
is what we show next.

Given any $\psi\in{\cal D}((\partial\Omega)^{c})$, we clearly have,
with the usual abuse of notation, that $\psi\chi_{\Omega}\in{\cal D}(\Omega)$
and $\psi\chi_{\Omega^{c}}\in{\cal D}(\overline{\Omega}^{c})$, from
which follows that $\left\langle \chi_{\Omega}\cdot h,\psi\right\rangle =\left\langle \chi_{\Omega}\cdot h,\psi\chi_{\Omega}\right\rangle +\left\langle \chi_{\Omega}\cdot h,\psi\chi_{\Omega^{c}}\right\rangle =0$, where
the first term is zero because $\left(\chi_{\Omega}\cdot h\right)|_{\Omega}=0$
(as follows from Lemma \ref{lem:restriction} and the hypothesis $h|_{\Omega}=0$)
and the second term is zero because $\left(\chi_{\Omega}\cdot h\right)|_{\overline{\Omega}^{c}}=0$
(as follows again by Lemma \ref{lem:restriction}, now read with $\overline{\Omega}^{c}$
instead of $\Omega$ and using the fact that $\chi_{\Omega}|_{\overline{\Omega}^{c}}=0$).

The linearity of ${\rm ext}$ is clear. On the other hand, given any
$f\in A_{p,q}^{s}(\Omega)$ and any $g\in A_{p,q}^{s}\Rn $ such
that $f=g|_{\Omega}$, using the fact that $\chi_{\Omega}$ is a pointwise
multiplier for $A_{p,q}^{s}\Rn $ we get that
\[
\|{\rm ext}f\|_{A_{p,q}^{s}\Rn }=\|\chi_{\Omega}\cdot g\|_{A_{p,q}^{s}\Rn }\lesssim\|g\|_{A_{p,q}^{s}\Rn },
\]
proving boundedness of ${\rm ext}$. To prove ${\rm ext}$ is indeed an extension we need to show that $({\rm ext}f)|_{\Omega}=f$, i.e., that $(\chi_{\Omega}\cdot g)|_{\Omega}=f$. Recalling that ${\cal S}(\R^n)$ is dense in $A^s_{p,q}(\R^n)$ for the ranges of $p$ and $q$ we are considering (e.g., \cite[\S2.3.3]{Triebel83ThFS}), consider a sequence $(\varphi_{j})_{j\in\N}\subset\mathcal{S}\Rn$ 
such that $\varphi_{j}\to g$ in $A_{p,q}^{s}\Rn$, so that, consequently,
$\chi_{\Omega}\cdot\varphi_{j}\to\chi_{\Omega}\cdot g$ in $A_{p,q}^{s}\Rn$
and, therefore, $(\chi_{\Omega}\cdot\varphi_{j})|_{\Omega}\to(\chi_{\Omega}\cdot g)|_{\Omega}$
in ${\cal D}'(\Omega)$. 
Then, given any $\varphi\in{\cal D}(\Omega$),
we have, by \cite[Lemma 4.2.1/1 (p.~144)]{RunstSickel},
\begin{align*}
\left\langle (\chi_{\Omega}\cdot g)|_{\Omega},\varphi\right\rangle  & =\lim_{j\to\infty}\left\langle (\chi_{\Omega}\cdot\varphi_{j})|_{\Omega},\varphi\right\rangle \\
= & \lim_{j\to\infty}\left\langle \chi_{\Omega}\cdot\varphi_{j},\widetilde{\varphi}\right\rangle =\lim_{j\to\infty}\int_{\R^n}\chi_{\Omega}(x)\varphi_{j}(x)\widetilde{\varphi}(x)\,\rd x\\
= & \lim_{j\to\infty}\left\langle \varphi_{j}|_{\Omega},\varphi\right\rangle =\left\langle g|_{\Omega},\varphi\right\rangle =\left\langle f,\varphi\right\rangle .
\end{align*}

That ${\rm supp}({\rm ext}f)\subset\overline{\Omega}$ for
any $f\in A_{p,q}^{s}(\Omega)$ follows again by Lemma \ref{lem:restriction} read
with $\overline{\Omega}^{c}$ instead of $\Omega$ and using the fact
that $\chi_{\Omega}|_{\overline{\Omega}^{c}}=0$. To see that ${\rm ext}:A_{p,q}^{s}(\Omega)\to A_{p,q,\overline{\Omega}}^{s}$ is an isomorphism, it is enough to show that its left inverse $|_\Omega:A_{p,q,\overline{\Omega}}^{s}\to A_{p,q}^{s}(\Omega)$ is a bijection, in which case $\rm ext$ is invertible and ${\rm ext}^{-1}=|_\Omega$. But $|_\Omega$ is surjective, since $|_\Omega\circ  {\rm ext}f=f$, for $f\in A_{p,q}^{s}(\Omega)$. Further, $|_\Omega$ is injective since, if $g\in A_{p,q,\overline{\Omega}}^{s}$ and $g|_\Omega  = 0$, then $g\in A_{p,q,\partial \Omega}^{s}$ and, as noted above, $A_{p,q,\partial \Omega}^{s}=\{0\}$.
\end{proof}

\begin{rem}[Relationship of Proposition \ref{prop:extbyzero} to previous results]
\label{rem:extbyzero}
The final statement of Proposition \ref{prop:extbyzero} generalises standard results about $p=2$ Sobolev spaces on Lipschitz domains, for which it is known that $\tH^s(\Omega)=H^s_{\overline\Omega}$ for all $s\in\R$ \cite[Theorem~3.29]{McLean} and that $|_\Omega$ is an isomorphism from $\tH^s(\Omega)=H^s_{\overline\Omega}$ to $H^s(\Omega)$ for $|s|<1/2$ (see, e.g., \cite[Corollary~3.29(ix) \& Lemma~3.31(iii)]{ChaHewMoi:13}). 


\end{rem}

\begin{prop}
\label{prop:tilde-subscript}
{Let $\Omega$,  $p$,
$q$ and $s$ be as in Proposition \ref{prop:extbyzero}, and suppose that $\partial\Omega$ is bounded. Then} 
\[
\tA_{p,q}^{s}(\Omega)=A_{p,q,\overline{\Omega}}^{s}.
\]
\end{prop}

\begin{proof}
Clearly, we only need to prove that $A_{p,q,\overline{\Omega}}^{s}\subset\tA_{p,q}^{s}(\Omega)$.

Given $f\in A_{p,q,\overline{\Omega}}^{s}$, by definition we have
that $f\in A_{p,q}^{s}\Rn $ and ${\rm supp}f\subset\overline{\Omega}$.

In particular, $f|_{\Omega}\in A_{p,q}^{s}(\Omega)$, and we claim
that it also holds that $f|_{\Omega}\in\mathring{A}_{p,q}^{s}(\Omega)$,
where $\mathring{A}_{p,q}^{s}(\Omega)$ stands for the closure of
${\cal D}(\Omega)$ in $A_{p,q}^{s}(\Omega)$. This claim follows
by applying \cite[Proposition 2.2(i)]{Ca:00},\footnote{{While \cite[Proposition 2.2]{Ca:00} is stated under the assumption that $\Omega$ is bounded, 
part (i) of \cite[Proposition 2.2]{Ca:00} holds also for unbounded $\Omega$, as long as $\partial \Omega$ is bounded. 
This can be proved by the same argument given in \cite{Ca:00}, except that one should replace ${\cal S}$ by ${\cal D}$ throughout (in the notation of \cite{Ca:00}).}}
noting that our hypothesis that $\dimA(\partial \Omega) <n$ guarantees that also $\overline{\dim}_{\rm M}(\partial \Omega)<n$, by  \eqref{eq:uMin}.


Consider then a sequence $(\varphi_{k})_{k}\subset{\cal D}(\Omega)$
converging to $f|_{\Omega}$ in $A_{p,q}^{s}(\Omega)$. We claim that
$(\widetilde{\varphi_{k}})_{k}$ converges to $f$ in $A_{p,q}^{s}\Rn $,
which will finish our proof. Since our $\partial\Omega$ is $A_{p,q}^{s}$-null
(as mentioned in the proof of Proposition \ref{prop:extbyzero}), from \cite[Lemma 4.4]{caetano2019density}
it follows that the restriction operator is an isometric isomorphism
from $A_{p,q,\overline{\Omega}}^{s}$ onto the space $RA_{p,q,\overline{\Omega}}^{s}:=\{f\in{\cal D}'(\Omega):\:f=g|_{\Omega}$
for some $g\in A_{p,q,\overline{\Omega}}^{s}\}$ endowed with the
corresponding quotient quasi-norm (which is in general stronger than
the quasi-norm in $A_{p,q}^{s}(\Omega)$). Therefore,
\begin{equation}
\|\widetilde{\varphi_{k}}-f\|_{A_{p,q}^{s}\Rn }=\|\varphi_{k}-f|_{\Omega}\|_{RA_{p,q,\overline{\Omega}}^{s}}.\label{eq:A-RA}
\end{equation}
On the other hand, given any $h\in RA_{p,q,\overline{\Omega}}^{s}$,
clearly $h\in A_{p,q}^{s}(\Omega)$ too, therefore, with ${\rm ext}$
the extension from Proposition \ref{prop:extbyzero}, ${\rm ext}\,h\in A_{p,q,\overline{\Omega}}^{s}$,
hence
\begin{equation}
\|h\|_{RA_{p,q,\overline{\Omega}}^{s}}\leq\|{\rm ext}\,h\|_{A_{p,q}^{s}\Rn }\lesssim\|h\|_{A_{p,q}^{s}(\Omega)}\label{eq:RA-A(Omega)}
\end{equation}
(so actually in our setting the quasi-norm in $RA_{p,q,\overline{\Omega}}^{s}$
 is equivalent to the restriction to this space of the quasi-norm in
$A_{p,q}^{s}(\Omega)$). From (\ref{eq:A-RA}), (\ref{eq:RA-A(Omega)})
and the definition of $(\varphi_{k})_{k}$, our claim follows.
\end{proof}


{The following theorem
follows from 
\cite[Corollary 4.17]{caetano2019density} in the case $s\neq 0$, and the above proposition in the case $s=0$.}

\begin{thm}
\label{cor:tildesubscriptGeneral}
Let $\Omega\subset\R^n$ be a thick domain, with $\partial\Omega$ bounded and $\dim_{\rm A}(\partial\Omega)<n$.
Then  
\[
\tA_{p,q}^{s}(\Omega)=A_{p,q,\overline{\Omega}}^{s},
\qquad s\in\R,\,\,1<p,q<\infty.
\]
\end{thm}

\begin{proof}
For $s\not=0$, the conclusion follows from \cite[Corollary 4.17]{caetano2019density},
because of the thickness of 
{$\Omega$ and the fact that $|\partial\Omega|=0$, as noted in the proof of Proposition \ref{prop:extbyzero}.}  
 (Note that in the Sobolev space case  $H^s_p(\R^n)$ this argument works for all $s\in \R$, using \cite[Corollary 4.18]{caetano2019density} instead of \cite[Corollary 4.17]{caetano2019density}.) 
For $s=0$, the conclusion follows from Proposition \ref{prop:tilde-subscript}.
\end{proof}

For convenience we state separately the conclusion of Theorem \ref{cor:tildesubscriptGeneral} for the {$n$-attractor case, highlighting also the Sobolev space case. As noted in the proof of Theorem \ref{cor:tildesubscriptGeneral}, in the Sobolev space setting the case $s=0$ is already covered by \cite[Corollary~4.18]{caetano2019density}, so Proposition \ref{prop:tilde-subscript} is not required.} 
\begin{cor}
\label{cor:tildesubscriptHs}
Let $\Gamma\subset\R^n$ be an $n$-attractor.
Then
\begin{align}
\label{eq:tildesubscript}
{\tA_{p,q}^{s}(\Gamma^{\circ})=A_{p,q,\Gamma}^{s}\qquad\mbox{and}\qquad \tA_{p,q}^{s}(\Gamma^{c})=A_{p,q,\overline{\Gamma^{c}}}^{s},
\qquad s\in\R,\,\,1<p,q<\infty,}
\end{align}
and, in particular,
\begin{align}
\label{eq:tildesubscriptHs}
\tH^{s}_{p}(\Gamma^{\circ})=H_{{p,}\Gamma}^{s}\qquad\mbox{and}\qquad \tH^{s}_{p}(\Gamma^{c})=H_{{p,}\overline{\Gamma^{c}}}^{s},\qquad s\in\R, \quad { 1<p<\infty}.
\end{align}
\end{cor}
\begin{proof}
Since $\overline{\Gamma^\circ}=\Gamma$ (which holds by Proposition \ref{prop:nset}\eqref{vi}), the results follow by applying Theorem \ref{cor:tildesubscriptGeneral} with $\Omega:=\Gamma^\circ$ and then with $\Omega:=\Gamma^c$. For both cases, the thickness of $\Omega$ was established in Theorem \ref{thm:thick}, 
the fact that $\dim_{\rm A}(\partial\Omega)<n$ holds by Proposition \ref{prop:nset}\eqref{v}, and the fact that $\partial\Omega=\partial\Gamma$ is bounded holds because $\Gamma$ (and hence also $\partial\Gamma$) is compact.
\end{proof}

\begin{rem}[Equality of Sobolev spaces on open $n$-sets]
\label{re:tildesubscript}
In the case when $\Gamma=\overline{\Omega}$, for an open $n$-set $\Omega$ (equivalently, see the discussion following \eqref{eq:nset}, $\Gamma$ is an $n$-set and $|\partial \Gamma|=0$), one can prove results of independent interest connected to those in \cite{caetano2019density} and Corollary \ref{cor:tildesubscriptHs} above (the proof of which relied on \cite[Corollary~4.18]{caetano2019density}), that do not require thickness. 

Suppose that $\Omega$ is an open $n$-set (which implies that $|\partial\Omega| = 0$). 
Then we claim that
\begin{equation}
\widetilde{H}_p^{s}(\overline{\Omega}^{c}) = H_{p,\Omega^{c}}^{s},\qquad 0\leq s\leq 1,\quad 1<p<\infty.\label{eq:s>0}
\end{equation}
To prove this, we note that the inclusion $\widetilde{H}^{s}_{p}(\overline{\Omega}^{c})\subset H_{p,\Omega^{c}}^{s}$ is trivial, so it remains to show that $H_{p,\Omega^{c}}^{s}\subset \widetilde{H}^{s}_{p}(\overline{\Omega}^{c})$. For this, we note that if $s\geq 0$ and $u\in H^s_{p,\Omega^c}$ then $u=0$ almost everywhere on $\Gamma=\overline\Omega$ with respect to Lebesgue measure. 
Then, for $0<s\leq 1$, 
$u$ is an element of 
$\ker([\rm Tr_{\overline{\Omega},0}]_\mu)$, in the notation of \cite{Hinz} (with $\mu$ denoting Lebesgue measure restricted to $\overline{\Omega}$), 
which equals $\widetilde{H}^{s}_{p}(\overline{\Omega}^{c})$ for $0<s\leq 1$, by \cite[Corollary 3.2]{Hinz} (applied with $d=n$ and $\alpha=s\in (0,1]$, in which case $m_t=0$ in \cite[(23)]{Hinz}). 
For $s=0$ the fact that $u\in \widetilde{H}_{p}^{0}(\overline{\Omega}^{c})$ follows from the density of $C^\infty_0(\overline{\Omega}^c)$ in $L_p(\overline{\Omega}^c)$ (see, e.g.,\ \cite{Adams}). 
From \eqref{eq:s>0} and
\cite[Lemma~4.14]{caetano2019density} 
it follows that
\begin{equation}
	\widetilde{H}_{p}^{s}(\Omega)=H_{p,\overline{\Omega}}^{s},\qquad -1\leq s\leq 0,\quad 1<p<\infty,\label{eq:-s<0}
\end{equation}
a result of the type \cite[Theorem 4.16]{caetano2019density}, where $\widetilde{H}^{s}_p(\Omega)=H_{p,\overline{\Omega}}^{s}$ is proven for $s<0$ under the assumption that $|\partial\Omega|=0$, $\Omega$ is $I$-thick, and that $(\overline{\Omega})^{\circ}=\Omega$.

Suppose, alternatively, that $\Omega$ is open and  $\overline{\Omega}^{c}$
is an open $n$-set. Then, instead of (\ref{eq:s>0}), noting that $\overline{\overline{\Omega}^{c}}^c=(\overline{\Omega})^{\circ}$, we obtain by similar reasoning that
\[
\widetilde{H}^{s}_p((\overline{\Omega})^{\circ})=H_{p,\overline{\Omega}}^{s},\qquad 0\leq s\leq 1,\quad 1<p<\infty.
\]
Suppose further that
$(\overline{\Omega})^{\circ}\setminus\Omega$ is $(-s,p')$-null (in the sense of \cite[Definition~2.1]{HewMoi:15}, meaning that $H^{-s}_{p',F}=0$ for every closed set $F\subset (\overline{\Omega})^{\circ}\setminus\Omega$), which holds in particular if $(\overline{\Omega})^{\circ}=\Omega$. Then,
by 
\cite[Proposition~2.11]{HewMoi:15} (applied with $F_1=\Omega^c$ and $F_2=((\overline{\Omega})^\circ)^c$, in the notation of \cite{HewMoi:15}) and \cite[Proposition~3.5]{caetano2019density}, it follows 
that $\widetilde{H}^{s}_p((\overline{\Omega})^\circ)=\widetilde{H}^{s}_p(\Omega)$. 
Hence, if $(\overline{\Omega})^{\circ}\setminus\Omega$ is $(-1,p')$-null (e.g., if $(\overline{\Omega})^{\circ}=\Omega$), which implies that $(\overline{\Omega})^{\circ}\setminus\Omega$ is also $(-s,p')$-null for every $s\leq 1$, then 
combining with 
the above 
gives
\begin{equation}
\widetilde{H}^{s}_p(\Omega)=	H_{p,\overline{\Omega}}^{s},\qquad 0\leq s\leq 1,\quad 1<p<\infty,\label{eq:s>0compl}
\end{equation}
a result of the type
\cite[Corollary 4.13]{caetano2019density}, which provides $\widetilde{H}^{s}_p(\Omega)=H_{p,\overline{\Omega}}^{s}$ for {$s\geq 0$} under the assumption that $|\partial\Omega|=0$ and that $\Omega$ is $E$-thick (which implies that $(\overline{\Omega})^{\circ}=\Omega$, see~\cite[Proposition 4.7(v)]{caetano2019density}).

Combining the above results,
we obtain the following:
If both $\Omega$ and $\overline{\Omega}^{c}$ are open $n$-sets (which implies that $|\partial\Omega|=0$) and $(\overline{\Omega})^{\circ}\setminus\Omega$ is $(-1,p')$-null 
(e.g., if $(\overline{\Omega})^{\circ}=\Omega$), then
\begin{equation}
	\widetilde{H}^{s}_p(\Omega)=H_{p,\overline{\Omega}}^{s},\qquad -1\leq s\leq 1,\quad 1<p<\infty,\label{eq:final}
\end{equation}
which should be compared with \cite[Corollary 4.18]{caetano2019density}, which provides $\widetilde{H}^{s}_p(\Omega)=H_{p,\overline{\Omega}}^{s}$ for all $s\in\R$ under the assumption that $|\partial\Omega|=0$ and $\Omega$ is thick (so that again $(\overline{\Omega})^{\circ}=\Omega$).

In Appendix \ref{sec:Appendix2} (see Theorem \ref{thm:exnset}) we construct examples of bounded open $n$-sets $\Omega$ with $\overline{\Omega}^\circ = \Omega$ that show that \eqref{eq:-s<0} does not hold for every open $n$-set $\Omega$ if $-1\leq s\leq 0$ is replaced by $-1\leq s\leq \epsilon$, for any $\epsilon>0$. We also, in Corollary \ref{cor:exnset}, exhibit examples  of bounded domains $\Omega$ with $\overline{\Omega}^\circ = \Omega$ and $\overline{\Omega}^c$ an open $n$-set, that show that \eqref{eq:s>0} does not hold for all open $n$-sets $\Omega$ 
if $0\leq s\leq 1$ is replaced by $-\epsilon\leq s\leq 1$ with $\epsilon>0$.
\end{rem}

\begin{rem} [Examples for which the above results apply]\label{rem:exIFSSobolev}
{Corollary \ref{cor:tildesubscriptHs} applies}, and hence the relations \eqref{eq:tildesubscript} and \eqref{eq:tildesubscriptHs} hold, whenever $\Gamma$ is an $n$-attractor, and hence for all the examples in Figure \ref{fig:2sets}. In particular, \eqref{eq:tildesubscript} and \eqref{eq:tildesubscriptHs} hold for the case of the standard Koch snowflake domain mentioned in Remark \ref{rem:exIFS}. This was proved already in \cite[Proposition~5.2]{caetano2019density}, for a more general class of snowflake domains, albeit not for the case $s=0$ for the general $A^s_{p,q}$ spaces, by proving that $\Omega$ is thick by a direct argument involving polygonal approximations (without using that $\Gamma$ is the attractor of an IFS).

{Corollary \ref{cor:tildesubscriptHs} also applies}, and hence the relations \eqref{eq:tildesubscript} and \eqref{eq:tildesubscriptHs} also hold, for the example \eqref{eq:Gammaeq}. In the Sobolev space case this follows alternatively,
but only for the restricted range 
$|s|\leq 1$,
by Remark \ref{re:tildesubscript} (with $\Omega=\Gamma^\circ= \bigcup_{k=0}^\infty (\alpha^{2k+1},\alpha^{2k})$, using the facts that 
$\Omega$ and $\overline{\Omega}^c$ are both open $n$-sets, {$(\overline\Omega)^\circ=\Omega$, and $\overline\Omega=\Gamma$, which hold} by Proposition \ref{prop:nset}),
and alternatively,
but only for the restricted range $|s|\leq 1/2$ and the case $p=2$, by \cite[Theorem~3.24]{ChaHewMoi:13}, and alternatively, but only for $s\in \N$, by \cite[\S11]{AdHe} (or see \cite[Lemma~3.20(i)]{ChaHewMoi:13} for the $p=2$ case).
\end{rem}

\begin{rem}[Example to which \eqref{eq:final} applies which is not thick] \label{rem:ntex}
An example of a domain $\Omega\subset \R$ to which Remark \ref{re:tildesubscript} applies but which is not $I$-thick (and so not thick), so that \cite[Corollary~4.18]{caetano2019density} does not apply to deduce \eqref{eq:final} for all $s\in \R$, is the set $\Omega=\bigcup_{k=1}^{\infty}\big(\frac{1}{2k+1},\frac{1}{2k}\big)$. For this example $\Omega = (\overline{\Omega})^\circ$ and one can show that $\Omega$  and $\overline{\Omega}^c$ are open $n$-sets, so that \eqref{eq:final} holds.  
For this example 
\eqref{eq:final} follows alternatively, but only for the restricted range $|s|\leq 1/2$ and the case $p=2$, from \cite[Theorem~3.24]{ChaHewMoi:13}, and alternatively, but only for $s\in \N$, by \cite[\S11]{AdHe}. 
{We note that for this example Proposition \ref{prop:tilde-subscript} does not apply because $\dim_{\rm A}(\partial\Omega)=1=n$ (see \cite[Theorem 2.1.1]{Fraser:21}).}
\end{rem}


\subsection{Interpolation {results}}
\label{sec:Interpolation}


Our next goal is to obtain {interpolation results involving} the function spaces considered above. We start with the
following auxiliary result, which establishes that the extension operator defined in the previous section via pointwise multiplication by a characteristic function enjoys a certain uniformity property. 
The term ``common extension operator'' in the following means that the distribution ${\rm ext}f$ does
not depend on the space $A_{p,q}^{s}(\Omega)$ where $f$ lives, provided $s$ and $q$ lie in ranges considered.
More precisely, temporarily denoting the extension operator (\ref{eq:extbyzero}) by ${\rm ext}_{p,q}^s:A_{p,q}^{s}(\Omega) \to A_{p,q}^{s}\Rn$, it means that if $s_{1},s_{2}$ lie in the range
considered for the parameter $s$ and $q_{1},q_{2}$ lie in the range considered for the parameter $q$ then ${\rm ext}_{p,q_1}^{s_1} f = {\rm ext}_{p,q_2}^{s_2} f$ for all $f\in A_{p,q_{1}}^{s_{1}}(\Omega)\cap A_{p,q_{2}}^{s_{2}}(\Omega)$. 
{We note that a discussion of extension operators and references to some relevant background literature on this topic can be found in \cite[Remark~3.3]{Tri08}.} 

\begin{prop}
\label{prop:common-pfixed} Let $\Omega$, $\delta$ and $p$
be as in Proposition \ref{prop:pointmult}. Then the map ${\rm ext}$
given by (\ref{eq:extbyzero}) is a common extension operator from
$A_{p,q}^{s}(\Omega)$ into $A_{p,q}^{s}\Rn $ as $s$ ranges over
the interval
\[
\left(\delta\left(\frac{1}{p}-1\right),\frac{\delta}{p}\right)
\]
and $q$ ranges over $[1,\infty)$ in the $F$ case and over $(0,\infty)$
in the $B$ case.
\end{prop}

\begin{proof}
We already know, as a by-product of the proof of Proposition \ref{prop:extbyzero},
that, for fixed $s$, $p$ and $q$, the product $\chi_{\Omega}\cdot g$
defining ${\rm ext}f$ does not depend on the $g$ chosen in $A_{p,q}^{s}\Rn $
and such that $f=g|_{\Omega}$. What we need to prove now is that,
for fixed $p$, it does also not depend on whatever choice we made
of $g\in A_{p,q}^{s}\Rn $ with $s$ and $q$ in the ranges given,
as long as $f=g|_{\Omega}$. Clearly, by what has been mentioned,
it will be enough to show that, for any $s_{1},s_{2}$ in the range
considered for the parameter $s$ and for any $q_{1},q_{2}$ in the
range considered for the parameter $q$, given any $f\in A_{p,q_{1}}^{s_{1}}(\Omega)\cap A_{p,q_{2}}^{s_{2}}(\Omega)$,
there exists $g\in A_{p,q_{1}}^{s_{1}}\Rn \cap A_{p,q_{2}}^{s_{2}}\Rn $
such that $f=g|_{\Omega}$ (also because the definition of pointwise
multiplication is considered at the level of the distributions, and
not of the particular $A$ spaces being considered --- 
cf.~\eqref{eq:product-def}).
And this is clearly true, due to the embeddings $A_{p,q_{1}}^{s_{1}}\Rn \hookrightarrow A_{p,q_{2}}^{s_{2}}\Rn $
either if $s_{1}=s_{2}$ and $q_{1}\leq q_{2}$ or if $s_{1}>s_{2}$.
\end{proof}

Our main interpolation result {in Theorem \ref{prop:wolff} below} is a partial 
generalisation   
of Triebel's result \cite[Theorem 4.17]{Tri08}\footnote{We note that, in the case of real interpolation, \cite[Theorem 4.17]{Tri08} covers the case where $q_0$, $q_1$ or $q$ equal $\infty$, which we do not consider here.}, which provides similar results, but with 
our assumption that $\dim_{\rm A}(\partial\Omega)<n$ replaced by the assumption that 
$\partial\Omega$ is a $d$-set for some $n-1\leq d<n$. 
In the case where $\partial\Omega$ is bounded, our assumption that $\dim_{\rm A}(\partial\Omega)<n$ is weaker than Triebel's $d$-set assumption (see the comments after \eqref{eq:uMin}, where we note that $\dimA(\partial \Omega)=d$ if $\partial \Omega$ is a compact $d$-set). But, if $\partial\Omega$ is unbounded, $\dim_{\rm A}(\partial\Omega)<n$ is not implied by Triebel's assumption that $\partial \Omega$ is a $d$-set with $n-1\leq d<n$. For,  as discussed at the end of Appendix \ref{sec:Appendix4}, if $d\in [0,n)$ there exists an unbounded $d$-set $F\subset \R^n$ with $\dimA(F)=n$, and, if $\Omega:= F^c$, then $\partial \Omega = F$.
 We recall from Remark \ref{rem:boundarydset} that Triebel's $d$-set assumption is not guaranteed in our main case of interest, when $\Omega$ is the interior or complement of an $n$-attractor, 
whereas the assumption that $\dim_{\rm A}(\partial\Omega)<n$ is guaranteed to hold in that case. 

The structure of our proof of {Theorem \ref{prop:wolff}} 
follows that of \cite[Theorem 4.17]{Tri08}. 
The key difference compared to {the proof
of \cite[Theorem 4.17]{Tri08}} 
is that, to obtain a common extension operator on an interval of smoothness parameters containing $s=0$, we appeal to Proposition \ref{prop:common-pfixed}, 
rather than to \cite[Theorem 4.10]{Tri08}, which requires the $d$-set assumption. 
{We also include more details than were provided in Triebel's proof, 
for the convenience of the non-expert reader.} 
In this theorem, as usual, $(\cdot,\cdot)_{\theta,q}$ denotes real interpolation and $[\cdot,\cdot]_\theta$ complex interpolation.
\begin{thm}
\label{prop:wolff} 
Let $\Omega\subset\R^n$ be a thick domain, with $\overline{\Omega}\neq \R^n$ and $\dim_{\rm A}(\partial\Omega)<n$.
Let $1<p<\infty$ and $0<\theta<1$.
\begin{enumerate}
\item Let $q_{0},q_{1},q\in(0,\infty)$. Let $s_1,s_2\in\R$ with $s_0\neq s_1$, and let 
\begin{equation}
s=(1-\theta)s_{0}+\theta s_{1}.\label{eq:s}
\end{equation}
Then
\begin{equation}
(B_{p,q_{0}}^{s_{0}}(\Omega),B_{p,q_{1}}^{s_{1}}(\Omega))_{\theta,q}=(F_{p,q_{0}}^{s_{0}}(\Omega),F_{p,q_{1}}^{s_{1}}(\Omega))_{\theta,q}=B_{p,q}^{s}(\Omega).\label{eq:4.43}
\end{equation}
\item Let $q_{0},q_{1}\in(1,\infty)$ 
and 
\[
\frac{1}{q}=\frac{1-\theta}{q_{0}}+\frac{\theta}{q_{1}}.
\]
Let $s_{0},s_{1}\in\mathbb{R}$ and let $s$ be given by (\ref{eq:s}). Then 
\begin{equation}
[B_{p,q_{0}}^{s_{0}}(\Omega),B_{p,q_{1}}^{s_{1}}(\Omega)]_{\theta}=B_{p,q}^{s}(\Omega)\label{eq:4.45}
\end{equation}
and 
\begin{equation}
[F_{p,q_{0}}^{s_{0}}(\Omega),F_{p,q_{1}}^{s_{1}}(\Omega)]_{\theta}=F_{p,q}^{s}(\Omega).\label{eq:4.46}
\end{equation}
\end{enumerate}
\end{thm}

\begin{proof}
We start with part 1 and divide the proof into several steps:

\textbf{Step 1.1.}

Since $B_{p,\min\{p,\kappa\}}^{t}(\Omega)\hookrightarrow F_{p,\kappa}^{t}(\Omega)\hookrightarrow B_{p,\max\{p,\kappa\}}^{t}(\Omega)$,
$t\in\R$, $0<\kappa<\infty$, --- check, e.g., \cite[Proposition 2.3.2/2.(iii), p. 47]{Triebel83ThFS}
and apply $|_{\Omega}$ --- it is sufficient to prove (\ref{eq:4.43})
for the $B$ spaces. Also, we assume that $s_0<s_1$, noting that once \eqref{eq:4.43} is proved for $s_0<s_1$ it follows also for $s_0>s_1$ by standard properties of interpolation spaces. 

\textbf{Step 1.2.}

By \cite[Theorem 2.4.2.(i)]{Triebel83ThFS}, one has
\[
(B_{p,q_{0}}^{s_{0}}\Rn ,B_{p,q_{1}}^{s_{1}}\Rn )_{\theta,q}=B_{p,q}^{s}\Rn 
\]
and notice that such an identity carries over to the corresponding spaces
on $\Omega$ for $s_{0},s_{1}$ in a given interval $(a,b)$ 
if, given $s_0,s_1\in(a,b)$, 
there exists an interval $(a',b')\subset(a,b)$ such that $s_{0},s_{1}\in(a',b')$
and a common linear and bounded extension operator
\begin{equation}
E:B_{p,\kappa}^{t}(\Omega)\to B_{p,\kappa}^{t}\Rn ,\quad t\in(a',b'),\;0<\kappa<\infty.\label{eq:E}
\end{equation}
The argument is as follows: it is then possible to extend the definition
of $E$ in an obvious way to a linear operator
\[
E:B_{p,q_{0}}^{s_{0}}(\Omega)+B_{p,q_{1}}^{s_{1}}(\Omega)\to B_{p,q_{0}}^{s_{0}}\Rn +B_{p,q_{1}}^{s_{1}}\Rn 
\]
such that $|_{\Omega}\circ E$ coincides with the identity operator
on $B_{p,q_{0}}^{s_{0}}(\Omega)+B_{p,q_{1}}^{s_{1}}(\Omega)$; using
the interpolation property \cite[\S1.6.2, Remark 2]{Triebel92ThFSII}, we get that $E$ maps $(B_{p,q_{0}}^{s_{0}}(\Omega),B_{p,q_{1}}^{s_{1}}(\Omega))_{\theta,q}$
linearly and boundedly into $(B_{p,q_{0}}^{s_{0}}\Rn ,B_{p,q_{1}}^{s_{1}}\Rn )_{\theta,q}$,
that is, into $B_{p,q}^{s}\Rn $. Thus $|_{\Omega}\circ E$
gives the inclusion $(B_{p,q_{0}}^{s_{0}}(\Omega),B_{p,q_{1}}^{s_{1}}(\Omega))_{\theta,q}$ $\hookrightarrow B_{p,q}^{s}(\Omega)$;
on the other hand, using the same interpolation property,
we get that $|_{\Omega}$ maps $(B_{p,q_{0}}^{s_{0}}\Rn ,$ $B_{p,q_{1}}^{s_{1}}\Rn )_{\theta,q}$
(that is, $B_{p,q}^{s}\Rn $) linearly and boundedly into $(B_{p,q_{0}}^{s_{0}}(\Omega),B_{p,q_{1}}^{s_{1}}(\Omega))_{\theta,q}$.
Thus $|_{\Omega}\circ E$ gives the inclusion $B_{p,q}^{s}(\Omega)\hookrightarrow(B_{p,q_{0}}^{s_{0}}(\Omega),B_{p,q_{1}}^{s_{1}}(\Omega))_{\theta,q}$.

There are three intervals $(a,b)$ to which the above argument can
be applied: the semi-infinite intervals $\R_{+} :=(0,\infty)$ and $\R_{-} :=(-\infty,0)$, and the finite interval $\left(\delta\left(\frac{1}{p}-1\right),\frac{\delta}{p}\right)$,
where $\delta\in (0,1]$ is as in Proposition \ref{prop:common-pfixed}.
In the case of $\R_{+}$, it follows from 
\cite[Theorem 4.4]{Tri08} with $(a',b')=(0,u)$,
where $u$ is chosen such that $s_{1}<u$, so that (\ref{eq:E}) holds
with $E$ the extension operator ${\rm ext}_{u}$ given in \cite[Theorem 4.4]{Tri08}.
In the case of $\R_{-}$, it follows from 
\cite[Theorem 4.7]{Tri08} with $(a',b')=(-\varepsilon^{-1},0)$,
where $0<\varepsilon<1$ is chosen such that $-\varepsilon^{-1}<s_{0}$,
so that (\ref{eq:E}) holds with $E$ the extension operator ${\rm ext}^{\varepsilon}$
given in \cite[Theorem 4.7]{Tri08}. Finally, in the third interval
pointed out above, it follows from Proposition \ref{prop:common-pfixed}
with $(a',b')=(a,b)$, so that (\ref{eq:E}) holds with $E$ the
extension operator ${\rm ext}$ from that proposition.

\textbf{Step 1.3.}

Up to now we have obtained the $B$ case of (\ref{eq:4.43}) for $s_{0},s_{1}$
in each of three intervals of real numbers, which have some overlap,
and this will allow us to extend the result to the whole of $\R$.
This is a consequence of reiteration and of the so-called Wolff's
interpolation theorem for real interpolation \cite[Theorem 1]{Wolff1982},
as we shall see.

We start by writing Wolff's
interpolation theorem in a form more suitable to
our context:

Let $-\infty<s_{0}<t_{2}<t_{3}<s_{1}<\infty$ and $0<\theta_{2}<\theta_{3}<1$
be such that
\begin{align}
\label{eq:theta23}
t_{2}=(1-\theta_{2})s_{0}+\theta_{2}s_{1}\quad\mbox{and}\quad t_{3}=(1-\theta_{3})s_{0}+\theta_{3}s_{1}.
\end{align}
Let $0<\eta_{2},\eta_{3}<1$ be such that
\begin{align}
\label{eq:eta23}
t_{2}=(1-\eta_{2})s_{0}+\eta_{2}t_{3}\quad\mbox{and}\quad t_{3}=(1-\eta_{3})t_{2}+\eta_{3}s_{1}.
\end{align}
Let $0<q_{2},q_{3}\leq\infty$. Assume that
\begin{equation}
(B_{p,q_{0}}^{s_{0}}(\Omega),B_{p,q_{3}}^{t_{3}}(\Omega))_{\eta_{2},q_{2}}=B_{p,q_{2}}^{t_{2}}(\Omega)\quad\mbox{and}\quad(B_{p,q_{2}}^{t_{2}}(\Omega),B_{p,q_{1}}^{s_{1}}(\Omega))_{\eta_{3},q_{3}}=B_{p,q_{3}}^{t_{3}}(\Omega).\label{eq:hypWolff}
\end{equation}
Then
\begin{equation}
(B_{p,q_{0}}^{s_{0}}(\Omega),B_{p,q_{1}}^{s_{1}}(\Omega))_{\theta_{2},q_{2}}=B_{p,q_{2}}^{t_{2}}(\Omega)\quad\mbox{and}\quad(B_{p,q_{0}}^{s_{0}}(\Omega),B_{p,q_{1}}^{s_{1}}(\Omega))_{\theta_{3},q_{3}}=B_{p,q_{3}}^{t_{3}}(\Omega).\label{eq:theWolff}
\end{equation}

To see that this is a particular reformulation of 
\cite[Theorem 1]{Wolff1982},  
notice that $B_{p,q_{0}}^{s_{0}}(\Omega)\cap B_{p,q_{1}}^{s_{1}}(\Omega)=B_{p,q_{1}}^{s_{1}}(\Omega)\subset B_{p,q_{3}}^{t_{3}}(\Omega)=B_{p,q_{2}}^{t_{2}}(\Omega)\cap B_{p,q_{3}}^{t_{3}}(\Omega)$
and, as follows by straightforward calculations\footnote{\label{fn:equivsystem}It may be helpful to notice first that (\ref{eq:etasthetas})
is equivalent to $\eta_{2}=\frac{\theta_{2}}{\theta_{3}}$ and $\eta_{3}=\frac{\theta_{3}-\theta_{2}}{1-\theta_{2}}$.},
\begin{equation}
\theta_{2}=\frac{\eta_{2}\eta_{3}}{1-\eta_{2}+\eta_{2}\eta_{3}}\quad\mbox{and}\quad\theta_{3}=\frac{\eta_{3}}{1-\eta_{2}+\eta_{2}\eta_{3}}.\label{eq:etasthetas}
\end{equation}

We claim now that whenever (\ref{eq:4.43}) holds in the $B$ case
for $s_{0},s_{1}$ in each of two intervals $(a,b)$ and $(a',b')$
such that $a<a'<b<b'$, then it also holds for $s_{0},s_{1}$ in the
union $(a,b')$. Of course, due to the assumption, what remains to
be shown is that it holds when $a<s_{0}\leq a'<b\leq s_{1}<b'$.

Consider $t_{2},t_{3}$ such that $a<s_{0}\leq a'<t_{2}<t_{3}<b\leq s_{1}<b'$
in the above formulation of Wolff's interpolation theorem. Since the
hypotheses in our claim imply that (\ref{eq:hypWolff}) holds with
$q_{2}=q_{3}=q$, then (\ref{eq:theWolff}) also holds under these
restrictions. In particular, this shows that (\ref{eq:4.43}) holds
in the $B$ case when $a<s_{0}\leq a'<s<b\leq s_{1}<b'$. In order to
show that it also holds in the remaining cases $a<s_{0}<s\leq a'<b\leq s_{1}<b'$
and $a<s_{0}\leq a'<b\leq s<s_{1}<b'$, just use reiteration. We illustrate
its use in the first of these two cases (but see Step 2.3 below, where
we illustrate its use for the second case when dealing with complex
interpolation): consider $t_{3}$ such that $a<s_{0}<s\leq a'<t_{3}<b\leq s_{1}<b'$
and $0<\eta_{2},\theta_{3}<1$ such that $s=(1-\eta_{2})s_{0}+\eta_{2}t_{3}$
and $t_{3}=(1-\theta_{3})s_{0}+\theta_{3}s_{1}$; from what we have
just obtained, we already have that $(B_{p,q_{0}}^{s_{0}}(\Omega),B_{p,q_{1}}^{s_{1}}(\Omega))_{\theta_{3},q}=B_{p,q}^{t_{3}}(\Omega)$;
on the other hand, since $s_{0},t_{3}\in(a,b)$, we have by hypothesis
that $(B_{p,q_{0}}^{s_{0}}(\Omega),B_{p,q}^{t_{3}}(\Omega))_{\eta_{2},q}=B_{p,q}^{s}(\Omega)$;
therefore, using reiteration, $B_{p,q}^{s}(\Omega)=(B_{p,q_{0}}^{s_{0}}(\Omega),(B_{p,q_{0}}^{s_{0}}(\Omega),B_{p,q_{1}}^{s_{1}}(\Omega))_{\theta_{3},q})_{\eta_{2},q}=(B_{p,q_{0}}^{s_{0}}(\Omega),B_{p,q_{1}}^{s_{1}}(\Omega))_{\theta_{3}\eta_{2},q}$,
where, as is easily seen\footnote{Note that we can use footnote \ref{fn:equivsystem} here too, since
$s$ and $\theta$ here are taking the roles, respectively, of $t_{2}$
and $\theta_{2}$ there.}, $\theta_{3}\eta_{2}=\theta$. Our claim is proved.

We use now this claim first with $(a,b)=\R_{-}$ and $(a',b')=\left(\delta\left(\frac{1}{p}-1\right),\frac{\delta}{p}\right)$
and afterwards with $(a,b)=\left(
{-\infty}
,\frac{\delta}{p}\right)$
and $(a',b')=\R_{+}$ to conclude the proof of part 1.

\medskip{}

We will be briefer in the proof of part 2, focusing more in the modifications
to be made to the arguments used previously.

\textbf{Step 2.1.}

We resort again to the use of $A$ to mean $B$ throughout or to mean
$F$ throughout in what follows. And when $s_{0}\not=s_{1}$ we shall
assume that $s_{0}<s_{1}$, since once (\ref{eq:4.45}) and
(\ref{eq:4.46}) are proved for $s_{0}<s_{1}$ they hold also
for $s_{0}>s_{1}$, by properties of interpolation spaces.

\textbf{Step 2.2.}

By \cite[Theorems 2.4.1.(d) and 2.4.2/1.(d), pp. 182, 185]{Triebel78ITFSDO},
one has
\[
[A_{p,q_{0}}^{s_{0}}\Rn ,A_{p,q_{1}}^{s_{1}}\Rn ]_{\theta}=A_{p,q}^{s}\Rn 
\]
and, arguing as in the proof of Step 1.2, such an identity carries over to the corresponding spaces
on $\Omega$ for $s_{0},s_{1}$ in a given interval $(a,b)$ if, given $s_{0},s_{1}\in (a,b)$, 
there exists an interval $(a',b')\subset(a,b)$ such that $s_{0},s_{1}\in(a',b')$
and a common linear and bounded extension operator
\begin{equation}
E:A_{p,\kappa}^{t}(\Omega)\to A_{p,\kappa}^{t}\Rn ,\quad t\in(a',b'),\;1<\kappa<\infty.\label{eq:E-1}
\end{equation}
Again arguing as in the proof of Step 1.2, 
this holds when $(a,b)$ is
one of the three intervals $\R_{+}$, $\R_{-}$ and $\left(\delta\left(\frac{1}{p}-1\right),\frac{\delta}{p}\right)$,
where $\delta\in (0,1]$ is as in Proposition \ref{prop:common-pfixed}.

\textbf{Step 2.3.}

As in Step 1.3, the fact that (\ref{eq:4.45}) and (\ref{eq:4.46})
hold when $s_{0},s_{1}$ belong to each of the above three intervals
of real numbers, which have some overlap, will allow us to extend
the result to the whole of $\R$. This is obvious when $s_{0}=s_{1}$
and finishes the proof of part 2 in this case. When $s_{0}<s_{1}$
--- which we assume hereafter --- it is again a consequence of reiteration
and of a so-called Wolff's interpolation theorem, now with
respect to complex interpolation \cite[Theorem 2]{Wolff1982}, as
we shall see.

We start by writing the relevant Wolff's interpolation theorem in a form more suitable to
our context:

Let $-\infty<s_{0}<t_{2}<t_{3}<s_{1}<\infty$, $0<\theta_{2}<\theta_{3}<1$ and $0<\eta_{2},\eta_{3}<1$ 
be such that \eqref{eq:theta23} and \eqref{eq:eta23} hold. 
Let 
\begin{equation}
\frac{1}{q_{2}}=\frac{(1-\theta_{2})}{q_{0}}+\frac{\theta_{2}}{q_{1}}\quad\mbox{and}\quad\frac{1}{q_{3}}=\frac{(1-\theta_{3})}{q_{0}}+\frac{\theta_{3}}{q_{1}}.\label{eq:q2q3complex}
\end{equation}
Assume that
\begin{equation}
[A_{p,q_{0}}^{s_{0}}(\Omega),A_{p,q_{3}}^{t_{3}}(\Omega)]_{\eta_{2}}=A_{p,q_{2}}^{t_{2}}(\Omega)\quad\mbox{and}\quad[A_{p,q_{2}}^{t_{2}}(\Omega),A_{p,q_{1}}^{s_{1}}(\Omega)]_{\eta_{3}}=A_{p,q_{3}}^{t_{3}}(\Omega).\label{eq:hypWolff-1}
\end{equation}
Then
\begin{equation}
[A_{p,q_{0}}^{s_{0}}(\Omega),A_{p,q_{1}}^{s_{1}}(\Omega)]_{\theta_{2}}=A_{p,q_{2}}^{t_{2}}(\Omega)\quad\mbox{and}\quad[A_{p,q_{0}}^{s_{0}}(\Omega),A_{p,q_{1}}^{s_{1}}(\Omega)]_{\theta_{3}}=A_{p,q_{3}}^{t_{3}}(\Omega).\label{eq:theWolff-1}
\end{equation}

To see that this is a particular reformulation of 
\cite[Theorem 2]{Wolff1982}, 
notice that $A_{p,q_{0}}^{s_{0}}(\Omega)\cap A_{p,q_{1}}^{s_{1}}(\Omega)$
is a dense subspace of both $A_{p,q_{2}}^{t_{2}}(\Omega)$ and $A_{p,q_{3}}^{t_{3}}(\Omega)$
and that (\ref{eq:etasthetas}) (as well as footnote \ref{fn:equivsystem})
holds. Notice also that it comes as a consequence of the setting above
that 
\begin{equation}
\frac{1}{q_{2}}=\frac{(1-\eta_{2})}{q_{0}}+\frac{\eta_{2}}{q_{3}}\quad\mbox{and}\quad\frac{1}{q_{3}}=\frac{(1-\eta_{3})}{q_{2}}+\frac{\eta_{3}}{q_{1}}.\label{eq:compatibility}
\end{equation}

Now a claim corresponding to the one made in Step 1.3 above can be
proved with respect to (\ref{eq:4.45}) and (\ref{eq:4.46}) instead
of (\ref{eq:4.43}) and afterwards the proof of part 2 finishes as
at the end of Step 1.3, now using this new claim and what has been
obtained in Step 2.2 above. For the convenience of the reader, we
spell out some differences in the proof of this new claim, in comparison
with the proof of the previous one in Step 1.3 (of course, now we are referring to the complex version of Wolff's
interpolation theorem and reiteration is being applied in the context of complex interpolation).

To prove the conclusion of the claim when $a<s_{0}\leq a'<s<b\leq s_{1}<b'$,
we consider the same $t_{2},t_{3}$ such that $a<s_{0}\leq a'<t_{2}<t_{3}<b\leq s_{1}<b'$
but now consider $q_{2}$ and $q_{3}$ as in (\ref{eq:q2q3complex}).
Due to (\ref{eq:compatibility}), the hypotheses in our claim imply
that (\ref{eq:hypWolff-1}) holds, so that (\ref{eq:theWolff-1}) also
holds under these restrictions. In order to deal with the remaining
cases $a<s_{0}<s\leq a'<b\leq s_{1}<b'$ and $a<s_{0}\leq a'<b\leq s<s_{1}<b'$,
again just use reiteration. We illustrate its use now in the second
of these two cases (while in Step 1.3 above we have illustrated its
use for the first case when dealing with real interpolation): consider
$t_{2}$ such that $a<s_{0}\leq a'<t_{2}<b\leq s<s_{1}<b'$ and $0<\eta_{3},\theta_{2}<1$
such that $t_{2}=(1-\theta_{2})s_{0}+\theta_{2}s_{1}$ and $s=(1-\eta_{3})t_{2}+\eta_{3}s_{1}$;
consider also $q_{2},q_{3}$ defined by $\frac{1}{q_{2}}=\frac{(1-\theta_{2})}{q_{0}}+\frac{\theta_{2}}{q_{1}}$
and $\frac{1}{q_{3}}=\frac{(1-\eta_{3})}{q_{2}}+\frac{\eta_{3}}{q_{1}}$;
from what we have just obtained, we already have that $[A_{p,q_{0}}^{s_{0}}(\Omega),A_{p,q_{1}}^{s_{1}}(\Omega)]_{\theta_{2}}=A_{p,q_{2}}^{t_{2}}(\Omega)$;
on the other hand, since $t_{2},s_{1}\in(a',b')$, we have by hypothesis
that $[A_{p,q_{2}}^{t_{2}}(\Omega),A_{p,q_{1}}^{s_{1}}(\Omega)]_{\eta_{3}}=A_{p,q_{3}}^{s}(\Omega)$;
however, we can make $s$ fit in the definitions of parameters in
the above Wolff's interpolation theorem, namely with $s$ taking the
role of $t_{3}$ over there, and consequently with $q$ coinciding
with the $q_{3}$ in (\ref{eq:q2q3complex}), which is the same as
our $q_{3}$ here, due to (\ref{eq:compatibility}); therefore, using
reiteration, $A_{p,q}^{s}(\Omega)=[[A_{p,q_{0}}^{s_{0}}(\Omega),A_{p,q_{1}}^{s_{1}}(\Omega)]_{\theta_{2}},A_{p,q_{1}}^{s_{1}}(\Omega),]_{\eta_{3}}=[A_{p,q_{0}}^{s_{0}}(\Omega),A_{p,q_{1}}^{s_{1}}(\Omega)]_{\theta_{2}(1-\eta_{3})+\eta_{3}}$,
where, as is easily seen\footnote{Notice that, again, we can use footnote \ref{fn:equivsystem}, since
$s$ and $\theta$ here are taking the roles, respectively, of $t_{3}$
and $\theta_{3}$ there.}, $\theta_{2}(1-\eta_{3})+\eta_{3}=\theta$.
\end{proof}

Theorem \ref{prop:wolff} implies that, under appropriate assumptions, the $F^s_{p,q}(\Omega)$ and $B^s_{p,q}(\Omega)$ spaces form interpolation scales, {in the following sense}. 

{
\begin{defn}
\label{d:IntScale}
Given an interval $I \subset \mathbb{R}$ and a family of quasi-Banach spaces $\{A_s\}_{s \in I}$, we say that $\{A_s\}_{s \in I}$ is an \emph{interpolation scale} with respect to a given interpolation method if, for all $s_0, s_1 \in I$ and all $\theta \in (0,1)$, interpolating the pair $(A_{s_0}, A_{s_1})$ with parameter $\theta$ yields $A_{(1-\theta)s_0 + \theta s_1}$, with equivalent quasi-norms.
\end{defn}}
\begin{cor}
\label{cor:InterpolationScale}
Let $\Omega$ be as in Theorem \ref{prop:wolff}.
Let $1<p<\infty$.
\begin{enumerate}
	\item[(i)] Let $q\in(0,\infty).$ Then $\{B_{p,q}^{s}(\Omega)\}_{s\in\R}$ is
	an interpolation scale with respect to real interpolation $(\cdot,\cdot)_{\theta,q}$.
	\item[(ii)] Let $q\in(1,\infty)$. Then $\{B_{p,q}^{s}(\Omega)\}_{s\in\R}$ and
	$\{F_{p,q}^{s}(\Omega)\}_{s\in\R}$ are both interpolation scales
	with respect to complex interpolation $[\cdot,\cdot]_{\theta}$.
\end{enumerate}
\end{cor}

\begin{proof}
We start with part (i), where we have to show that, given any $s_{0},s_{1}\in\R$
and $0<\theta<1$,
\[
(B_{p,q}^{s_{0}}(\Omega),B_{p,q}^{s_{1}}(\Omega))_{\theta,q}=B_{p,q}^{(1-\theta)s_{0}+\theta s_{1}}(\Omega).
\]
If $s_{0}=s_{1}$, it follows from general properties of real interpolation.
In the opposite case, again by general properties of real interpolation
it is enough to deal with the case when $s_{0}<s_{1}$. Finally, the
result for this case is an immediate consequence of part 1 of Theorem
\ref{prop:wolff}.

As for part (ii), i.e., that, given $s_{0},s_{1}\in\R$ and $0<\theta<1$,
\[
[B_{p,q}^{s_{0}}(\Omega),B_{p,q}^{s_{1}}(\Omega)]_{\theta}=B_{p,q}^{(1-\theta)s_{0}+\theta s_{1}}(\Omega)\quad\mbox{and}\quad[F_{p,q}^{s_{0}}(\Omega),F_{p,q}^{s_{1}}(\Omega)]_{\theta}=F_{p,q}^{(1-\theta)s_{0}+\theta s_{1}}(\Omega),
\]
this is an immediate consequence of part 2 of Theorem \ref{prop:wolff}.
\end{proof}

From Corollary \ref{cor:InterpolationScale} one can derive interpolation results about the spaces $\tA_{p,q}^{s}(\Omega)$ by duality.

\begin{cor}
\label{cor:InterpolationScaleTilde}
{Let $\Omega$ be as in Theorem \ref{prop:wolff}.}
Let $1<p,q<\infty$. Then $\{\tB_{p,q}^{s}(\Omega)\}_{s\in\R}$ is an interpolation scale with respect to real interpolation $(\cdot,\cdot)_{\theta,q}$ and $\{\tB_{p,q}^{s}(\Omega)\}_{s\in\R}$ and $\{\tF_{p,q}^{s}(\Omega)\}_{s\in\R}$ are both interpolation scales with respect to complex interpolation $[\cdot,\cdot]_{\theta}$. 
\end{cor}
\begin{proof}
This follows from Corollary \ref{cor:InterpolationScale}, standard duality properties of interpolation spaces (e.g.\ \cite[Theorems 1.11.2 and 1.11.3]{Triebel78ITFSDO}) and the fact that, as mentioned in \S\ref{sec:FunctionSpaceNotation},
{$\tA^{-s}_{p',q'}(\Omega)$ is a realisation of the dual space of $A^s_{p,q}(\Omega)$} for $s\in \R$ and $1<p,q<\infty$. 
\end{proof}

For convenience we state separately the results for the {$p=2$} Sobolev space case. We recall that in this Hilbert space setting the real and complex interpolation methods produce the same space (see, e.g., \cite{InterpolationCWHM,InterpolationE_CWHM}).
\begin{cor}
\label{cor:InterpolationSobolev}
{Let $\Omega$ be as in Theorem \ref{prop:wolff}.} 
Then 
$\{H^{s}(\Omega)\}_{s\in\R}$ and $\{\tH^{s}(\Omega)\}_{s\in\R}$ are interpolation scales with respect to both real interpolation $(\cdot,\cdot)_{\theta,2}$ and complex interpolation $[\cdot,\cdot]_{\theta}$. 
\end{cor}

{Finally, we specialise to the $n$-attractor case.
\begin{cor}
\label{cor:Interpolationnset}
Let $\Gamma\subset\R^n$ be an $n$-attractor and 
let either $\Omega=\Gamma^\circ$ or $\Omega=\Gamma^c$. 
Then the conclusions of Theorem \ref{prop:wolff} and Corollaries \ref{cor:InterpolationScale},  \ref{cor:InterpolationScaleTilde} and \ref{cor:InterpolationSobolev} apply. 
\end{cor}
}
\begin{proof}
{It suffices to show that the conditions of Theorem \ref{prop:wolff} are satisfied. For this, the thickness of $\Omega$ holds by Theorem \ref{thm:thick}, the fact that $\overline{\Omega}\neq \R^n$ holds because $\Gamma$ is compact and $\Gamma^\circ$ is non-empty, and the fact that $\dim_{\rm A}(\partial\Omega)<n$ holds by Proposition \ref{prop:nset}\eqref{v},\eqref{vi} (see the proof of Corollary \ref{cor:newDt}). 
}
\end{proof}

\begin{rem}
As mentioned before Theorem \ref{prop:wolff}, 
under the additional assumption that $\partial\Omega$ is a $d$-set for some $n-1\leq d< n$, the conclusions of 
Theorem \ref{prop:wolff}, and hence also Corollaries \ref{cor:InterpolationScale}-\ref{cor:Interpolationnset},
follow alternatively by \cite[Theorem~4.17]{Tri08}.
Our results 
therefore extend those provided by \cite[Theorem~4.17]{Tri08} to cases where $\dim_{\rm A}(\partial\Omega)<n$ but $\partial\Omega$ is not a $d$-set for any $n-1\leq d<n$ (cf.\ Remark \ref{rem:boundarydset}). 

We remark that the results of \cite[Theorem~4.17]{Tri08} apply in particular to the class of classical snowflake domains $\Omega$ considered in \cite[\S5.1]{caetano2019density}, since for these domains thickness of $\Omega$ is provided by \cite[Proposition~5.2]{caetano2019density}, and 
{the required $d$-set property of $\partial\Omega$ was proved in \cite[Proposition~5.3 and Remark~5.4]{caetano2019density}.}
 To the best of our knowledge only one of these domains is an $n$-attractor --- the Koch snowflake (see Figure \ref{fig:2sets}). 
\end{rem}

\begin{rem}
To compare the results of Corollary \ref{cor:InterpolationSobolev} to other cases, we note that if 
$\Omega\subset\R^n$ is a Lipschitz hypograph or a Lipschitz domain then $\{H^{s}(\Omega)\}_{s\in\R}$ and $\{\tH^{s}(\Omega)\}_{s\in\R}$ are known to be interpolation scales 
with respect to both real and complex interpolation, 
and, more generally, if $\Omega$ is an $(\epsilon,\delta)$ locally uniform domain, for some $\epsilon, \delta>0$, then this same conclusion holds for the reduced range of spaces $\{H^{s}(\Omega)\}_{s\geq 0}$ and $\{\tH^{s}(\Omega)\}_{s\leq 0}$ (see, e.g., \cite[Corollaries 4.7 and 4.10]{InterpolationCWHM}). 
Still more generally (cf.\ \cite[Proposition~1, p119]{wallin1991trace}), if $\Omega$ is an $n$-set, the same conclusion holds for the slightly more reduced range of spaces $\{H^{s}(\Omega)\}_{s>0}$ and $\{\tH^{s}(\Omega)\}_{s<0}$: using the duality and common extension arguments used before (see, in particular, the arguments in Step 1.2 of the proof of Theorem \ref{prop:wolff} above), the assertion follows from \cite[Proposition 2.13]{bechtel:19}, with a reference to \cite{Rychkov:00}. By similar arguments, it also follows for that more restricted range of spaces if, alternatively, $\Omega$ is an I-thick domain with $\overline{\Omega} \not= \mathbb{R}^n$ and $|\partial \Omega|=0$ (cf.\ \cite[Theorem 4.4]{Tri08}).
\end{rem}

\subsection{Piecewise-constant approximation results}
\label{sec:Approximation}

In this section we apply the interpolation results obtained in the previous section to prove best approximation error estimates in fractional Sobolev spaces (in the case $p=2$) for the approximation of functions on $n$-attractors by piecewise-constant functions. 
We first prove a result about 
meshes on general domains, then specialise this to the case of $n$-attractors. 


\begin{defn}[Mesh]
\label{def:Mesh} 
{
Let $\Omega\subset\R^n$ be a domain.
Given some countable (e.g.\ finite) index set $J$ 
we say that $\cT=\{T_j\}_{j\in J}$ is a mesh of $\Omega$ 
if, for each $j\in J$, $T_j\subset\Omega$ is a bounded domain, 
$T_j\cap T_{j'}= \emptyset$ for $j'\neq j$, 
and 
\[ \Big| \Omega \setminus \bigcup_{j\in J} T_j\Big|=0.\]
}
\end{defn}

Given $h>0$ and a mesh $\cT_h=\{T_j\}_{j\in J}$ of $\Omega$ such that $\diam(T_j)\leq h$ for $j\in J$, let $V_h$ denote the set of piecewise constant functions on $\cT_h$ that lie in $L_2(\Omega)$, and let $P_h:L_2(\Omega)\to V_h$ denote the $L_2$-orthogonal projection onto $V_h$. We can then prove the following approximation result, in which the constant $C$ depends only on the dimension $n$\footnote{It is possible, using similar methods of argument, to prove versions of Proposition \ref{prop:Approx} and Corollary \ref{cor:Interpolation} for piecewise polynomial approximation of arbitrary degree $m\geq 0$, extending the results here for the piecewise constant case $m=0$; see \cite{PiecewisePolynomial} for details.}. 
{We recall that in our notation $H^1(\Omega)=\{u|_\Omega:u\in H^1(\R^n)\}$ so the existence of the $U$ in the statement of the following proposition is guaranteed without any regularity assumption on $\Omega$.} 

\begin{prop}
\label{prop:Approx}
Let $\Omega\subset\R^n$ be a domain.
Let $h>0$ 
and suppose that $\cT_h=\{\Omega_j\}_{j\in J}$ is a 
mesh of $\Omega$
such that $\diam(\Omega_j)\leq h$ for $j\in J$.
Then, for any $u\in H^1(\Omega)$, and any $U\in H^1(\R^n)$ such that $U|_\Omega = u$, 
\begin{align}
\label{eq:PoincareConsequenceGeneral1}
\|u-P_hu\|_{L_2(\Omega)} \leq Ch \|\nabla U\|_{L_2(\R^n)},
\end{align} 
and 
\begin{align}
\label{eq:PoincareConsequenceGeneral2}
\|u-P_hu\|_{L_2(\Omega)} \leq Ch\|u\|_{H^1(\Omega)},
\end{align} 
where
\[C=\frac{3^{1+n/2}\sqrt{n}}{\pi}.\]
\end{prop}

\begin{proof}
{The main tool in the proof is the Poincar\'e inequality. However, rather than applying the Poincar\'e inequality on the mesh elements themselves, we instead apply the Poincar\'e inequality on an overlapping grid of cubes based on a tessellation of $\R^n$, which allows us to obtain an upper bound for the summed contributions to the best approximation error from the mesh elements.}

\begin{figure}
\centering
\includegraphics[width=60mm]{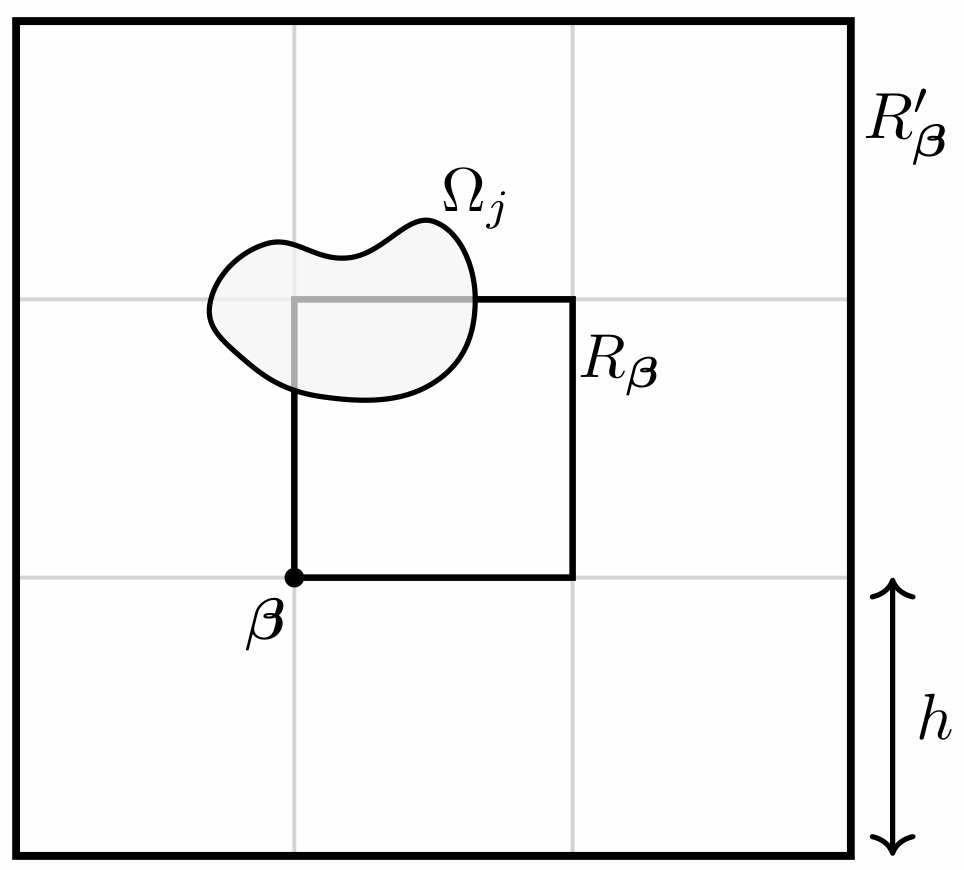}
\caption{{Covering the mesh $\cT_h$ by overlapping cubes $R_{\bbeta}'$. 
If $\Omega_j\in \cT_h$ intersects the cube $R_{\bbeta}$, for some $\bbeta \in \Z^n$, then $\Omega_j$ is contained in the larger cube $R_{\bbeta}'$.
}}
\label{f:Cubes}
\end{figure}

For a non-empty bounded measurable set $E\subset \R^n$ and $v\in L_1(E)$ let $\pi_E v$ denote the average of $v$ over $E$, i.e.\ $\pi_Ev = \frac{1}{|E|}\int_E v$, which for $v\in L_2(E)$ represents the best $L_2$-approximation of $v$ on $E$ by a constant.  
We recall the standard fact that on any open cube $Q$ of side-length $r$ 
we have the Poincar\'e inequality (e.g., 
\cite{PayneWeinberger60,Bebendorf03,NazRep14}, \cite[\S5.8.1]{Evans2010}) 
\begin{equation} \label{eq:PoincareBall}
\left\|u-\pi_Q u\right\|_{L_2(Q)} \leq \frac{\sqrt{n}r}{\pi}\|\nabla u\|_{L_2(Q)}, \quad u\in W^1(Q),
\end{equation}
where $W^1(Q)=\{u\in L_2(Q): \nabla u\in L_2(Q)\}$.  
Now let 
$h$, $\Omega$, $\cT_h$, $u$ and $U$ be as in the hypothesis of the {proposition}. 
For each $\bbeta\in \Z^n$ let $R_{\bbeta}=[0,h]^n + h\bbeta$. 
If $j\in J$  
and $\Omega_j\cap R_{\bbeta}\neq \emptyset$ then, since $\diam(\Omega_j)\leq h$, $\Omega_j$ is contained in $R_{\bbeta}'$, the open cube of side length $3h$ formed by taking the interior of the union of  $R_{\bbeta}$ with all its neighbouring cubes 
{--- see Figure \ref{f:Cubes} for an illustration.} 
Hence, by \eqref{eq:PoincareBall} with $r=3h$, 
\begin{align*}
\label{}
\sum_{\substack{j\in J\\ \Omega_j\cap R_\bbeta \neq \emptyset}}\|u-\pi_{\Omega_j}u\|^2_{L_2(\Omega_j)}
\leq \sum_{\substack{j\in J\\ \Omega_j\cap R_\bbeta \neq \emptyset}}\|u-\pi_{R_{\bbeta}'}U\|^2_{L_2(\Omega_j)}
& = \sum_{\substack{j\in J\\ \Omega_j\cap R_\bbeta \neq \emptyset}}\|U-\pi_{R_{\bbeta}'}U\|^2_{L_2(\Omega_j)}\\
& \leq \|U-\pi_{R_{\bbeta}'}U\|^2_{L_2(R_{\bbeta}')}\\
&\leq \frac{9nh^2}{\pi^2}\|\nabla U\|^2_{L_2(R_{\bbeta}')}.
\end{align*}
Then 
\begin{align*}
 \label{}
\|u - P_hu\|_{L_2(\Omega)}^2
 = \sum_{j\in J}\|u-\pi_{\Omega_j}u\|^2_{L_2(\Omega_j)}
 &\leq \sum_{\bbeta\in \Z^n}
\sum_{\substack{j\in J\\ \Omega_j\cap R_\bbeta \neq \emptyset}}\|u-\pi_{\Omega_j}u\|^2_{L_2(\Omega_j)}\\
&\leq 
\frac{9nh^2}{\pi^2}
\sum_{\bbeta\in \Z^n}
\|\nabla U\|^2_{L_2(R_{\bbeta}')}\\
&= \frac{3^{n+2}nh^2}{\pi^2} \sum_{\bbeta\in \Z^n} \|\nabla U\|^2_{L_2(R_{\bbeta})},
\end{align*}
where we used the fact that the contribution of each cube $R_\bbeta$ gets counted $3^n$ times in the penultimate sum, 
{because each cube $R_\bbeta$ is contained in exactly $3^n$ of the enlarged cubes (see Figure \ref{f:Cubes}).} 
From this,  \eqref{eq:PoincareConsequenceGeneral1} follows by taking square roots. Clearly \eqref{eq:PoincareConsequenceGeneral2} is a consequence of \eqref{eq:PoincareConsequenceGeneral1}.
\end{proof}

The following corollary follows from Proposition \ref{prop:Approx} by interpolation theory and a standard duality argument (cf.\ \cite[Theorem~4.1.33]{sauter-schwab11}, \cite[Lemma A.1]{BEMfract}). 
Here, and henceforth, $\widetilde{u}$, for $u\in L_2(\Omega)$, denotes the extension of $u$ by zero from $\Omega$ to $\R^n$ and, for $h>0$, $\widetilde V_h := \{\widetilde u:u\in V_h\}$.
 We recall that, for any domain $\Omega\subset\R^n$, extension by zero is a unitary isomorphism from $H^0(\Omega)=L_2(\Omega)$ onto $\tH^0(\Omega)$ (its inverse being the restriction operator). 
We also recall from \cite[Corollary~4.9]{InterpolationCWHM} (cf.\ also the proof of Corollary \ref{cor:InterpolationScaleTilde} above) that $\{H^s(\Omega)\}_{0\leq s\leq 1}$ is an interpolation scale if and only if $\{\tH^s(\Omega)\}_{-1\leq s\leq 0}$ is an interpolation scale. 
\begin{cor}
\label{cor:Interpolation}
Under the assumptions of Proposition \ref{prop:Approx}, and with $C$ denoting the constant from that proposition, we have that, 
for $-1\leq s_1\leq 0 \leq s_2\leq 1$,
\begin{align}
\label{eq:Hsboundh_dn}
\inf_{u_h\in \widetilde{V}_h}
\|\widetilde{u}-u_h\|_{\tH^{s_1}(\Omega)}\leq
\|\widetilde{u}-\widetilde{P_hu}\|_{\tH^{s_1}(\Omega)}\leq (Ch)^{s_2-s_1}\|u\|_{H^{s_2}(\Omega)},\qquad u\in H^{s_2}(\Omega).
\end{align}
Suppose that additionally $\{H^s(\Omega)\}_{0\leq s\leq 1}$ is an interpolation scale (equivalently, that $\{\tH^s(\Omega)\}_{-1\leq s\leq 0}$ is an interpolation scale). 
Then, for $-1\leq s_1\leq s_2\leq 0$, 
\begin{align}
\label{eq:Hsboundh_dn2}
\inf_{u_h\in \widetilde{V}_h}
\|u-u_h\|_{\tH^{s_1}(\Omega)}\leq c h^{s_2-s_1}\|u\|_{\tH^{s_2}(\Omega)},\qquad u\in \tH^{s_2}(\Omega),
\end{align}
for some $c>0$ independent of $u$ and $h$. 
\end{cor}
\begin{proof}
For the first part we adopt the argument used in the proof of \cite[Lemma~A.1]{ScreenBEM}. Proposition \ref{prop:Approx} provides the bound \eqref{eq:Hsboundh_dn} for the case $s_1=0$ and $s_2=1$. And for the case $s_1=s_2=0$ \eqref{eq:Hsboundh_dn} follows from the fact that $P_h$ is an orthogonal projection in $L_2(\Omega)$.  
From these two cases we can deduce \eqref{eq:Hsboundh_dn} for the case $s_1=0$, $0\leq s_2\leq 1$ by interpolation applied to the composition of extension by zero with the operator $I-P_h$. For details we refer to the argument in \cite[Eqn~(46)]{ScreenBEM} and the surrounding sentences, noting that this argument does not require $\{H^s(\Omega)\}_{0\leq s\leq 1}$ to be an interpolation scale, using only that 
$H^{s_2}(\Omega)$ is embedded in $(L_2(\Omega),H^1(\Omega))_{s_2,2}$ and that the norm of the embedding is $\leq 1$ (as is guaranteed by \cite[Lemma~4.2]{InterpolationCWHM}).
Finally, one can then deduce \eqref{eq:Hsboundh_dn} for $-1\leq s_1\leq 0 \leq s_2\leq 1$ using the standard duality argument at the end of the proof of \cite[Lemma~A.1]{ScreenBEM}.

For the second part, under the additional assumption that $\{H^s(\Omega)\}_{0\leq s\leq 1}$ is an interpolation scale (equivalently, that $\{\tH^s(\Omega)\}_{-1\leq s\leq 0}$ is an interpolation scale), 
\eqref{eq:Hsboundh_dn2} follows for $-1\leq s_1\leq s_2\leq 0$ by applying interpolation to $I-P^{s_1}_h$, where $P^{s_1}_h$ denotes orthogonal projection in $\tH^{s_1}(\Omega)$ onto $\widetilde{V}_h$, using \eqref{eq:Hsboundh_dn} with $s_2=0$, combined with the trivial estimate
\begin{align}
\label{eq:Hs1bound}
\inf_{u_h\in \widetilde{V}_h}
\|u-u_h\|_{\tH^{s_1}(\Omega)}  = 
\|u-P^{s_1}_h u\|_{\tH^{s_1}(\Omega)}   
\leq  
\|u\|_{\tH^{s_1}(\Omega)}, \qquad u\in\tH^{s_1}(\Omega).
\end{align}

\end{proof}

Now let $\Gamma$ be an $n$-attractor. 
Then for each $0<h\leq \diam(\Gamma)$ the index set $L_h$ introduced in \eqref{eq:LhDef} defines a mesh of $\Omega:=\Gamma^\circ$, based on the self-similar structure of $\Gamma$, that satisfies the assumptions of Proposition \ref{prop:Approx}.
Furthermore, the interpolation scale assumption of Corollary \ref{cor:Interpolation} holds in this case by the results in \S\ref{sec:Interpolation}. Hence we obtain the following result, in which the notation $\Omega_\bm$ is defined as in \S\ref{sec:nsets}, i.e.\ for a vector index $\bm\in L_h$ we define $\Omega_\bm:=s_\bm(\Omega)=(\Gamma_\bm)^\circ$. 
\begin{cor}
\label{cor:approx}
Let $\Gamma\subset\R^n$ be an $n$-attractor {and} $\Omega:=\Gamma^\circ$.
Given $0<h\leq {\diam(\Gamma)}$, let $L_h$ be as defined in \eqref{eq:LhDef}, 
and let $\widetilde{V}_h$ be the space of piecewise constant functions on the mesh $\cT_h:=\{\Omega_\bm\}_{\bm\in L_h}$ of $\Omega$, 
extended by zero to functions on $\R^n$.
Then \eqref{eq:Hsboundh_dn} holds for $-1\leq s_1\leq 0 \leq s_2\leq 1$, and \eqref{eq:Hsboundh_dn2} holds for $-1\leq s_1\leq s_2\leq 0$, for some $c>0$ independent of $u$ and $h$. 
\end{cor}
\begin{proof}
We first prove that $\{\Omega_\bm\}_{\bm\in L_h}$ is a 
mesh of $\Omega$ in the sense of Definition \ref{def:Mesh}. For $\bm\in L_h$, $\Omega_\bm$ is by definition a bounded domain with $\diam(\Omega_\bm)\leq h$. Also, that 
{$\Omega_\bm\subset\Omega$} and 
$\Omega_\bm \cap \Omega_{\bm'}= \emptyset$ for $\bm\neq{\bm'}$ holds by Proposition \ref{prop:nset}\eqref{viii}. Finally, the fact that 
{$|\Omega\setminus \bigcup_{\bm\in L_h(\Omega)}\Omega_{\bm}|=0$} 
follows from 
{Proposition \ref{prop:nset}\eqref{iv} and \eqref{vi}} and the results in \cite[Section 5.1.1]{joly2024high}.
The statement of the corollary then follows from Corollary \ref{cor:Interpolation}, using the fact that the interpolation scale assumption is satisfied (see Corollary \ref{cor:InterpolationSobolev}). 
\end{proof}

\section{Application to screen scattering problems}
\label{sec:BVPsBIEs}


Our main application in this paper of the results in \S\ref{sec:nsets}-\S\ref{sec:FunctionSpaces} is to the problem of time-harmonic acoustic scattering of an incident wave $u^i$ propagating in $\R^{n+1}$ ($n=1,2$) by a sound-soft planar screen $\Gamma$. {Following \cite{ScreenPaper} we assume that $\Gamma$ is a 
 bounded} subset
of the hyperplane $\Gamma_\infty=\R^{n}\times\{0\}$; the propagation domain is then $D:=\R^{n+1}\setminus{\overline{\Gamma}}$. 
Our particular focus is on the case where $\Gamma=\tGamma\times\{0\}$ for some $n$-attractor $\tGamma\subset\R^n$ {(in which case $\Gamma=\overline{\Gamma}$)}. For brevity we shall describe this case simply by saying that ``$\Gamma\subset\Gamma_\infty$ is an $n$-attractor''. Using the results of \S\ref{sec:nsets}-\S\ref{sec:FunctionSpaces} we shall show that the classical integral equation formulation of the scattering problem is well-posed, and provide error estimates for its numerical approximation using a piecewise-constant Galerkin method. 

Before stating the problem, we fix some notation. As in \cite{CoercScreen2,ScreenPaper,BEMfract,HausdorffBEM} we define function spaces on the hyperplane $\Gamma_\infty := \R^{n}\times\{0\}$ and subsets of it (for example, the {closure $\overline{\Gamma}\subset\Gamma_\infty$ and relative interior $\Gamma^\circ\subset \Gamma_\infty$ of the screen $\Gamma$}) by associating $\Gamma_\infty$ with $\R^{n}$ and $\Gamma$ with the set $\tilde{\Gamma}\subset\R^n$ such that $\Gamma=\tilde{\Gamma}\times\{0\}$ and applying the definitions in \S\ref{sec:FunctionSpaces}, so that $H^s(\Gamma_\infty) := H^s(\R^{n})$, $H^s(\Gamma^{\circ}) := H^s(\tilde{\Gamma}^{\circ})$, $H^s_{\partial\Gamma}=H^s_{\partial\tilde{\Gamma}}$, etc..
For domains $\Omega\subset\R^{n+1}$ (e.g.\ the exterior domain $D=\R^{n+1}\setminus{\overline{\Gamma}}$)
we work with the classical Sobolev spaces\footnote{Our notation follows \cite{McLean}. Given a domain $\Omega\subset \R^{n+1}$, $W^1(\Omega){:=W^{1,2}(\Omega):=\{u\in L_2(\Omega):\nabla u\in L_2(\Omega)\}}$ is the Sobolev space whose norm is defined ``intrinsically'', via integrals over $\Omega$, while $H^1(\Omega)$ is defined (see \S\ref{sec:FunctionSpaces}) ``extrinsically'' as the set of restrictions to $\Omega$ of functions in $H^1(\R^{n+1})$. These spaces coincide if $\Omega$ is Lipschitz (e.g., \cite[Theorem 3.16]{McLean}), in particular if $\Omega = \R^{n+1}$, but not in general, in particular not, in general, 
for $\Omega = D$.}  $W^1(\Omega)$ and $W^1(\Omega,\Delta)$, normed by $\|u\|_{W^1(\Omega)}^2=\|u\|_{L_2(\Omega)}^2+\|\nabla u\|_{L_2(\Omega)}^2$ and $\|u\|_{W^1(\Omega,\Delta)}^2=\|u\|_{L_2(\Omega)}^2+\|\nabla u\|_{L_2(\Omega)}^2+ \|\Delta u\|_{L_2(\Omega)}^2$ respectively, and their ``local'' versions $W^{1,{\rm loc}}(\Omega)$ and $W^{1,{\rm loc}}(\Omega,\Delta)$, defined as the sets of measurable functions on $\Omega$ whose restrictions to any bounded domain $\Omega'\subset\Omega$ are in $W^{1}(\Omega')$ or $W^1(\Omega',\Delta)$ respectively. We denote by $\gamma^\pm:W^1(U^\pm)\to H^{1/2}(\Gamma_\infty)$ and $\dn^\pm:W^1(U^\pm,\Delta)\to H^{-1/2}(\Gamma_\infty)$ the standard Dirichlet and Neumann trace operators from the upper and lower half spaces $U^\pm:=\{x\in\R^{n+1}:\pm x_{n+1}>0\}$ onto the hyperplane $\Gamma_\infty$, where
the normal vector is assumed to point into $U^+$ in the case of the Neumann trace. Explicitly, the traces are the extension by density of $\gamma^\pm (u)(x):=\lim_{\substack{x'\to x\\x'\in U^\pm}}u(x')$ and $\dn^\pm u(x) := \lim_{\substack{x'\to x\\x'\in U^\pm}}\frac{\partial u}{\partial x_{n+1}}(x')$ for $u\in C^\infty_0(\R^{n+1})|_{U^\pm}$ 
and $x\in \Gamma_\infty$.
We note that if $u\in W^{1}(\R^{n+1})$ then $\gamma^+(u|_{U^+})=\gamma^-(u|_{U^-})$ \cite[\S2.1]{ScreenPaper}.
Finally,
we denote by $C^\infty_{0,\Gamma}$ the set of functions in $C^\infty_0(\R^{n+1})$ that equal one in a neighbourhood of $\Gamma$.

{In the case} 
that $\Gamma\subset\Gamma_\infty$ has non-empty (relative) interior $\Gamma^\circ$, 
we let 
$\cS:\tH^{-1/2}(\Gamma^\circ)\to C^2(D)\cap W^{1,{\rm loc}}(\R^{n+1})$ denote the (acoustic) single-layer potential operator, defined 
e.g.\ in \cite[\S2.2]{ScreenPaper}, 
which for a density $\phi\in L_2(\Gamma^\circ)$ has the integral representation
 \begin{align}
 \label{eq:SLPrep}
\cS\phi(\bx)=\int_{\Gamma^\circ}\Phi(\bx,\by)\phi(\by)\,\rd s(\by), \qquad \bx\in D,
\end{align}
where 
$\rd s$ 
denotes the usual (Lebesgue) surface measure and
 $\Phi(\bx,\by):=\re^{\ri k |\bx-\by|}/(4\pi |\bx-\by|)$ ($n=2$),  $\Phi(\bx,\by):=(\ri/4)H^{(1)}_0(k|\bx-\by|)$ ($n=1$), with $H_0^{(1)}$ the Hankel function of the first kind of order zero (e.g., \cite[Equation (9.1.3)]{AbramowitzStegun}). 
We let $S:\tH^{-1/2}(\Gamma^\circ)\to H^{1/2}(\Gamma^\circ)$ denote the single-layer boundary integral operator defined by 
$S\phi:=  \gamma^\pm(\sigma\cS\phi|_{U^\pm})|_{\Gamma^\circ}$, for $\phi\in \tH^{-1/2}(\Gamma^\circ)$, where $\sigma \in C^\infty_{0,\Gamma}$ is arbitrary. For $\phi\in L_\infty(\Gamma^\circ)$ we have the integral representation
\begin{align}
\label{eq:SLPrep2}
S\phi(\bx)=\int_{\Gamma^\circ}\Phi(\bx,\by)\phi(\by)\,\rd s(\by), \qquad \bx\in \Gamma^\circ.
\end{align}
Furthermore, $S$ 
is both continuous and coercive. 

\begin{lem}[{\cite{CoercScreen2}, \cite[\S2.2]{ScreenPaper}}] 
\label{lem:coer}
The sesquilinear form $a(\cdot,\cdot)$ on $\tH^{-1/2}(\Gamma^\circ)\times \tH^{-1/2}(\Gamma^\circ)$ defined by
\begin{align}
\label{eqn:Sesqui}
a(\phi,\psi):=\langle S\phi,\psi\rangle_{H^{1/2}(\Gamma^\circ) \times \tH^{-1/2}(\Gamma^\circ)},\qquad \phi, \psi\in \tH^{-1/2}(\Gamma^\circ),
\end{align}
is continuous and coercive: there exist constants $C_a, \alpha>0$, 
depending only on $k$ and $\diam(\Gamma)$, 
s.t.
\begin{equation} \label{eq:ContCoer}
|a(\phi,\psi)| \leq C_a\|\phi\|_{\tH^{-1/2}(\Gamma^\circ)}\, \|\psi\|_{\tH^{-1/2}(\Gamma^\circ)}, \quad |a(\phi,\phi)|\geq \alpha \|\phi\|_{\tH^{-1/2}(\Gamma^\circ)}^2, \quad \phi,\psi\in \tH^{-1/2}(\Gamma^\circ).
\end{equation}
\end{lem}


\subsection{The scattering problem}
\label{sec:Scattering} 
The scattering problem we consider, stated as Problem \ref{prob:BVP} below, is for the scattered field $u$, which is assumed to satisfy the Helmholtz equation %
\begin{align}
\label{eqn:HE}
\Delta u + k^2 u = 0, 
\end{align}
in $D:=\R^{n+1}\setminus {\overline{\Gamma}}$,
for some wavenumber $k>0$,
and the Sommerfeld radiation condition
\begin{align}
\label{eqn:SRC}
\pdone{u(x)}{r} - \ri k u(x) = o(r^{-n/2}), \qquad r:=|x|\to\infty, \text{ uniformly in } \hat x:=x/|x|.
\end{align}
We assume that the incident wave $u^i$ is an element of 
\[ W^i:=\{ u^i\in W^{1,{\rm loc}}(\R^{n+1}): \Delta u^i+k^2u^i=0 \text{ in some neighbourhood of }\Gamma\}.\]
For instance, $u^i$ might be a plane wave $u^i(x)=\re^{\ri k \vartheta\cdot x}$ for some $\vartheta\in\R^{n+1}$, $|\vartheta|=1$.

\begin{prob}
\label{prob:BVP}
Let $\Gamma\subset\Gamma_\infty$ be {bounded} 
with non-empty (relative) interior $\Gamma^\circ$. 
Given $k>0$ and $u^i\in W^i$, 
find $u\in C^2\left(D\right)\cap  W^{1,\mathrm{loc}}(D)$ satisfying
\rf{eqn:HE} in $D$, \rf{eqn:SRC},
and the boundary condition
\begin{align}\label{a1bc}
\gamma^\pm(\sigma u|_{U^\pm})|_{\Gamma^\circ}&=g:=-\gamma^\pm (\sigma u^i|_{U^\pm})|_{\Gamma^\circ}\in H^{1/2}(\Gamma^\circ),
\end{align}
{for some (and hence every) $\sigma\in C^\infty_{0,\Gamma}$.}
\end{prob}


%



Imposing the boundary condition by restriction \eqref{a1bc} to the (relative) interior of $\Gamma$ is the approach taken in many classic references on screen problems, e.g. \cite[Eqn.~(1.2)]{Ste:87}, \cite[Eqn.~(2.1)]{HoMaSt:96}, \cite[Eqn.~(1.1)]{StWe84} (for the analogous Laplace problem), \cite[Eqn.~(2.3)]{Ha-Du:90} (for the analogous wave equation problem). 
Alternatively, the boundary condition can be imposed in a weak sense.
\begin{prob}
\label{prob:BVPw}
Let $\Gamma\subset\Gamma_\infty$ be {bounded}.
Given $k>0$ and $u^i\in W^i$, 
find $u\in C^2\left(D\right)\cap  W^{1,\mathrm{loc}}(D)$ satisfying
\rf{eqn:HE} in $D$, \rf{eqn:SRC},
and the boundary condition
\begin{align}\label{bcweak}
u+u^i\in W^{1,\mathrm{loc}}_0(D).
\end{align}
\end{prob}

It is standard that Problem \ref{prob:BVPw} is well-posed, {for every bounded $\Gamma\subset \Gamma_\infty$} (see, e.g. \cite{ScreenPaper,HausdorffDomain}). {(Note that it is the compact set $\overline{\Gamma}$, rather than the bounded set $\Gamma$, that determines the solution to Problem \ref{prob:BVPw}.) The well-posedness of Problem \ref{prob:BVP} is characterised in the following result which is a sharpening of \cite[Theorem 6.2(a)]{ScreenPaper} and reduces to  \cite[Theorem 3.2]{BEMfract} in the case that $\Gamma=\Gamma^\circ$.
\begin{thm} \label{thm:wp} If $\Gamma\subset\Gamma_\infty$ is {bounded} 
with non-empty (relative) interior $\Gamma^\circ$, then Problem \ref{prob:BVP} has exactly one solution if and only if 
\begin{equation} \label{eq:ES}
\tH^{-1/2}(\Gamma^\circ)=H^{-1/2}_{\overline\Gamma}.
\end{equation} 
If \eqref{eq:ES} holds then Problems \ref{prob:BVP} and \ref{prob:BVPw} have the same unique solution. If \eqref{eq:ES} does not hold then the solution to Problem \ref{prob:BVPw} is a solution to Problem \ref{prob:BVP}, but Problem \ref{prob:BVP} has more than one solution. 
\end{thm}
\begin{proof} It is shown in \cite[Theorem 6.2(a)]{ScreenPaper} that Problem \ref{prob:BVP} has exactly one solution if and only if \eqref{eq:ES} holds and
$H^{1/2}_{\partial \Gamma}=\{0\}$. Now 
$$
\tH^{-1/2}(\Omega)\subset H^{-1/2}_{\overline{\Omega}} \subset H^{-1/2}_{\overline\Gamma},
$$ 
 where $\Omega:=\Gamma^\circ$. Thus \eqref{eq:ES} implies that $\tH^{-1/2}(\Omega)=H^{-1/2}_{\overline{\Omega}}$, so that, by \cite[Lemma 4.15(ii)]{caetano2019density}, $H^{1/2}_{\partial \Omega}=\{0\}$. Thus, if \eqref{eq:ES} holds and $\psi\in H^{1/2}_{\partial \Gamma}$, then  $\psi \in H^{-1/2}_{\partial \Gamma}\subset H^{-1/2}_{\overline\Gamma}=H^{-1/2}_{\overline{\Omega}}$, so that $\supp(\psi)\subset \overline{\Omega}\cap \partial \Gamma\subset \partial \Omega$, so that $\psi\in H^{1/2}_{\partial \Omega}=\{0\}$. Thus \eqref{eq:ES} implies that $H^{1/2}_{\partial \Gamma}=\{0\}$, i.e.\ the requirement that $H^{1/2}_{\partial \Gamma}=\{0\}$ in  \cite[Theorem 6.2(a)]{ScreenPaper} is redundant. (In the case $\Gamma=\Gamma^\circ$ this was noted already in  \cite[Theorem 3.2]{BEMfract}.) The remaining claims of the theorem follow from the observation  in \cite[Lemma 3.14]{ScreenPaper} that (regardless of whether or not \eqref{eq:ES} holds) the solution of Problem \ref{prob:BVPw} is a solution of Problem \ref{prob:BVP}.
\end{proof}
}
 
Combining this theorem with Corollary \ref{cor:tildesubscriptHs} we obtain the first part of the following result. The remainder of the theorem, a reformulation of Problems \ref{prob:BVP} and \ref{prob:BVPw} as a {boundary integral equation (BIE)}  for the jump in the normal derivative of $u$ across $\Gamma$, follows from \cite[Theorem 3.29]{ScreenPaper} applied with (in the notation of \cite{ScreenPaper}) $V^-:= H^{-1/2}(\Gamma^\circ)$.


\begin{thm}[{Consequence of Corollary  \ref{cor:tildesubscriptHs}, Theorem \ref{thm:wp}, and \cite[Theorem 3.29]{ScreenPaper}}] 
\label{thm:Closed}
Let $\Gamma\subset\Gamma_\infty$ be an $n$-attractor. 
%
Then Problem \ref{prob:BVP}
has a unique solution $u$, which coincides with the unique solution of Problem \ref{prob:BVPw}, and %
satisfies
\begin{align}
\label{eqn:Rep}
u(x )= -\cS\phi(x), \qquad x\in D,
\end{align}
where $\phi = \dn^+(\sigma u|_{U^+})-\dn^-(\sigma u|_{U^-}) \in \tH^{-1/2}(\Gamma^\circ)$ (with $\sigma \in C^\infty_{0,\Gamma}$ arbitrary) is the unique solution of 
\begin{equation}
\label{eqn:BIE}
S\phi = -g.
\end{equation}
\end{thm}

\begin{rem}[Relating \eqref{eqn:BIE} to the BIE in \cite{HausdorffBEM}]
The fact that $\tH^{-1/2}(\Gamma^\circ)=H^{-1/2}_{\overline{\Gamma}}$ means that the BIE \eqref{eqn:BIE} coincides with the BIE \cite[Eqn.~(34)]{HausdorffBEM}, which is stated with $\tH^{-1/2}(\Gamma^\circ)$ replaced by $H^{-1/2}_{\overline{\Gamma}}$. 
\end{rem}

\begin{rem}[The case when $\Gamma^\circ$ is an open $n$-set] \label{rem:nset} By the results in Remark \ref{re:tildesubscript}, in particular equation \eqref{eq:-s<0}, Theorem \ref{thm:Closed} holds by the same arguments if the assumption that $\Gamma$ is an $n$-attractor is replaced by the weaker assumption that $\Gamma^\circ$ is an open $n$-set and $\overline{\Gamma^\circ}=\overline{\Gamma}$. These weaker assumptions hold if $\Gamma$ is an $n$-attractor by Proposition \ref{prop:nset}. They also hold in the classical case that $\Gamma = \Gamma^\circ$ is a bounded Lipschitz domain and $\overline{\Gamma}=\overline{\Gamma^\circ}$ is its closure.
\end{rem}

\begin{rem}[Examples of bounded sets $\Gamma$ for which \eqref{eq:ES} fails] \label{rem:fails} While Theorem \ref{thm:Closed} and Remark \ref{rem:nset} exhibit cases  where \eqref{eq:ES} hold, it is easy to construct bounded sets $\Gamma\subset \Gamma_\infty$, with non-empty relative interior $\Gamma^\circ$, for which \eqref{eq:ES} fails. For example, suppose that  $\Gamma:=\Gamma_1 \cup \Gamma_2$, where $\Gamma_1$ and $\Gamma_2$ are disjoint compact sets,  $\Gamma^\circ_1\neq \emptyset$,
$\Gamma_2^\circ=\emptyset$, and $H^{-1/2}_{\Gamma_2}\neq \{0\}$ (these last two properties hold, for example, if $n-2<\dimH(\Gamma_2)<n$; see, e.g., \cite[Theorem 2.12]{HewMoi:15}). Then $H^{-1/2}_{\Gamma_2}\subset H^{-1/2}_\Gamma$, $\tH^{-1/2}(\Gamma^\circ)\subset H^{-1/2}_{\Gamma_1}$, and $H^{-1/2}_{\Gamma_1}\cap H^{-1/2}_{\Gamma_2}=\{0\}$, so that $H^{-1/2}_\Gamma \not\subset\tH^{-1/2}(\Gamma^\circ)$. 
A concrete example in the case $n=2$ would be where $\Gamma_1$ is a closed disc and $\Gamma_2$ is a Koch curve, with $\Gamma_1\cap \Gamma_2=\emptyset$.
It is more challenging to construct examples where $\Gamma^\circ$ is the interior of $\overline{\Gamma}$ and is dense in $\overline{\Gamma}$, as is the case, by Proposition \ref{prop:nset},  if $\Gamma$ is an $n$-attractor. We construct such examples in Appendix \ref{sec:Appendix2}; see Corollary \ref{cor:exnset}.
\end{rem}

The following regularity result follows from Corollaries \ref{cor:tildesubscriptHs} and \ref{cor:Interpolationnset} and results in \cite{HausdorffDomain}.

\begin{prop} \label{rem:reg}
Let $\Gamma\subset\Gamma_\infty$ be an $n$-attractor. Where $\phi\in \tH^{-1/2}(\Gamma^\circ)$ is the unique solution of \eqref{eqn:BIE},  it holds,  for some $\varepsilon\in (0,1/2]$ independent of the incident field $u^i\in W^i$,  
that $\phi\in \tH^s(\Gamma^\circ)$, for $-1/2\leq s<-1/2+\varepsilon$.
\end{prop}
\begin{proof}
Recall that  $\tH^s(\Gamma^\circ):= \tH^s(\tilde \Gamma^\circ)\subset H^s(\R^n)$, let $\varphi:=-\phi \otimes \delta \in {\cal S}'(\R^{n+1})$, so that $\supp(\varphi)\subset\Gamma$,  and note that, for $-1<s<0$, $\phi\in H^s(\R^n)$ if and only if $\varphi\in H^{s-1/2}(\R^{n+1})$, by \cite[Lemma 3.39]{McLean}. It follows that $\varphi\in H^{-1}(\R^{n+1})$ and that, where $u$ is as in Theorem \ref{thm:wp},  $u(x)=-\cS\phi(x) = \cA\varphi(x)$, $x\in D$, where $\cA \varphi$ is as defined in  \cite[(3.7)]{HausdorffDomain}. But also, by \cite[Theorem 3.4]{HausdorffDomain}, $u(x) = \cA \widetilde \varphi(x)$, $x\in D$, where $\widetilde \varphi\in H^{-1}(\R^{n+1})$, with $\supp(\widetilde \varphi)\subset \Gamma$, is the unique solution of \cite[(3.15)]{HausdorffDomain}, with $g$ in that equation given by \cite[(3.6)]{HausdorffDomain}. Since $\cA\varphi = \cA\widetilde \varphi$ almost everywhere in $\R^{n+1}$, it follows from \cite[(3.9)]{HausdorffDomain} that $\varphi=\widetilde \varphi$. Further, \cite[Proposition 3.19, Remark  3.20]{HausdorffDomain} give, for some $\varepsilon\in (0,1/2]$ independent of the incident field $u^i\in W^i$, that $\widetilde \varphi\in H^{s-1/2}(\R^{n+1})$, for $-1/2\leq s<-1/2+\varepsilon$. (Note that \cite[Proposition 3.19]{HausdorffDomain} relies on our interpolation result, Corollary \ref{cor:Interpolationnset}, as noted in \S\ref{sec:app}.) Thus, for  $-1/2\leq s<-1/2+\varepsilon$, $\phi\in H^s_{\Gamma}=\tH^s(\Gamma^\circ)$, by Corollary \ref{cor:tildesubscriptHs}.
\end{proof}

\subsection{The {boundary element method}}
\label{sec:BEM}

To solve the {boundary integral equation (BIE)} \eqref{eqn:BIE} numerically, we follow the {boundary element method (BEM)} approach described in \cite[\S5]{HausdorffBEM}.\footnote{Our notation in what follows is slightly different to that in \cite[\S5]{HausdorffBEM} because here we work with (relatively) open elements $\Omega_\bm$ (in line with classical FEM/BEM) whereas \cite[\S5]{HausdorffBEM} works with closed elements $\Gamma_\bm$ (because for the more general fractal screen problems discussed in \cite{HausdorffBEM} the elements may have empty interior). But in our context it makes no difference because $\Gamma_\bm=\overline{\Omega_\bm}$ and $|\partial\Omega_\bm|=|\partial\Gamma_\bm|=0$, where $|\cdot|$ denotes Lebesgue surface measure on $\Gamma_\infty$.}
Let $\Gamma$ be an $n$-attractor, and for $0<h\leq 
{\diam(\Gamma)}$ let $L_h$ and $\widetilde{V}_h$ be defined as in Corollary \ref{cor:approx}. 
Our numerical method for solving \rf{eqn:BIE} uses $\widetilde{V}_h$ as the approximation space in a Galerkin method. 
That is, we seek $\phi_h\in \widetilde{V}_h$ such that (with $a$ defined by \rf{eqn:Sesqui})
\begin{align}
\label{eqn:Variational}
a(\phi_h,\psi_h)
=-\langle g,\psi_h\rangle_{H^{1/2}(\Gamma^\circ)\times \tH^{-1/2}(\Gamma^\circ)}, \qquad \psi_h\in \widetilde{V}_h.
\end{align}

Let $N:= \# L_h$ and let $\{\bm(i)\}_{i=1}^N$ be an indexing of the elements of $L_h$. Let $\{\phi_i\}_{i=1}^N$ be the canonical $L_2(\Gamma^\circ)$-orthonormal basis for $\widetilde{V}_h$, i.e.\ (with $\chi_{\Omega_{\bm(i)}}$ the characteristic function of $\Omega_{\bm(i)}$)
\[ \phi_i = |\Omega_{\bm(i)}|^{-1/2}\chi_{\Omega_{\bm(i)}}, \quad i=1,\ldots,N.\]
Then $\phi_h = \sum_{i=1}^N c_i \phi_i$, where  
the coefficient vector $\vec{c}=(c_1,\ldots,c_N)^T\in\C^N$ satisfies the system
\begin{equation} \label{eq:cvec}
A \vec{c} = \vec{b},
\end{equation}
where the entries of $A$ and $b$ are given, for $i,j\in \{1,\ldots,N\}$, by
\begin{align}
\label{eqn:GalerkinElements}
A_{ij}&=|\Omega_{\bm(i)}|^{-1/2}|\Omega_{\bm(j)}|^{-1/2}\int_{\Omega_{\bm(i)}}\int_{\Omega_{\bm(j)}} \Phi(x,y) \, \rd s(y)\rd s(x),
\quad
b_{i}&=-|\Omega_{\bm(i)}|^{-1/2}\int_{\Omega_{\bm(i)}} g(x)\, \rd s(x).
\end{align}


The following basic convergence result is a direct consequence of \cite[Theorem~5.1]{HausdorffBEM} and its proof.

\begin{thm}
\label{thm:Convergence}
Let $\Gamma\subset\Gamma_\infty$ be an $n$-attractor. Then for each $0<h\leq 
{\diam(\Gamma)}$ the variational problem \rf{eqn:Variational}
has a unique solution $\phi_h\in \widetilde{V}_h$, satisfying the quasioptimality estimate
\begin{align}
\label{eq:quasiopt}
\|\phi-\phi_h\|_{H^{-1/2}(\Gamma_\infty)} \leq \frac{C_a}{\alpha}\inf_{v_h\in \widetilde{V}_h}\|\phi-v_h\|_{H^{-1/2}(\Gamma_\infty)},
\end{align}
where $\phi\in \tH^{-1/2}(\Gamma^\circ)$ denotes the solution of \rf{eqn:BIE}, 
{and $C_a$, $\alpha$ are as in Lemma \ref{lem:coer}.}  
Furthermore, $\|\phi - \phi_h\|_{H^{-1/2}(\Gamma_\infty)} \to 0$ as $h\to 0$. 
\end{thm}

\subsection{Galerkin error estimates}
\label{sec:ErrorEstimates}

We now use the best approximation error results proved in \S\ref{sec:Approximation} to give an error bound for the Galerkin approximation to the BIE \eqref{eqn:BIE} when the solution
$\phi$ is sufficiently smooth. 
\begin{thm}
\label{thm:ConvergenceRate}
Let $\Gamma$, $\phi$, $\phi_h$ be as in Theorem \ref{thm:Convergence}. Suppose that $\phi\in \tH^{s}(\Gamma^\circ)$ for some $-1/2<s<0$. 
Then, 
for $0<h\leq 
{\diam(\Gamma)}$,
\begin{align}
\label{}
\|\phi-\phi_h\|_{H^{-1/2}(\Gamma_\infty)} \leq c
h^{s+1/2}\|\phi\|_{\tH^{s}(\Gamma^\circ)},
\end{align}
{for some constant $c>0$ independent of $h$ and $\phi$.}
\end{thm}

\begin{proof}
Combine the quasioptimality result \eqref{eq:quasiopt} with the best approximation error estimate \eqref{eq:Hsboundh_dn2}, which holds by Corollary \ref{cor:approx}.
\end{proof}

 Note that, for some $\varepsilon \in (0,1/2]$, $\phi\in \tH^{s}(\Gamma^\circ)$ for $-1/2<s<-1/2+\varepsilon$, by Proposition \ref{rem:reg}. 
Determining the largest $\varepsilon$ for which this holds
appears to be an open problem for a general $n$-attractor. In the special cases where $\Gamma$ is an interval ($n=1$) or a square ($n=2$), it is known that $\phi\in \tH^{s}(\Gamma^\circ)$ for all $-1/2<s<0$ 
by classical regularity results for first kind boundary integral equations, see \cite{StWe84,Ste:87,ErStEl90}. We note, however, that, even in these cases, generically $\phi\not\in \tH^0(\Omega)=L_2(\Omega)$, 
because $\phi$ has inverse square root singularities at the endpoints ($n=1$) or edges ($n=2$) of $\Omega$, 
see \cite{StWe84,Ste:87,ErStEl90}. 

\appendix

{
\section{Table of notation}
\label{sec:Notation}


\begin{longtable}{ll}
	\multicolumn{2}{l}{\textbf{Sets and measures on $\R^n$}}\\[4pt]
	$B(x,r)$, $B^\circ(x,r)$ & closed, resp.\ open, ball of radius $r$ centred at $x$\\
	$|E|$ & Lebesgue measure of $E \subset \mathbb{R}^n$\\
	$\mathcal{H}^d$ & $d$-dimensional Hausdorff measure\\
	$\partial\Omega$, $\overline{\Omega}$, $\Omega^\circ$, $\Omega^c$ & boundary, closure, interior, complement of $\Omega$\\
	$\operatorname{dist}(y,E)$ & distance from $y$ to the set $E$\\
	$\operatorname{diam}(E)$ & diameter of $E$\\
	$\chi_\Omega$ & characteristic function of $\Omega$\\[6pt]

	\multicolumn{2}{l}{\textbf{Dimensions}}\\[4pt]
	$\dim_{\mathrm{H}}(E)$ & Hausdorff dimension of $E\subset\R^n$\\
	$\dim_{\mathrm{A}}(E)$ & Aikawa/Assouad dimension of $E$ \\
	$\dim_{\mathrm{M}}(E)$, $\underline{\dim}_{\mathrm{M}}(E)$, $\overline{\dim}_{\mathrm{M}}(E)$ & Minkowski dimension of $E$, lower and upper variants \\[6pt]

	\multicolumn{2}{l}{\textbf{Properties of subsets of $\R^n$}}\\[4pt]
	domain & non-empty open set\\
	$d$-set & set satisfying \eqref{eq:dset}\\
	$n$-set, open $n$-set & set satisfying \eqref{eq:nset}, resp.\ open set satisfying \eqref{eq:nset}\\
	porous & Definition \ref{def:Porous}\\
	$E$-porous &  Definition \ref{def:EPorous} \\
		thick & Both $I$-thick and $E$-thick (Definition \ref{def:Thick})\\
	$h$-set & Definition \ref{def:hset}\\
	uniformly porous & porous and an $h$-set (Definition \ref{def:up})\\
	$\mathcal{D}^t$ & (variant of) Frazier--Jawerth class (see \eqref{eq:D_t})\\
	[6pt]
	
	\multicolumn{2}{l}{\textbf{IFS on $\R^n$ and $n$-attractors}}\\[4pt]
	IFS & iterated function system\\
	$\{s_1,\ldots,s_M\}$ & IFS of contracting similarities\\
	$\rho_m$ & contraction ratio of $s_m$\\
	$\rho_{\min}$, $\rho_{\max}$ & $\min_m \rho_m$, $\max_m \rho_m$\\
	$\Gamma$ & attractor of the IFS\\
	$d$ & similarity dimension (solution of $\sum \rho_m^d = 1$ --- see \eqref{eq:dfirst})\\
	OSC & open set condition (see \eqref{oscfirst})\\
	$n$-attractor & IFS attractor satisfying OSC with similarity dimension $n$ (Definition \ref{def:1})\\
	$h_0$ & $\operatorname{diam}(\Gamma)$\\
	$I_\ell$ & set of multi-indices $\{1,\ldots,M\}^\ell$\\
	$\boldsymbol{m} = (m_1,\ldots,m_\ell)$ & multi-index in $I_\ell$\\
	$s_{\boldsymbol{m}}$, $\Gamma_{\boldsymbol{m}}$, $\Omega_{\boldsymbol{m}}$ & $s_{m_1}\circ\cdots\circ s_{m_\ell}$,\; $s_{\boldsymbol{m}}(\Gamma)$,\; $s_{\boldsymbol{m}}(\Omega)$, for $\bm \in I_\ell$\\
	$L_h$ & index set for decomposition at scale $h$ (see \eqref{eq:LhDef})\\[6pt]

	\multicolumn{2}{l}{\textbf{Function spaces and distributions}}\\[4pt]
	$\mathcal{S}(\mathbb{R}^n)$, $\mathcal{S}'(\mathbb{R}^n)$ & Schwartz functions, tempered distributions\\
	$\mathcal{D}(\Omega) = C_0^\infty(\Omega)$, 	$\mathcal{D}'(\Omega)$ & test functions on $\Omega$, 
distributions on $\Omega$\\
	$\langle\cdot,\cdot\rangle$ & duality pairing on $\cD'(\Omega)\times \cD(\Omega)$ (also that on $\mathcal{S}'(\mathbb{R}^n)\times \mathcal{S}(\mathbb{R}^n)$)	\\
	$H^s(\mathbb{R}^n)$  ($=H^s_2(\R^n)$) & 
	$L_2$ fractional Sobolev/Bessel potential space of order $s$\\
	$H^s(\Omega)$ & restriction space $\{u|_\Omega : u \in H^s(\mathbb{R}^n)\}$\\
	$\widetilde{H}^s(\Omega)$ & $\overline{C_0^\infty(\Omega)}^{H^s(\mathbb{R}^n)}$\\
	$H^s_E$ & $\{u \in H^s(\mathbb{R}^n) : \operatorname{supp} u \subset E\}$\\
	$H^s_p(\mathbb{R}^n)$ 
($= F^s_{p,2}(\mathbb{R}^n)$ for $1<p<\infty$)	
		& $L_p$ fractional Sobolev/Bessel potential space of order $s$\\
	$B^s_{p,q}(\mathbb{R}^n)$ & Besov space\\
	$F^s_{p,q}(\mathbb{R}^n)$ & Triebel--Lizorkin space\\
	$A^s_{p,q}$ & either $B^s_{p,q}$ or $F^s_{p,q}$\\
	$A^s_{p,q}(\Omega)$, $\widetilde{A}^s_{p,q}(\Omega)$, $A^s_{p,q,E}$ & restriction, tilde, supported-on-$E$ variants (see \S\ref{sec:FunctionSpaceNotation})\\
$E$ is $A^s_{p,q}$-null & $A^s_{p,q,E}=\{0\}$	\\
	$\mathring{A}_{p,q}^{s}(\Omega)$ & $\overline{C^\infty_0(\Omega)}^{A_{p,q}^{s}(\Omega)}$\\
	$W^1(\Omega)$ ($=W^{1,2}(\Omega)$) & $\{u \in L_2(\Omega) : \nabla u \in L_2(\Omega)\}$\\
	$W^{1,\mathrm{loc}}(\Omega)$ & local version of $W^1(\Omega)$\\
	$V_h$, $\widetilde{V}_h$ & space of piecewise constants on a mesh, resp.\ extended by zero
	\\[6pt]
	
	\multicolumn{2}{l}{\textbf{Operators}}\\[4pt]
	$|_\Omega$ & (distributional) restriction to $\Omega$\\
		$\operatorname{ext}$ & (distributional) extension by zero (see \eqref{eq:extbyzero})\\
	$P_h$ & $L_2$-orthogonal projection onto $V_h$\\
	$(\cdot,\cdot)_{\theta,q}$ & real interpolation with parameters $\theta$, $q$\\
	$[\cdot,\cdot]_\theta$ & complex interpolation with parameter $\theta$\\[6pt]
	
	\multicolumn{2}{l}{\textbf{Scattering application (\S\ref{sec:BVPsBIEs}) }}\\[4pt]
	$\Gamma_\infty$ & $\mathbb{R}^n \times \{0\}$ (hyperplane)\\
	$D$ & $\mathbb{R}^{n+1} \setminus \overline{\Gamma}$ (propagation domain)\\
	$U^\pm$ & upper/lower half-spaces $\{x\in\R^{n+1}:\pm x_{n+1}>0\}$\\
	$k$ & wavenumber\\
	$u^i$, $W^i$ & incident wave, space of incident waves\\
	$\Phi(x,y)$ & fundamental solution of Helmholtz equation\\
	$\gamma^\pm$, $\partial_n^\pm$ & Dirichlet and Neumann trace operators from $U^\pm$ to $\Gamma_\infty$\\
$C^\infty_{0,\Gamma}$ & functions in $C^\infty_0(\R^{n+1})$ equal to one in a neighbourhood of $\Gamma$	
	\\
	$\mathcal{S}$ & single-layer potential operator (see \eqref{eq:SLPrep})\\
	$S$ & single-layer boundary integral operator (see \eqref{eq:SLPrep2})\\
	$C_a$, $\alpha$ & continuity and coercivity constants for $S$ (see \eqref{eq:ContCoer})\\
	BIE & boundary integral equation\\
	BEM & boundary element method\\
\end{longtable}
}

\section{Connection to Caetano et al.~(2025)}
\label{sec:DomainComment}

{
In our earlier paper \cite{HausdorffDomain} we referenced a number of results from the current paper. But the notation in \cite{HausdorffDomain} differs slightly from that used here, so for the benefit of readers of \cite{HausdorffDomain} we provide an explanation of the connection between the two papers. 

In \cite{HausdorffDomain} the notation $\IH^t(\Gamma)$, $t\in \R$, is used to denote a certain family of trace spaces on a compact $d$-set $\Gamma\subset\R^n$ for general $0<d\leq n$ (see \cite[\S3.2]{HausdorffDomain} for definitions). In the special case where $\Gamma$ is an $n$-attractor (in the terminology of the current paper), 
let $\Omega:=\Gamma^\circ$. Then for $s\geq0$ the space $\IH^s(\Gamma)$ of \cite{HausdorffDomain} coincides with the restriction space $H^s(\Omega)$ of the current paper, with equal norms, because $\tr$ coincides with $|_\Gamma$ and $|\partial\Gamma|=|\partial\Omega|=0$ (cf.\ the proof of \cite[Theorem 3.9]{HausdorffDomain}). Hence, by duality (see \cite[\S3(b)]{HausdorffDomain} and \S\ref{sec:FunctionSpaceNotation}), for $s\leq 0$ the space $\IH^s(\Gamma)$ of \cite{HausdorffDomain} coincides with the tilde space $\tH^s(\Omega)$ of the current paper, again with equal norms. With these facts in mind, 
the specific results from the current paper used in \cite{HausdorffDomain} are justified as follows:
\begin{itemize}
\item The claim in \cite[Remark~3.14(ii)]{HausdorffDomain} that if $\Gamma$ is an $n$-attractor then $\tH^1(\Gamma^c)=H^1_{\overline{\Gamma^c}}$ is proved in Corollary \ref{cor:tildesubscriptHs} of the current paper. 
\item The claim in {the proof of} \cite[Proposition~3.19]{HausdorffDomain} that if $\Gamma$ is an $n$-attractor then $\{\IH^s(\Gamma)\}_{s\geq 0}$ and $\{\IH^s(\Gamma)\}_{s\leq 0}$ are interpolation scales is proved in 
{Corollary} 
\ref{cor:Interpolationnset} 
of the current paper. 
\item The result \cite[(4.19)]{HausdorffDomain} {(with $t=1$ in the notation of \cite{HausdorffDomain})} used in the proof of \cite[Theorem 4.5]{HausdorffDomain} in the case where $\Gamma$ is an $n$-attractor is equivalent to the claim that, 
for $0<s<1$ and $0<h\leq \diam(\Gamma)$,
\[ \inf_{u_h\in \widetilde V_h}\| u-u_h\|_{\tH^{-1}(\Omega)}\leq c\,h^{1-s}\| u\|_{\tH^{-s}(\Omega)}, \qquad u\in \tH^{-s}(\Omega),
\]
where $\widetilde V_h$ is a space of piecewise constant functions (extended by zero so as to be defined on $\R^n$) on a certain 
mesh of $\Omega$. 
This is proved in Corollary \ref{cor:approx} of the current paper (where $\widetilde V_h$ is also defined), 
{being equivalent to \eqref{eq:Hsboundh_dn2} with $s_1=-1$.} 
\end{itemize} 

In relation to the results of \cite{HausdorffDomain} we mention that in recent work (\cite[Remark 3.3]{Hinz}) it has been shown (in the notation of \cite{HausdorffDomain}) that $\mathrm{ker}(\mathrm{tr}_\Gamma)=\tH^1(\Gamma^c)$ whenever $\Gamma$ is an $n$-set. This means that \cite[Assumption 3.12]{HausdorffDomain} can be replaced simply by the assumption that $\Gamma$ is a $d$-set for some $n-2<d\leq n$.}

\section{Thickness, porosity, $h$-sets, and uniform porosity}
\label{sec:Appendix}

We recall the following definitions from \cite{Tri08}.

\begin{defn}[Porosity (ball condition), {\cite[Definition~3.4(i) \& Remark~3.5]{Tri08}}]
\label{def:Porous}
A closed set $\Gamma\subset\R^n$ is said to be \emph{porous}, or to \emph{satisfy the ball condition}, if there exists $\eta\in(0,1)$ such that for every $x\in \R^n$ and $r\in (0,1)$ there exists $y\in \R^n$ such that 
\begin{align}
\label{eq:PorousDefn}
B(y,\eta r)\subset B(x,r) \quad \text{and} \quad B(y,\eta r)\cap \Gamma = \emptyset.
\end{align}
\end{defn}

References to literature on the relevance of porosity in fractal geometry and function spaces can be found in \cite[Remark~3.17]{Tri08}. 

\begin{defn}[$E$-porosity, {\cite[Definition~3.16(i)]{Tri08}}]
\label{def:EPorous}
A domain $\Omega\subsetneqq\R^n$ is said to be \emph{$E$-porous} if there exists $\eta\in(0,1)$ such that for every $x\in \partial\Omega$ and $r\in (0,1)$ there exists $y\in \R^n$ such that 
\begin{align}
\label{eq:EPorousDefn}
B(y,\eta r)\subset B(x,r) \quad \text{and} \quad B(y,\eta r)\cap \overline{\Omega}= \emptyset.
\end{align}
\end{defn}

As mentioned in \cite[Remark~3.5 \& Remark~3.17]{Tri08}, the boundary $\Gamma$ of an $E$-porous domain $\Omega$ is porous. 

In the next definition, for a set $S\subset\R^n$, we denote by $\mathcal Q(S)$ the set of the (open) cubes contained in $S$ and with the edges parallel to the Cartesian axes; and for any $Q\in\mathcal Q(S)$ we denote by $l(Q)$ the length of its edges.
\begin{defn}[Thickness, {\cite[Definition~3.1(ii)--(iv), Remark~3.2]{Tri08}}]
\label{def:Thick}
Let $\Omega\subsetneqq \R^n$ be a domain.\begin{enumerate}[(i)]
\item $\Omega$ is said to be {\em$E$-thick (exterior thick)} if for any choice
of $c_{1},c_{2},c_{3},c_{4}>0$ and $j_{0}\in\N$ there are $c_{5},c_{6},c_{7},c_{8}>0$
such that for any $j\in\N$, $j\geq j_{0}$, and any \emph{interior
cube} $Q^{i}\in \mathcal{Q}(\Omega)$
with
\[
c_{1}2^{-j}\leq l(Q^{i})\leq c_{2}2^{-j}\quad\mbox{ and }\quad{\rm c_{3}2^{-j}\leq dist}(Q^{i},\partial\Omega)\leq c_{4}2^{-j},
\]
there exists an \emph{exterior cube} $Q^{e}\in\mathcal{Q}(\Omega^{c})$ with
\[
c_{5}2^{-j}\leq l(Q^{e})\leq c_{6}2^{-j}\quad\mbox{ and }\quad c_{7}2^{-j}\leq{\rm dist}(Q^{e},\partial\Omega)\leq{\rm dist}(Q^{i},Q^{e})\leq c_{8}2^{-j}.
\]

\item $\Omega$ is said to be {\em$I$-thick (interior thick)} if for any choice
of $c_{1},c_{2},c_{3},c_{4}>0$ and $j_{0}\in\N$ there are $c_{5},c_{6},c_{7},c_{8}>0$
such that for any $j\in\N$, $j\geq j_{0}$, and any exterior cube
$Q^{e}\in\mathcal{Q}(\Omega^{c})$ with
\[
c_{1}2^{-j}\leq l(Q^{e})\leq c_{2}2^{-j}\quad\mbox{ and }\quad c_{3}2^{-j}\leq{\rm dist}(Q^{e},\partial\Omega)\leq c_{4}2^{-j},
\]
there exists an interior cube $Q^{i}\in\mathcal{Q}(\Omega)$ with
\[
c_{5}2^{-j}\leq l(Q^{i})\leq c_{6}2^{-j}\quad\mbox{ and }\quad c_{7}2^{-j}\leq{\rm dist}(Q^{i},\partial\Omega)\leq{\rm dist}(Q^{e},Q^{i})\leq c_{8}2^{-j}.
\]
\item $\Omega$ is said to be {\em thick} if it is both $E$-thick and $I$-thick.
\end{enumerate}
\end{defn}


It is straightforward to prove that E-porosity implies E-thickness - see, e.g., \cite[Proposition 3.18]{Tri08}. 
A converse to this implication is provided by the following result, which may be of independent interest.  
Its proof is similar to the proof of \cite[Proposition 5.10]{caetano2019density}.
\begin{prop}
\label{prop:ThickPorous}
If $\Omega\subsetneqq\R^n$ is an $E$-thick domain and $\partial\Omega$ is porous (satisfies the ball condition), then $\Omega$ is $E$-porous.
\end{prop}

\begin{proof}
We want to prove that there exists $0<\gamma<1$ such that, for any
$x\in\partial\Omega$ and any $0<\ell<1$, there exists $y$ such
that $B(y,\gamma\ell)\subset B(x,\ell)$ and $B(y,\gamma\ell)\cap\overline{\Omega}=\emptyset$.

Following the assumption that $\partial\Omega$ satisfies the ball
condition, let $\eta\in(0,1)$ be such that, for any ball $B(x,r)$
centered at $x\in\partial\Omega$ and with radius $r\in(0,1)$, there
exists a ball $B(y,\eta r)$ such that 
\[
B(y,\eta r)\subset B(x,r)\quad\text{and}\quad B(y,\eta r)\cap\partial\Omega=\emptyset.
\]
Following \cite[Remark 5.9]{caetano2019density}, we assume that the
above $\eta$ has been chosen such that we even have that 
\[
\dist(B(y,\eta r),\partial\Omega)\geq\eta r.
\]

Fix $j_{0}=1$, $c_{1}=2\eta/\sqrt{n}$, $c_{2}=4\eta/\sqrt{n}$,
$c_{3}=\eta$ and $c_{4}=2$ in the definition of $E$-thickness applied
to $\Omega$. Consider then the constants $c_{5}$, $c_{6}$, $c_{7}$
and $c_{8}$ that come out from that definition and set $c=1/(\sqrt{n}\,c_{6}+c_{8}+1)$.

Let $x$ be any point of $\partial\Omega$ and consider $B(x,\ell)$
for any $0<\ell<1$. Applying the ball condition assumption of $\partial\Omega$
to $r=c\ell$, there exists $B(y,\eta c\ell)$ such that 
\[
B(y,\eta c\ell)\subset B(x,c\ell)\quad\text{and}\quad\dist\big(B(y,\eta c\ell),\partial\Omega\big)\geq\eta c\ell.
\]

One of the following two situations must happen: 
\[
B(y,\eta c\ell)\subset\Omega\quad\text{or}\quad B(y,\eta c\ell)\subset\overline{\Omega}^{c}.
\]
In the second case we have what we want with $\gamma=\eta c$. In
the first case, start by considering the open cube $Q(y,\eta c\ell/\sqrt{n})$
with edges parallel to the Cartesian axes, centred at $y$ and with
side length $2\eta c\ell/\sqrt{n}$. Being inside $B(y,\eta c\ell)$,
it is also contained in $\Omega$. Consider also $j\in\N$ such that
$2^{-j}<c\ell\leq2^{-(j-1)}$ and observe that (where $l$ stands
for the side length of the corresponding cube)
\[
c_{1}2^{-j}\leq l\big(Q(y,\eta c\ell/\sqrt{n})\big)\leq c_{2}2^{-j}\quad\mbox{ and }\quad c_{3}2^{-j}\leq{\rm dist}\big(Q(y,\eta c\ell/\sqrt{n}),\partial\Omega\big)\leq c_{4}2^{-j}
\]
for the constants $c_{1}$, $c_{2}$, $c_{3}$ and $c_{4}$ given
above. Then such $Q(y,\eta c\ell/\sqrt{n})$ is a so-called interior
cube with respect to the $E$-thick domain $\Omega$, therefore there
exists a so-called exterior open cube $Q^{e}\subset\Omega^{c}$ (actually,
$Q^{e}\subset\overline{\Omega}^{c}$) such that 
\[
c_{5}c\ell/2\leq l(Q^{e})<c_{6}c\ell\quad\mbox{ and }\quad c_{7}c\ell/2\leq{\rm dist}(Q^{e},\partial\Omega)\leq{\rm dist}\big(Q(y,\eta c\ell/\sqrt{n}),Q^{e}\big)<c_{8}c\ell.
\]
For $z\in Q^{e}$ it holds by the choice of $c$ that 
\begin{align*}
|x-z| & \leq\diam(Q^{e})+\dist\big(Q(y,\eta c\ell/\sqrt{n}),Q^{e}\big)\\
 & \quad+\sup_{w\in Q(y,\eta c\ell/\sqrt{n})}|w-x|<(\sqrt{n}\,c_{6}+c_{8}+1)c\ell=\ell,
\end{align*}
so that $Q^{e}\subset B(x,\ell)$. Let $y'$ be the centre of $Q^{e}$
and consider $B(y',\frac{c_{5}c\ell}{4})$, so that we have
\[
B\left(y',\frac{c_{5}c\ell}{4}\right)\subset B\left(y',\frac{l(Q^{e})}{2}\right)\subset Q^{e}\subset B(x,\ell)\quad\mbox{and}\quad B\left(y',\frac{c_{5}c\ell}{4}\right)\cap\overline{\Omega}=\emptyset.
\]
The proof is completed by choosing $\gamma=\min\{\eta c,\frac{c_{5}c}{4}\}$.
\end{proof}

A key notion we have introduced in \S\ref{sec:nsets} is that of a $d$-set. A related, more general, notion is that of an $h$-set. Our definition below is essentially that of the originator Bricchi \cite[Definition 2.2]{bricchi2003complements} {(see also \cite{bricchi2002existence})}, except that, to fit with the applications we make: i) our $F$ need not be compact; ii) we require that $h(0)=0$.
\begin{defn}[$h$-set, {\cite[Definition 2.2]{bricchi2003complements}}] \label{def:hset} Suppose that $h:[0,1]\to [0,1]$ is continuous and strictly increasing with $h(0)=0$ and $h(1)=1$. Then a closed set $F\subset \R^n$ is said to be an {\em $h$-set} if there exists a Radon measure $\mu$ on $\R^n$, with $\supp(\mu)=F$, and constants $c_1,c_2>0$ such that
\begin{equation} \label{eq:mucond}
c_1 h(r) \leq \mu(B(x,r)\cap F) \leq c_2 h(r), \qquad x\in F, \quad 0<r< 1.
\end{equation} 
\end{defn}
Clearly, for $d\in (0,n]$, if $F$ is a $d$-set then it is an $h$-set with $h(r)=r^d$ and $\mu=\cH^d|_F$. Conversely, if $F$ is an $h$-set with $h(r)=r^d$, for some $d\in (0,n]$, so the measure $\mu$ satisfies \eqref{eq:mucond} with $h(r)=r^d$, then $\mu$ is equivalent to $\cH^d|_{F}$ (e.g., \cite[Theorem 3.4]{Triebel97FracSpec}), so that $F$ is a $d$-set. 
\begin{defn}[Uniformly porous, {\cite[p.~85]{Tri08}}] \label{def:up} We say that a closed set $F\subset \R^n$ is {\em uniformly porous} if it is porous (Definition \ref{def:Porous}) and is an $h$-set, for some Radon measure $\mu$ and $h:[0,1]\to[0,1]$ that satisfy the conditions of Definition \ref{def:hset}.
\end{defn}
Since $d$-sets, with $0<d<n$, are porous \cite[Proposition 4.3]{Cae02} and, as noted above, are examples of $h$-sets, they are also examples of uniformly porous sets. That $d$-sets are porous can be seen also from the following extension of \cite[Proposition 4.3]{Cae02}. This is taken from \cite[Proposition 9.18]{Tri01}, except for the last sentence which is a straightforward exercise.
\begin{prop}[{\cite[Proposition 9.18]{Tri01}}] \label{prop:hsetp} Suppose that $F\subset \R^n$ is compact and is an $h$-set.  Then $F$ is porous (and so uniformly porous) if and only if there are positive numbers $c$ and $\lambda\in (0,n)$ such that
\begin{equation} \label{eq:po}
h(2^{-j}) \leq c2^{(n-\lambda)\ell} h(2^{-j-\ell}), \qquad j\in \N_0, \quad \ell\in \N_0.
\end{equation}  
This last condition is equivalent, since $h$ is increasing, to the requirement that, for some $C>0$ and the same $\lambda>0$,
\begin{equation} \label{eq:po2}
\frac{h(R)}{h(r)} \leq C\left(\frac{R}{r}\right)^{n-\lambda}, \qquad 0<r\leq R\leq 1.
\end{equation}
\end{prop}

\section{Bounded sets $\Gamma$ for which $C_0^\infty(\Gamma^\circ)$ is not dense in $H^s_{p,\Gamma}$}
\label{sec:Appendix2}

In this appendix we construct, for $n\geq 2$, examples of bounded sets $\Gamma\subset \R^n$, with $\Gamma^\circ=(\overline{\Gamma})^\circ$ and $\overline{\Gamma^\circ}=\overline{\Gamma}$, for which $\tH^{s}_p(\Gamma^\circ)  \subsetneqq H^{s}_{p,\Gamma}$ for certain ranges of $s$ and $p$, in particular for $s=-1/2$ and $p=2$ as promised in Remark \ref{rem:fails}. In detail, we construct, for $1<p<\infty$ and $0<s\leq n/p$, a bounded open $n$-set $\Omega=\overline{\Omega}^\circ$ with $\tH^s_p(\Omega)\subsetneqq H^s_{p,\overline{\Omega}}$ and then deduce, by applying \cite[Lemma 4.14]{caetano2019density}, that $U:= \overline{\Omega}^c$ satisfies $U=\overline{U}^\circ$ and $\tH^{-s}_p(U)\subsetneqq H^{-s}_{p,\overline{U}}$, with $\overline{U}^c=\Omega$ an open $n$-set. 

Theorem \ref{thm:exnset} that follows is our main result which leads to the above examples. This theorem shows, additionally, that \eqref{eq:-s<0} is sharp in one direction, precisely that \eqref{eq:-s<0} does not hold for all open $n$-sets $\Omega$ if the range $-1\leq s\leq 0$ is extended to $-1\leq s\leq \epsilon$, for some $\epsilon>0$.
Our construction of $\Omega$ in the proof of Theorem \ref{thm:exnset} is, in part, inspired by that in the proof of \cite[Theorem~11.5.5]{AdHe} where, for $m\in \N$ and $1<p\leq n$, a compact set $\Gamma\subset \R^n$ is constructed such that $\overline{\Gamma^\circ}=\Gamma$ and $\Gamma$ is not $(m,p)$-stable in the language of \cite{AdHe}, equivalently  $\tH^{m}_p(\Gamma^\circ)  \subsetneqq H^{m}_{p,\Gamma}$. Thus, except that  \cite{AdHe} does not aim to construct $\Gamma$ so that $\Gamma^\circ$ is an $n$-set\footnote{\label{foot:AB} But note that Lemma \ref{lem:geom} shows, where $K$ is the compact set constructed in the proof of \cite[Theorem~11.5.5]{AdHe}, that $\overline{K^\circ}=K$, as observed previously above Lemma 3.20 in \cite{ChaHewMoi:13}. It also shows that if the parameters $r_k$ and $R_k$ in that proof are chosen with $r_k\leq R_k/3$, $k\in \N$, then $K$ is an $n$-set and $K^\circ$ an open $n$-set.}, \cite[Theorem~11.5.5]{AdHe} establishes our result for the special case that $s\in \N$ with $1\leq s<n$ and $1<p\leq n/s$, but not for any non-integral values of $s$. 
If $s\in \N$ and $p>n$ then every compact set $\Gamma\subset \R^n$ is $(s,p)$-stable, see \cite[Theorem 2.9(a)]{Po:72}. To the best of our knowledge, Theorem \ref{thm:exnset}, applied with $0<s<1$, also provides the first examples of compact sets $\Gamma=\overline{\Omega}$ that are not $(s,p)$-stable with $p>n$.\footnote{The paragraph above Theorem 11.5.5 in \cite{AdHe} might be read as suggesting that every compact set $\Gamma$ is $(s,p)$-stable, {\em for every $s>0$}, provided $p>n$. But we think that Adams and Hedberg intended the additional condition $s\in \N$ in this paragraph; that would be consistent with Remark 1 after \cite[Theorem 11.5.4]{AdHe}.}

\begin{thm} \label{thm:exnset}
Suppose that $n\in \N$, $n\geq 2$, $1<p<\infty$ and $0<s\leq n/p$. Then there exists a bounded open $n$-set $\Omega\subset \R^n$ such that  $\Omega=\overline{\Omega}^\circ$ and $\tH^s_p(\Omega)\subsetneqq H^s_{p,\overline{\Omega}}$.
\end{thm}
As noted above, this theorem has the following corollary which, setting $s=1/2$ and $p=2$, provides the example promised in Remark \ref{rem:fails}. This result also shows that \eqref{eq:s>0} is sharp in one direction, precisely that \eqref{eq:s>0} does not hold for all open $n$-sets $\Omega$ if the range $0\leq s\leq 1$ is extended to $-\epsilon\leq s\leq 1$, for some $\epsilon>0$. 
\begin{cor} \label{cor:exnset} 
Suppose that  $n\in \N$, $n\geq 2$, $1<p<\infty$ and $0<s\leq n/p'$. Then, where $\Omega\subset \R^n$ is the open $n$-set of Theorem \ref{thm:exnset} read with $p'$ instead of $p$, $U:=\overline{\Omega}^c$ has the properties that  $U=\overline{U}^\circ$ and $\tH^{-s}_p(U)\subsetneqq H^{-s}_{p,\overline{U}}$, with $\overline{U}^c$ an open $n$-set. For $R>\max_{x\in \overline{\Omega}}|x|$ the bounded domain $U_R:= \{x\in U:|x|<R\}$ has the same properties. 
\end{cor}
We defer the proofs of the above theorem and corollary, and first prove a lemma that captures the largest part of the construction of the set $\Omega$ that will be used in the proof of Theorem \ref{thm:exnset}. Recall that $B(x,r)=\{y\in \R^n:|x-y|\leq r\}$, for $x\in \R^n$ and $r>0$. Let $B^\circ(x,r):= (B(x,r))^\circ$ denote the corresponding open ball, i.e., $B^\circ(x,r)=\{y\in \R^n:|x-y|< r\}$, and let $B_r^\circ:= B^\circ(0,r)$.
\begin{lem} \label{lem:geom}
Suppose that $\Gamma\subset B_1^\circ$ is a compact $d$-set with $0<d<n$. Choose a sequence of points $(x_k)_{k\in \N}\subset \Gamma$ and a positive sequence $(R_k)_{k\in \N}\subset (0,\infty)$ such that 
\begin{equation} \label{eq:ni}
B(x_k,R_k)\subset B_1^\circ, \quad  B(x_k,R_k)\cap B(x_{k'},R_{k'})=\emptyset, \qquad  k,k'\in \N, \quad k'\neq k,
\end{equation}
and such that 
\begin{equation} \label{eq:dense}
\Gamma\cap\bigcup_{k=1}^\infty B(x_k,R_k) \mbox{ is dense in $\Gamma$.}
\end{equation}
Choose a sequence $(r_k)_{k\in \N}\subset (0,\infty)$ such that $r_k<R_k$, $k\in \N$, and set
$$
\Omega := B_1^\circ \setminus \overline{\bigcup_{k=1}^\infty B(x_k,r_k)}.
$$
Then $\Omega = (\overline{\Omega})^\circ$. If, moreover,
\begin{align}
 \label{eq:rkRk}
 r_k\leq R_k/3,
 \qquad k\in \N,
 \end{align}
 then $\Omega$ is an $n$-set.
\end{lem} 
\begin{proof} 
To see, first of all, that there exist sequences $(x_k)$ and $(R_k)$ satisfying \eqref{eq:ni} and \eqref{eq:dense}, let $Y:=\{y_k:k\in \N\}$ be any dense subset of $\Gamma$ (e.g., $Y=\cup_{j\in \N}Y_j$, where $Y_j\subset \Gamma$ is the set of centres of the balls of a finite cover of $\Gamma$ by balls of radius $1/j$). Note that $Y$ is infinite since $\dimH(\Gamma)>0$, and we may assume, without loss of generality, that the elements of the sequence $(y_k)_{k\in \N}$ are distinct. Using the notation $\Gamma_k:= \Gamma\setminus \cup_{j=1}^{k} B(x_j,R_j)$, for $k\in \N$, choose $x_1=y_1$ and, for $k=1,2,\ldots$, first choose $R_k>0$ so that $\Gamma_k\neq\emptyset$, in which case $\Gamma\cap B(x,\epsilon)\subset \Gamma_k$, for some $x\in \Gamma$ and $\epsilon>0$, so that, by \eqref{eq:dset}, $\Gamma_k$ contains infinitely many points of $\Gamma$ and so infinitely many elements of $Y$, and then choose $x_{k+1}=y_{i}$, where $i$ is the smallest natural number such that $y_{i}\in \Gamma_k$. With this construction $Y\subset \cup_{k=1}^\infty B(x_k,R_k)$, so that \eqref{eq:dense} holds.

The definitions in the statement of the lemma imply that $R_k<1$, for $k\in \N$, $R_k\to 0$ as $k\to\infty$, and $\Omega=K^\circ$, where
\begin{equation} \label{eq:K}
K:= \overline{B_1^\circ} \setminus \bigcup_{k=1}^\infty B^\circ(x_k,r_k).
\end{equation}
 Let
 $\gamma$ denote the set of  limit points of the set $S:=\{x_k:k\in \N\}$, noting that $\gamma\cap S=\emptyset$, indeed $\gamma \cap B(x_k,r_k)=\emptyset$ for $k\in \N$. Note that $K$ and
  $\gamma\subset \Gamma$  are compact and that, if $(r_k')_{k\in \N}\subset (0,\infty)$ and $r_k'\to 0$ as $k\to\infty$, then
\begin{equation} \label{eq:sumeq}
  \overline{\bigcup_{k=1}^\infty B(x_k,r'_k)} = \gamma \cup \bigcup_{k=1}^\infty B(x_k,r'_k).
\end{equation}
Thus
\begin{equation} \label{eq:Omega}
\Omega = B_1^\circ\setminus \left(\gamma \cup \bigcup_{k=1}^\infty B(x_k,r_k)\right).
\end{equation}
It follows from \eqref{eq:K}--\eqref{eq:Omega} that 
$$
K = \Omega \cup \partial B^\circ_1 \cup \gamma \cup \bigcup_{k=1}^\infty \partial B(x_k,r_k).
$$
From this it follows that   $\overline{\Omega}=K$. For certainly $\overline{\Omega}\subset K$ and $\Omega \cup \partial B^\circ_1\subset \overline{\Omega}$.
 Further, by \eqref{eq:ni} and the definition of $\Omega$,
\begin{equation} \label{eq:ann}
B^\circ(x_k,R_k)\setminus B(x_k,r_k) \subset \Omega, \qquad k\in \N,
\end{equation}
so that
$\partial B(x_k,r_k)\subset \overline{\Omega}$, for $k\in \N$. Further, if $x\in \gamma$ then there exists a subsequence $(x_{k_j})$ of the sequence $(x_k)$ such that $x_{k_j}\to x$ as $j\to\infty$, and then, where $y_j$ denotes any point in $\partial B(x_{k_j},r_{k_j})$, $y_j\in \overline{\Omega}$ and $y_j\to x$ as $j\to\infty$, so that $x\in \overline{\Omega}$.
Thus $\overline{\Omega}=K$, so that $\Omega=\overline{\Omega}^\circ$. 

To prove that $\Omega$ is an $n$-set under the additional assumption \eqref{eq:rkRk}, we first note that it suffices to prove the $n$-set condition \eqref{eq:nset} for $0<r\leq r_0$ for some $0<r_0<1$, and then it follows for all $0<r\leq 1$ with a modified constant. Hence we can assume that $0<r\leq 1$ is replaced by $0<r\leq r_0$ in the $n$-set condition \eqref{eq:nset}, where 
$$
r_0:={\rm dist}\left(\overline{{J}},\partial B^\circ_1\right)/3, \quad \mbox{with } \quad J :=  \bigcup_{k=1}^\infty B(x_k,R_k).
$$
 
Now let $x\in\Omega$ and $0<r\leq r_0$. If $B(x,r)$ intersects $\partial B^\circ_1$ then it doesn't intersect $\overline{J}$, 
and in that case the $n$-set condition follows from the fact that $B^\circ_1$ is an open $n$-set. On the other hand, if $B(x,r)$ does not intersect $\partial B^\circ_1$ (so that $B(x,r)\subset B^\circ_1$) then we claim that \eqref{eq:rkRk} implies that
\begin{align}
\label{eq:nsetcond}
|B(x,r)\cap\Omega|\geq \frac{c_n}{2}r^n,
\end{align}
where $c_n$ is the volume of the unit ball in $\R^n$. To prove this we note, using \eqref{eq:Omega} and that $|\gamma|=0$ as $\gamma \subset \Gamma$ and $\Gamma$ is a $d$-set with $d<n$, that .
\[ |B(x,r)\cap \Omega| = |B(x,r)\cap (B^\circ_1\setminus J)| +
\sum_{k=1}^\infty |B(x,r)\cap (B(x_k,R_k)\setminus B(x_k,r_k))|,\]
and also that
$
B(x,r)\cap (B^\circ_1\setminus J)=B(x,r)\cap J^c$.
Then we claim that 
\begin{align}
\label{eq:ballcond}
|B(x,r)\cap (B(x_k,R_k)\setminus B(x_k,r_k))| \geq \frac{1}{2}|B(x,r)\cap B(x_k,R_k)|, \qquad k\in \N,
\end{align}
from which it follows 
that
\begin{eqnarray*} |B(x,r)\cap \Omega| &\geq & |B(x,r)\cap J^c| +
\frac{1}{2}\sum_{k=1}^\infty |B(x,r)\cap B(x_k,R_k)|\\
& = &  |B(x,r)\cap J^c| +
\frac{1}{2}|B(x,r)\cap J| \geq \frac{1}{2}|B(x,r)|=\frac{c_n}{2}r^n. 
\end{eqnarray*}
To verify \eqref{eq:ballcond} we consider three cases, illustrated in Figure \ref{fig:balls}:
\begin{figure}[t!]
\centering
\subfl{Case 1}{\includegraphics[height = 43mm]{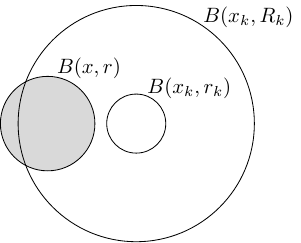}}
\hspace{1mm}
\subfl{Case 2}{\includegraphics[height = 43mm]{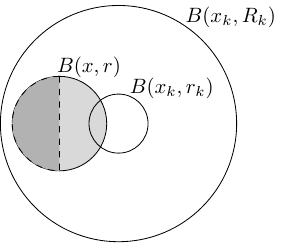}}
\hspace{1mm}
\subfl{Case 3}{\includegraphics[height = 43mm]{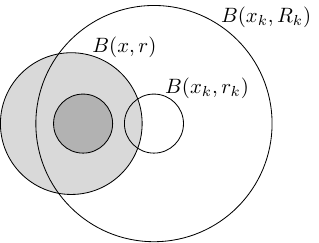}}
\caption{The three cases considered in the proof of the $n$-set property in Lemma \ref{lem:geom}.}
\label{fig:balls}
\end{figure}
\begin{itemize}
\item Case 1: If $B(x,r)\cap B(x_k,r_k)= \emptyset$ then trivially  \eqref{eq:ballcond} holds.
\item Case 2: If $B(x,r)\cap B(x_k,r_k)\neq\emptyset$ and $B(x,r)\subset B(x_k,R_k)$ then $|B(x,r)\cap B(x_k,R_k)|=|B(x,r)|$, and since $x\in B(x_k,r_k)^c$ at least half of $B(x,r)$ lies in $B(x_k,R_k)\setminus B(x_k,r_k)$, so \eqref{eq:ballcond} holds.
\item Case 3: If $B(x,r)\cap B(x_k,r_k)\neq\emptyset$ and $B(x,r)\cap B(x_k,R_k)^c\neq \emptyset$ then by \eqref{eq:rkRk} the set $B(x,r)\cap (B(x_k,R_k)\setminus B(x_k,r_k))$ contains a ball of radius $r_k$. Hence $|B(x,r)\cap (B(x_k,R_k)\setminus B(x_k,r_k))|\geq |B(x,r)\cap B(x_k,r_k)|$, so that 
\begin{eqnarray*}
|B(x,r)\cap B(x_k,R_k)| & = & |B(x, r)\cap B(x_k, r_k)|+|B(x, r)\cap  (B(x_k, R_k)\setminus B(x_k, r_k))|\\
&\leq & 
2|B(x,r)\cap (B(x_k,R_k)\setminus B(x_k,r_k))|,
\end{eqnarray*}
which implies \eqref{eq:ballcond}.
\end{itemize}
\end{proof}

\begin{proof}[Proof of Theorem \ref{thm:exnset}]
For $E\subset \R^n$ let $\Cp_{s,p}(E)$ denote the $(s,p)$-capacity of $E$ as defined in \cite[Definition 2.2.6]{AdHe}. For $s>0$ and $1<p<\infty$, let $\varphi_{s,p}(r):= \Cp_{s,p}(B^\circ_r)$, for $r>0$. Note that  $\varphi_{s,p}$ is non-decreasing on $(0,\infty)$. Further, for $s>0$ and $1<p\leq n/s$, $\varphi_{s,p}(r)\to 0$ as $r\to 0$; see \cite[Props.~5.1.2, 5.1.4]{AdHe} for the case $sp<n$, \cite[Proposition~5.1.3]{AdHe} for the case $sp=n$.

Suppose now that $n\geq 2$, $1<p<\infty$, and $0<s\leq n/p$. Clearly $\tH^s_p(\Omega)\subset H^s_{p,\overline{\Omega}}$ for all domains $\Omega$.
An ingredient in demonstrating, for particular sets $\Omega$, that 
\begin{equation}\label{eq:neq}
\tH^s_p(\Omega)\neq H^s_{p,\overline{\Omega}}
\end{equation}
is the characterisation, in the case that $\Omega$ is a bounded domain with $\Omega=\overline{\Omega}^\circ$, that \eqref{eq:neq} holds if and only if $\overline{\Omega}$ is not $(s,p)$-stable in the sense of \cite[Definition 11.5.2]{AdHe} (note that, if $E=\overline{\Omega}$ and $\overline{\Omega}^\circ=\Omega$, then Adams and Hedberg's $L_0^{s,p}(E^\circ)$ is our $\tH^s_p(\Omega)$ and their $L_0^{s,p}(E)$ is our $H^s_{p,\overline{\Omega}}$; see \cite[Equations (11.5.2), (9.1.1), Definition 11.5.1]{AdHe}). Importantly, this in turn holds, by the comments below \cite[Theorem~11.5.4]{AdHe} (cf.~\cite[Theorem 2.7]{Po:72}),
if
\begin{equation} \label{eq:cap}
\Cp_{s,p}(G\setminus\overline{\Omega})<\Cp_{s,p}(G\setminus\Omega),
\end{equation} 
for some domain $G\subset \R^n$.   (Indeed, if $0<s\leq 1$ then, applying additionally \cite[Theorem~11.5.4]{AdHe}, we see that $\overline{\Omega}$ is not $(s,p)$-stable if and only if \eqref{eq:cap} holds for some domain $G\subset \R^n$.)

To make use of these results, define $\Omega$ as in Lemma \ref{lem:geom}, choosing the parameters as follows. First choose the $d$-set $\Gamma$, with 
\begin{equation} \label{eq:drange}
n-sp<d<n.
\end{equation}
Since $\Gamma$ is a $d$-set, $\cH^d(\Gamma)>0$ and there exists $c_2>0$ such that 
\begin{equation} \label{eq:rhdset}
\cH^d(\Gamma\cap B(x,r))\leq c_2 r^d, \qquad x\in \Gamma, \quad 0<r\leq 1.
\end{equation} 
Next, choose the points $(x_k)_{k\in \N}\subset \Gamma$ and the radii $(R_k)_{k\in \N}\subset (0,\infty)$ as in Lemma \ref{lem:geom} but with the additional constraint that   
\begin{equation} \label{eq:r1}
c_2\sum_{k=1}^\infty R_k^d\ <\cH^d(\Gamma).
\end{equation}
To see the point of this constraint, note that, by \eqref{eq:dense} and \eqref{eq:sumeq},
$$
\Gamma = \overline{\Gamma\cap\bigcup_{k=1}^\infty B(x_k,R_k)} = \gamma \cup \left(\Gamma\cap\bigcup_{k=1}^\infty B(x_k,R_k)\right),
$$
where $\gamma$ is as defined in the proof of Lemma \ref{lem:geom}.
Thus, using \eqref{eq:rhdset} and \eqref{eq:r1}, and recalling from the proof of Lemma \ref{lem:geom} that $R_k<1$ for $k\in \N$,
$$
\cH^d(\gamma) \geq \cH^d(\Gamma)-\sum_{k=1}^\infty \cH^d(\Gamma \cap B(x_k,R_k)) \geq \cH^d(\Gamma)-c_2\sum_{k=1}^\infty R_k^d>0,
$$
so that $\dimH(\gamma)=d$. Since $d$ satisfies \eqref{eq:drange} it follows (see, e.g., \cite[Theorem~3.8(ii)]{HewMoi:15}, noting that the capacity $\overline{\mathrm{cap}}_{s,p}$ in \cite{HewMoi:15} coincides with our $\Cp_{s,p}$ by \cite[Definition 3.1, Remark 3.2]{HewMoi:15}) that $\Cp_{s,p}(\gamma)>0$. Finally, choose the radii $(r_k)_{k\in \N}$ so that $0<r_k\leq R_k/3$ and so that
\begin{equation} \label{eq:rconst}
\sum_{k=1}^\infty \varphi_{s,p}(r_k)< \Cp_{s,p}(\gamma).
\end{equation}

With these choices of the parameters defining $\Omega$, it follows by Lemma \ref{lem:geom} that $\Omega$ is an $n$-set and that $\Omega =(\overline{\Omega})^\circ$. Further, from the proof of the lemma we have that $\overline{\Omega} = K$, where $K$ is defined by \eqref{eq:K}. To complete the proof, set $G:= B^\circ_1$, so that $G$ is open. Then $G\setminus \Omega \supset \gamma$, so that  $\Cp_{s,p}(G\setminus \Omega) \geq  \Cp_{s,p}(\gamma)$. On the other hand,
$$
G\setminus \overline{\Omega}=B^\circ_1\setminus K =  \bigcup_{k=1}^\infty B^\circ(x_k,r_k)
$$
so that
$$
\Cp_{s,p}(G\setminus \overline{\Omega})\leq  \sum_{k=1}^\infty \Cp_{s,p}(B^\circ(x_k,r_k)) = \sum_{k=1}^\infty \varphi_{s,p}(r_k)< \Cp_{s,p}(\gamma),
$$
by \eqref{eq:rconst}.
Thus \eqref{eq:cap} holds, which implies \eqref{eq:neq}.
\end{proof}

\begin{proof}[Proof of Corollary \ref{cor:exnset}]
Suppose that  $n\in \N$, $n\geq 2$, $1<p<\infty$ and $0<s\leq n/p'$. 
Let $\Omega\subset \R^n$ be the open $n$-set of Theorem \ref{thm:exnset} read with $p'$ instead of $p$, with $\overline{\Omega}^\circ = \Omega$ and $\tH^s_{p'}(\Omega)\subsetneqq H^s_{p',\overline{\Omega}}$. Let $U:=\overline{\Omega}^c$, and
note that $\overline{\Omega}^\circ = \Omega$ implies that $\overline{U}=\Omega^c$ and $(\overline{U})^\circ=U$. Clearly, $\tH^{-s}_p(U)\subset H^{-s}_{p,\overline{U}}$, but, applying  \cite[Lemma 4.14]{caetano2019density}, we see that  $\tH^{-s}_p(U)\neq H^{-s}_{p,\overline{U}}$.

Suppose now that $R>\max_{x\in \overline{\Omega}}|x|$ and let $U_R:=U\cap B^\circ_R$. Clearly, $\overline{U_R}^\circ = U_R$ and $\overline{U_R}^c=\Omega\cup (\R^n\setminus B(0,R))$, which is the union of two open $n$-sets and so an open $n$-set. Also,  $\tH_p^{-s}(U_R)\subset H^{-s}_{p,\overline{U_R}}$. To see that $\tH_p^{-s}(U_R)\neq H^{-s}_{p,\overline{U_R}}$, suppose that  $\tH_p^{-s}(U_R)= H^{-s}_{p,\overline{U_R}}$ and $u\in H^{-s}_{p,\overline{U}}\setminus \tH^{-s}_p(U)$. Choose $\chi\in C_0^\infty(\R^n)$ such that $\supp(\chi)\subset B^\circ_R$ and $\chi=1$ in a neighbourhood of $\overline{\Omega}$. Then $\chi u\in  H^{-s}_{p,\overline{U_R}}=\tH_p^{-s}(U_R)\subset \tH_p^{-s}(U)$ and $\supp((1-\chi)u)\subset U$, so that also $(1-\chi) u\in \tH_p^{-s}(U)$, so that $u\in \tH_p^{-s}(U)$, a contradiction.
\end{proof}



\section{{A direct proof of membership of the class $\cD^t$
}}
\label{sec:Appendix3}

{In 
Corollary \ref{cor:newDt} 
we proved that if $\Gamma\subset\R^n$ is an $n$-attractor and $\Omega$ is either $\Gamma^\circ$ or $\Gamma^c$ then $\Omega$ is a member of the 
class $\cD^t$ for all $0<t<n-\dimH(\partial\Gamma)$. The proof relied on 
{Theorem} 
\ref{lem:newDt} and the equality of the Hausdorff and Aikawa/Assouad dimension of $\partial\Gamma$ in this case. 
In Remark \ref{rem:direct} we claimed that using more direct arguments one can prove a slightly weaker statement (which is nonetheless sufficient for our later arguments), namely that $\Omega$ belongs to $\cD^t$ for 
$0<t<n-d$, where $d\in [n-1,n)$ is the dimension of a $d$-set containing $\partial\Gamma$ (the existence of which is known by Proposition \ref{prop:nset}\eqref{x}). In this appendix we provide details of these arguments.}


\begin{lem}
	\label{lem:integral}
	Let $\emptyset\not=A\subset \R^n$ be measurable, bounded and such that $\partial A$
	is contained in a compact $d$-set for some $0\leq d<n$. Then we
	have for all $0\leq t<n-d$ that
	\[
	\int_{A}{\rm dist}(y,\partial A)^{-t}\,\rd y<\infty.
	\]
\end{lem}

\vspace*{-4ex}

\begin{proof}
	We split the above integral as
	\[
	\int_{A\cap(\partial A)^{(1)}}{\rm dist}(y,\partial A)^{-t}\,\rd y\;+\;\int_{A\setminus(\partial A)^{(1)}}{\rm dist}(y,\partial A)^{-t}\,\rd y,
	\]
	where, 
	for $E\subset \R^n$ and $\epsilon>0$, $E^{(\epsilon)}:=\{y\in \R^n:\dist(y,E)<\epsilon\}$ stands for the $\epsilon$ neighbourhood of $E$.
	Clearly, the second integral is bounded above by
	$|A|$ and the first one is bounded above by
	\begin{equation}
		\int_{(\partial A)^{(1)}}{\rm dist}(y,\partial A)^{-t}\,\rd y.\label{eq:partial Omega}
	\end{equation}
	Therefore we just have to prove that the integral in \eqref{eq:partial Omega} is finite
	whenever $0\leq t<n-d$. The case $t=0$ being trivial, we assume
	$0<t<n-d$ in what follows.
	
	Let $B$ be the compact $d$-set, with $0\leq d<n$, containing $\partial A$
	mentioned in the statement of this lemma. We clearly have, for any
	$y\in\R^n$, that ${\rm dist}(y,\partial A)\geq{\rm dist}(y,B)$, from
	which it follows that $(\partial A)^{(1)}\subset B^{(1)}$ and that the integral
	in (\ref{eq:partial Omega}) is bounded above by
	\begin{equation}
		\int_{B^{(1)}}{\rm dist}(y,B)^{-t}\,\rd y.\label{eq:Gamma'}
	\end{equation}
	
	We finish our proof by showing that the latter integral is finite.
	We follow essentially the argument in the first part of the proof
	of \cite[Proposition 16.5]{Tri01}, but we give the details (which are omitted
	there). 
	
	We start by taking advantage of the following consequence of \cite[Remark 2.9]{CC08}
	(see also \cite[Lemma 2.17]{Moura:01}): each neighbourhood $B^{(r^{-1})}$ of $B$ intersects $\approx r^{d}$ cubes of any given regular tessellation
	of $\R^n$ by cubes of sides parallel to the axes and side length $r^{-1}$,
	for any given $r\geq1$, with equivalence constants independent of
	$r$ and of the tessellation. In the above reference the case $d=0$
	is not included, but the result for this case is a trivial one, since a
	compact $0$-set must be finite.
	
	Given such an $r$ and such a tessellation, and writing 
	\[
	B^{(r^{-1})}\subset\bigcup_{i\in I}Q_{i},
	\]
	where $I$ indexes the cubes $Q_{i}$ of the given tessellation which
	have non-void intersection with $B^{(r^{-1})}$, we then have that
	\[
	|B^{(r^{-1})}|\leq\sum_{i\in I}|Q_{i}|=\#I\times r^{-n}\leq cr^{d-n},
	\]
	where $\#$ denotes the cardinality of a set, from which the finiteness of (\ref{eq:Gamma'}) follows, using also
	the fact that $|B|=0$, and where $\widehat{B}_{j}:=B^{(2^{-j/(n-d)})}\setminus B^{(2^{-(j+1)/(n-d)})}$,
	$j\in\No$,
	\begin{align*}
		\int_{B^{(1)}}{\rm dist}(y,B)^{-t}\,\rd y  &=\int_{B^{(1)}\setminus B}{\rm dist}(y,B)^{-t}\,\rd y
		=\sum_{j=0}^{\infty}\int_{\widehat{B}_{j}}{\rm dist}(y,B)^{-t}\,\rd y\\
		& \leq\sum_{j=0}^{\infty}|\widehat{B}_{j}|2^{(j+1)t/(n-d)}
		\leq\sum_{j=0}^{\infty}|B^{(2^{-j/(n-d)})}|2^{(j+1)t/(n-d)}\\
		& \leq c\sum_{j=0}^{\infty}2^{-j}2^{jt/(n-d)}<\infty.
	\end{align*}

\vspace*{-6ex}

\end{proof}

Before stating the main result of this section in Lemma \ref{lem:Dt}, we first present an elementary lemma, which is a slight sharpening of \cite[Lemma\ 9.2]{Fal}. 
Specifically, the statement of \cite[Lemma 9.2]{Fal} assumes that the $\{E_j\}$ are open and disjoint, but the proof of \cite[Lemma\ 9.2]{Fal}, 
which we repeat here for convenience, 
uses neither the openness nor the disjointness of the $\{E_j\}$, but only requires the $\{E_j\}$ to contain interior balls $\{\tilde{B}_j\}$ that are disjoint, 
hence our more general statement. 
\begin{lem}[{\cite[Lemma\ 9.2]{Fal}}]
	\label{lem:Intersection}
	Let $r>0$ and $0<a_1\leq a_2$ be given, and let $\{E_j\}$ be a collection of subsets of $\R^n$ such that each $E_j$ contains an open ball $\tilde{B}_j$ of radius $a_1r$ and is contained in a closed ball $B_j$ of radius $a_2r$, such that the balls $\{\tilde{B}_j\}$ are pairwise disjoint. Then any closed ball $B$ of radius $r$ intersects at most $(1+2a_2)^na_1^{-n}$ of the closures $\{\overline{E_j}\}$.	
\end{lem}

\vspace*{-4ex}

\begin{proof}
	If $\overline{E_j}$ intersects $B$ then $\overline{E_j}$ is contained in the closed ball of radius $(1+2a_2)r$ concentric with $B$. Suppose that $q$ of the sets $\{\overline{E_j}\}$ intersect $B$. Summing the volumes of the corresponding interior balls $\tilde{B}_j$, it follows that $qa_1^n r^n \leq (1+2a_2)^nr^n$,
	giving the stated bound for $q$.
\end{proof}

\begin{lem}
	\label{lem:Dt}
	Let $\Gamma\subset\R^n$ be an $n$-attractor,  so that $\partial \Gamma$ is contained in a compact $d$-set, with $n-1\leq d<n$, and let $\Omega$ denote either $\Gamma^\circ$ or $\Gamma^c$.
	Then \eqref{eq:D_t} 
	holds, 
	i.e., $\Omega$ belongs to the Frazier-Jawerth class $\cD^t$, 
for all 
$0< t<n-d$. 
\end{lem}

\vspace*{-4ex}

\begin{proof}
	Note that the existence of the compact $d$-set in the statement
	is guaranteed by part \eqref{x} of Proposition \ref{prop:nset}. 
	{Note also that it suffices to prove the claimed result for $\Omega=\Gamma^\circ$, and then the analogous result for $\Omega=\Gamma^c$ follows automatically because \eqref{eq:D_t} only involves $\partial\Omega$ (not $\Omega$ itself), and $\partial(\Gamma^\circ)=\partial(\Gamma^c)$ (by part \eqref{v} of Proposition \ref{prop:nset}).} 
	Therefore we assume henceforth that $\Omega:=\Gamma^\circ$. 
	We claim that by the IFS structure and by Lemma \ref{lem:Intersection} 
	there exists $r_{0}\in(0,1]$ 
	and ${N_0}\in\N$ such that for any $x\in\partial\Omega$
	and any $0<r\leq r_{0}$ there exists $N\in \{1,2,\ldots,N_0\}$ and a collection of domains $E_1,\ldots,E_N$ such that (i) $B(x,r)\subset \cup_{j=1}^N \overline{E_j}$, (ii) for each $j$, $E_j$ is similar to $\Omega$, $\diam(E_j)\leq r$, and $E_j\cap \partial\Omega=\emptyset$, and (iii) for $j\neq k$, $|\overline{E_j} \cap \overline{E_k}|=0$. 
	To see this, we recall from the proof of Proposition \ref{prop:nset}\eqref{x} that there exist $\ell\in\N$ and $\bm\in I_\ell$ such that $U:=s_\bm^{-1}(\Omega)$ is an open neighbourhood of $\Gamma=\overline{\Omega}$ such that $U$ is similar to $\Omega$ and $\overline{U}$ equals the union of the closures of finitely many similar copies of $\Omega$ (these are the sets $\{\Omega^\bm_{\bm'}\}_{\bm'\in I_\ell}$), none of which intersect $\partial\Omega$, and whose closures intersect pairwise in sets of zero measure. 
	{Furthermore, it is easy to see that one can choose $\ell$ and $\bm$ to ensure that $\dist(\partial\Omega,\partial U)\geq 1$.} 
By construction, the diameter of any of the sets $\Omega^\bm_{\bm'}$ lies in the interval ${[}(\rho_{\rm min}^\ell/\rho_\bm) h_0, (\rho_{\rm max}^\ell/\rho_\bm) h_0{]}$, where $\rho_\bm := \diam(s_\bm(\Gamma))/h_0$. {Provided that} ${0<}r\leq (\rho_{\rm min}^\ell/\rho_\bm) h_0$, one can further decompose each of these sets into a finite collection of sets $E$, each a similar copy of $\Omega$ whose diameter lies in the interval $(\rho_{\rm min}r,r]$ (cf.\ the $L_h$ decomposition introduced in \S\ref{sec:nsets}). Since each such set $E$ is similar to $\Omega$, there exist $a_1,a_2>0$ such that any such $E$ contains an open ball of radius $a_1 r$ and is contained in a ball of radius $a_2 r$.
Then 
	Lemma \ref{lem:Intersection} can be applied to complete the proof of the claim, which holds with $r_0=\min((\rho_{\rm min}^\ell/\rho_\bm) h_0,1)$ and {$N_0$ any integer greater than or equal to} $(1+2a_2)^na_1^{-n}$. 
	
	Thus, for any $t\geq0$ and
	$0<r\leq r_{0}$,
	\[
	\int_{B(x,r)\setminus\partial\Omega}{\rm dist}(y,\partial\Omega)^{-t}\,\rd y\leq\sum_{j=1}^{N}\int_{E_j}{\rm dist}(y,\partial\Omega)^{-t}\,\rd y\leq\sum_{j=1}^{N}\int_{E_j}{\rm dist}(y,\partial E_j)^{-t}\,\rd y,
	\]
	where we used the fact that $E_j\cap \partial\Omega=\emptyset$,
	which implies that ${\rm dist}(y,\partial\Omega)\geq{\rm dist}(y,\partial E_j)$
	for $y\in E_j$. Next we note that, for each $j$
	and each $0\leq t\leq n$, by the similarity of $E_j$
	and $\Omega$, and the fact that $\diam(E_j)\leq r$, and where $h_{0}:={\rm diam}(\Gamma)$,
	\[
	\int_{E_j}{\rm dist}(y,\partial E_j)^{-t}\,\rd y=\left(\frac{{\rm diam}(E_j)}{h_{0}}\right)^{n-t}\int_{\Omega}{\rm dist}(y,\partial\Omega)^{-t}\,\rd y\leq\left(\frac{r}{h_{0}}\right)^{n-t}\int_{\Omega}{\rm dist}(y,\partial\Omega)^{-t}\,\rd y,
	\]
	from which it follows that
	\[
	\sup_{0<r\leq r_{0}}r^{t-n}\int_{B(x,r)\setminus\partial\Omega}{\rm dist}(y,\partial\Omega)^{-t}\,\rd y\leq {N_0}h_{0}^{t-n}\int_{\Omega}{\rm dist}(y,\partial\Omega)^{-t}\,\rd y.
	\]
	We note also that if $r_0<1$ and {$0\leq t\leq n$} then, since $\dist(\partial\Omega,\partial U)\geq 1$ {and $\partial(U\setminus\partial\Omega) = \partial U \cup \partial\Omega$},
	\[
	\sup_{r_0<r\leq 1}r^{t-n}\int_{B(x,r)\setminus\partial\Omega}{\rm dist}(y,\partial\Omega)^{-t}\,\rd y\leq r_0^{t-n}\int_{U\setminus\partial\Omega}{\rm dist}(y,{\partial(U\setminus\partial\Omega}))^{-t}\,rd y.
	\]
	Combining these two estimates with 
	Lemma \ref{lem:integral} (applied with $A=\Omega$ and $A=U\setminus\partial\Omega$ respectively, noting that the latter satisfies the hypotheses of Lemma \ref{lem:integral} because $\partial(U\setminus\partial\Omega) = \partial U \cup \partial\Omega$ and $U$ is similar to $\Omega$) 
	we get (\ref{eq:D_t}) for any $0\leq t<n-d$, 
{and therefore that $\Omega\in \cD^t$ for $0<t<n-d$.}
\end{proof}

\section{{On the definitions of the Aikawa and Assouad dimensions}}
\label{sec:Appendix4}

The notions of the Aikawa and Assouad dimension of a subset of $\R^n$ play an important role in the arguments we make in \S\S\ref{sec:nsets}-\ref{sec:FunctionSpaces}. In this appendix we comment on the equivalence of the different definitions of these notions that appear in the literature that we cite, viz., \cite{LeTu:13,Fraser:14,Fraser:21} in relation to the Assouad dimension,  \cite{Aikawa:91,LeTu:13} in relation to the Aikawa dimension. As in the main text we will denote the Aikawa dimension, as defined by \eqref{eq:IADef}, by $\dim_{\rm A}(F)$. We will denote the Assouad dimension, as defined by \eqref{eq:IASDef} below, by $\dim_{\rm AS}(F)$. It follows from \cite[Theorem 1.1]{LeTu:13} (given the equivalence of the different notions of the Aikawa and Assouad dimensions that we comment on in this appendix) that $\dim_{\rm A}(F)=\dim_{\rm AS}(F)$, for all non-empty $F\subset \R^n$.  

We focus first on the Aikawa dimension $\dim_{\rm A}(F)$ of  a non-empty set $F\subset \R^n$, for which we have given the definition \eqref{eq:IADef}, that
\begin{align}
\label{eq:IADefr}
\dim_{\rm A}(F) := \inf I_{\rm A}(F) \;\mbox{ where  }\; I_{\rm A}(F):=\left\{s>0:\sup_{x\in F}\sup_{r>0}r^{-s}\int_{B(x, r)} \operatorname{dist}(y, F)^{s-n} \, \rd y <\infty\right\}.
\end{align}
This definition agrees with that in Aikawa \cite{Aikawa:91}, made for the case that $F\subset \R^n$ is closed with empty interior\footnote{That is, $\dimA(F)=d(F)$, when $F$ is closed with empty interior, where $d(F)$ is defined for $F$ in this class in the introduction of \cite{Aikawa:91}.}, if we use, as indicated below \eqref{eq:IADef}, the convention that, if $\dist(y,F)=0$, then $\dist(y,F)^{s-n} := +\infty$ if $s<n$, $\dist(y,F)^{s-n}:=1$ if $s=n$. With this convention, $n\in I_A(F)$, so that $\dimA(F)\leq n$, for every $F\subset \R^n$. Further, equality holds, i.e., $\dimA(F)=n$, if $|\overline{F}|>0$, in which case $\{y\in\ B(x,r):\dist(y,F)=0\}$ has positive measure for some $x\in F$ and $r>0$. 

Observe also that if $s>\dimA(F)$ then $s\in I_{\rm A}(F)$, so that either $I_{\rm A}(F)=(\dimA(F),\infty)$ or $I_{\rm A}(F)=[\dimA(F),\infty)$.  This is seen by the following elementary argument. Suppose that $s>\dimA(F)$. Then there exists $s'\in I_A(F)$ such that $\dimA(F)<s'<s$. Noting that for each $x\in F$ and each $r>0$ we have that $\dist(y, F)/r\leq 1$ for all $y\in B(x,r)$, it follows that $(\dist(y, F)/r)^s\leq (\dist(y, F)/r)^{s'}$, which implies that also $s\in I_{\rm A}(F)$.

 For every $r_0\in (0,\infty]$ and non-empty $F\subset \R^n$, let 
\begin{align}
\label{eq:IADefr2}
I^{r_0}_{\rm A}(F):=\left\{s>0:\sup_{x\in F}\sup_{0<r<r_0}r^{-s}\int_{B(x, r)} \dist(y, F)^{s-n} \, \rd y <\infty\right\},
\end{align}
noting that $I^\infty_{\rm A}(F)=I_{\rm A}(F)$. Then the definition of the Aikawa dimension of \cite[Definition 3.2]{LeTu:13}, specialised to the metric space $\R^n$ equipped with Lebesgue measure, is that\footnote{Strictly speaking, in the definition of $I^{\diam(F)}_{\rm A}(F)$ in \cite[Definition 3.2]{LeTu:13}, the closed ball $B(x,r)$ is replaced by the corresponding open ball $B^\circ(x,r)$. But it is easy to see that this change does not affect the membership of $I^{\diam(F)}_{\rm A}(F)$.}
\begin{equation} \label{eq:AikL}
\dimA(F) := \inf I^{\diam(F)}_{\rm A}(F)= \inf \left\{s>0:\sup_{x\in F}\sup_{0<r<\diam(F)}r^{-s}\int_{B(x, r)} \dist(y, F)^{s-n} \, \rd y <\infty\right\},
\end{equation}
for all non-empty $F\subset \R^n$. 
It was stated {without proof} in \cite[\S3.3]{LeTu:13} that the requirement $r<\diam(F)$ is unnecessary in the $\R^n$ setting, so that the definitions \eqref{eq:IADefr} and 
\eqref{eq:AikL} 
coincide. 
{Here we provide a justification of this.  Clearly \eqref{eq:IADefr} and 
\eqref{eq:AikL} coincide if $\diam(F)=\infty$. In the case that $F$ is bounded the coincidence of \eqref{eq:IADefr} and 
\eqref{eq:AikL} holds by the following lemma applied with $r_0=\diam(F)$.}

\begin{lem} \label{lem:equiv} $I^{r_0}_{\rm A}(F)=I_{\rm A}(F)$, for every non-empty bounded $F\subset \R^n$ and $r_0>0$.
\end{lem}

\vspace*{-4ex}

\begin{proof}
Suppose that $\emptyset \neq F\subset \R^n$ is bounded and $r_0>0$. Clearly,  $I_{\rm A}(F)\subset I^{r_0}_{\rm A}(F)$. Conversely, suppose that $t\in I^{r_0}_{\rm A}(F)$, so that there exists $c_{t,n}>0$ such that 
	\begin{equation}\label{1}
\int_{B(x,r)} \mathrm{dist}(y, F)^{t-n} \, \rd y \leq c_{t,n} r^t, \qquad x\in F, \quad 0<r<r_0.
	\end{equation}
To show that also $t\in I_{\rm A}(F)$ we have to show that \eqref{1} also holds, with possibly a different constant $c_{t,n}$, for all $x \in F$ and all $r \geq r_0$. Since, as discussed above the lemma, $[n,\infty)\subset I_{\rm A}(F)$, for every non-empty $F$, we may assume that $0<t<n$.
	
So suppose that $x\in F$ and $r\geq r_0$ and let $J_0 \in \mathbb{N}$ be the smallest number such that
	\[
	2^{J_0-2}\frac{r_0}{3} \geq \diam(F),
	\]
	and $J \in \mathbb{N}\setminus\{1\}$ be the unique number such that
	\[
	2^{J-1} \frac{r_0}{3} \leq r < 2^J \frac{r_0}{3}.
	\]
	Note that 
	\begin{equation}\label{3}
			\text{if $y \notin B\left(x, 2^{j-1} \frac{r_0}{3}\right)$, with $j \geq J_0$, then}
			\dist(y, F) \geq 2^{j-2} \frac{r_0}{3}.
	\end{equation}
	Indeed, if $y \notin B\left(x, 2^{j-1} r_0/3\right)$ and $z \in F$ we have, from $|y - x| - |y - z| \leq |x - z| \leq \diam(F)$,
	that
	\[
	|y - z| \geq |y - x| - \diam(F) > 2^{j-1} \frac{r_0}{3} - 2^{J_0-2} \frac{r_0}{3} \geq 2^{j-2} \frac{r_0}{3},
	\]
	 from which \eqref{3} follows. Note also that
	 the compact set
	\[
	F_{r_0/3} := \left\{ y \in \mathbb{R}^n : \mathrm{dist}(y, F) \leq \frac{r_0}{3} \right\}
	\]
	 is covered by $\bigcup_{x \in F} B^\circ(x, \frac{r_0}{2})$, so there exists $N \in \mathbb{N}$ and $x_1, \ldots, x_N \in F$ such that
	\[
	F_{r_0/3} \subset \bigcup_{i=1}^N B\left(x_i, \frac{r_0}{2}\right).
	\]
	Using these observations we will complete the proof by showing that the integral on the left hand side of \eqref{1} is bounded by $c'_{t,n}r^t$, for some constant $c'_{t,n}$ that depends only on $t$, $n$, $r_0$, and $F$.
	
	Now
	\begin{align*}
	\int_{B(x, r)} \mathrm{dist}(y, F)^{t-n} \rd y &\leq \int_{B\left(x, 2^J \frac{r_0}{3}\right)} \mathrm{dist}(y, F)^{t-n} \rd y\\
	&
	= \int_{B\left(x, 2^J \frac{r_0}{3}\right) \cap F_{r_0/3}} \mathrm{dist}(y, F)^{t-n} \rd y + \int_{B\left(x, 2^J \frac{r_0}{3}\right) \cap F_{r_0/3}^c} \mathrm{dist}(y, F)^{t-n} \rd y\\
&	=: I + II.
	\end{align*}
	Further,
	\[
	I \leq \int_{F_{r_0/3}} \mathrm{dist}(y, F)^{t-n} \rd y \leq \sum_{i=1}^N \int_{B\left(x_i, \frac{r_0}{2}\right)} \mathrm{dist}(y, F)^{t-n} \rd y \leq  N c_{t,n} 2^{-t} r^t \leq N c_{t,n} r^t.
	\]
	We split the analysis of $II$ into the cases: 
	(i) $J \leq J_0$; (ii) $J > J_0$.
In case (i), where $\omega_n$ is the volume of the unit ball in $\R^n$, 
	\begin{align*}
	II \leq \int_{B\left(x, 2^{J_0} \frac{r_0}{3}\right)} \left( \frac{r_0}{3} \right)^{t-n} \rd y = \omega_n \left(2^{J_0} \frac{r_0}{3}\right)^n \left( \frac{r_0}{3} \right)^{t-n} \leq \omega_n 2^{J_0n} 3^{-t} r^t \leq \omega_n 2^{J_0n}  r^t.
	\end{align*}
	In case (ii),
	\begin{align*}
	II &= \int_{B\left(x, 2^{J_0} \frac{r_0}{3}\right) \cap F_{r_0/3}^c} \mathrm{dist}(y, F)^{t-n} \rd y + \sum_{j=J_0+1}^J \int_{\left(B\left(x, 2^j \frac{r_0}{3}\right) \setminus B\left(x, 2^{j-1} \frac{r_0}{3}\right)\right) \cap F_{r_0/3}^c} \mathrm{dist}(y, F)^{t-n} \rd y\\
	&=: III + IV.
	\end{align*}
	Arguing as for $II$ in case (i), $III \leq \omega_n 2^{J_0n} r^t$. On the other hand, using \eqref{3} we get that
	\begin{align*}
	IV &\leq \sum_{j=J_0+1}^J \int_{B\left(x, 2^j \frac{r_0}{3}\right)} \left( 2^{j-2} \frac{r_0}{3} \right)^{t-n} \rd y\\
&	=  2^{2n-2t} \left( \frac{r_0}{3} \right)^t \omega_n \sum_{j=J_0+1}^J 2^{jt} =  2^{2n-2t} \left( \frac{r_0}{3} \right)^t \omega_n \frac{2^{(J+1)t} - 2^{(J_0+1)t}}{2^t - 1}\\
&	\leq \omega_n \frac{2^{2n}}{2^t - 1} 2^{-2t} \left( \frac{r_0}{3}  \right)^t 2^{(J+1)t} \leq \omega_n \frac{2^{2n}}{2^t - 1} r^t.
\end{align*}

	Putting everything together, we have shown that \eqref{1} holds with $c_{t,n}$ replaced by
	$$
	c_{t,n}' := Nc_{t,n} + \omega_n \left(2^{J_0n} + \frac{2^{2n}}{2^t-1}\right),
	$$
	which depends only on $t, n, r_0$ and $F$.
\end{proof}

Let us turn now to the Assouad dimension, for which we adopt the definition from  \cite[\S2.1]{Fraser:21}, that, for non-empty $F\subset \mathbb{R}^n$, the Assouad dimension of $F$ is
\begin{equation} \label{eq:IASDef}
\dim_{\rm AS}(F):= \inf I_{\rm AS}(F) \mbox{ where }
I_{\rm AS}(F) := \left\{s>0:\sup_{x\in F}\sup_{R>0}\sup_{0<r<R} (r/R)^{s}N_r(B(x,R)\cap F)< \infty\right\}.
\end{equation}
In the above equation $N_r(E)$, for $r>0$ and bounded $E\subset \R^n$, is the smallest number of open sets of diameter $\leq r$ required for a cover of $E$. We use in the main text results from  \cite{Fraser:21} (based on this definition), but also results from \cite{Fraser:14} based on what is described in \cite[\S2.2]{Fraser:21} as a ``subtly different definition'', that $\dim_{\rm AS}(F):= \inf I'_{\rm AS}(F)$, where 
\begin{equation} \label{eq:IASDef2}
I'_{\rm AS}(F) := \left\{s>0:\sup_{x\in F}\sup_{0<r<R\leq \rho} (r/R)^{s}N_r(B(x,R)\cap F)< \infty, \mbox{ for some } \rho>0\right\}.
\end{equation}
 It is easy to see that\footnote{Clearly, $I_{\rm AS}(F)\subset I'_{\rm AS}(F)$ for every non-empty $F$. Suppose that $F$ is bounded, $\rho_0:= \diam(F)$, and $s\in I'_{\rm AS}(F)$. Then, where $C_t$ denotes, for $t>0$, the condition $\sup_{x\in F}\sup_{0<r<R\leq t} (r/R)^{s}N_r(B(x,R)\cap F)< \infty$ in \eqref{eq:IASDef2}, the condition $C_\rho$ holds for some $\rho>0$. This implies, by the compactness of $\overline{F}$, that $C_t$ holds for some $t>\rho_0$, which implies that $s\in I_{\rm AS}(F)$.}, if $F$ is bounded, $I_{\rm AS}(F)= I'_{\rm AS}(F)$ so that, as noted in  \cite[\S2.2]{Fraser:21}, these definitions coincide. This is enough for the applications that we make in the main text of the results from \cite{Fraser:14}. These different definitions need not coincide if $F$ is unbounded. For example, as pointed out in  \cite[\S2.2]{Fraser:21}, the definition $\dim_{\rm AS}(F):= \inf I_{\rm AS}(F)$ assigns the set of integers the dimension 1, whereas  $\dim_{\rm AS}(F):= \inf I'_{\rm AS}(F)$ assigns the dimension 0. 

The definition of the Assouad dimension in \cite[Definition 3.1]{LeTu:13}, specialised to the $\R^n$ setting, is a variant on that in \cite{Fraser:21}, namely that $\dim_{\rm AS}(F):= \inf I''_{\rm AS}(F)$, where
\begin{equation} \label{eq:IASDef3}
I''_{\rm AS}(F) := \left\{s>0:\sup_{\stackrel{E\subset F}{\emptyset\neq E \,\mbox{\footnotesize bounded}}}\sup_{0<r<\diam(E)/2} (r/\diam(E))^{s}N'_r(E)< \infty\right\}
\end{equation}
and $N'_r(E)$, for $r>0$ and bounded $E\subset \R^n$, is the smallest number of open balls of radius $\leq r$ required for a cover of $E$. It is easy to see that $I_{\rm AS}(F)=I''_{\rm AS}(F)$, for every non-empty $F\subset \R^n$, since, if $\emptyset\neq E\subset \R^n$ is bounded, then $E\subset B(x,\diam(E))$, for every $x\in E$, so that the supremum taken over all bounded $E\subset F$ in \eqref{eq:IASDef3} can be replaced by one over all $B(x,R)\cap F$ with $x\in F$ and $R>0$, and $N_{2r}(E)\leq N'_{r}(E) \leq N_{r}(E)$. Thus the Assouad definitions of \cite{Fraser:21} and \cite{LeTu:13} coincide, and coincide with the definition in \cite{Fraser:14} when $F$ is bounded. Since \cite[Theorem 1.1]{LeTu:13} proves the equivalence of the Aikawa and Assouad dimensions, we have, moreover, as noted earlier, that $\dimA(F)=\dim_{\rm AS}(F)$, for every non-empty $F\subset \R^n$, with $\dimA(F)$ and $\dim_{\rm AS}(F)$ defined, respectively, by \eqref{eq:IADefr} and \eqref{eq:IASDef}. Recognising that these dimensions coincide we use $\dimA(F)$ in the main text to denote either the Aikawa or the Assouad dimensions, given by \eqref{eq:IADefr} or \eqref{eq:IASDef}.

In the case that $F\subset \R^n$ is bounded, one important result we use, that follows from   \cite[Equation (3.1)]{Fraser:21} and \cite[Proposition 2.8]{Fraser:14}, is that $\dimA(F)=d$ if $F$ is a $d$-set. This result need not hold if $F$ is unbounded. The example $F=x+\alpha\Z^n$, with $x\in \R^n$ and $\alpha>0$, illustrates this. It is a $d$-set with $d=0$ but $\dim_{\rm AS}(F)=n$, since $N_{\alpha/2}(B(x,R)\cap \alpha \Z^n)$ grows asymptotically like $R^n$ as $R\to\infty$. More generally, if $F_0\subset \R^n$ is a compact $d$-set, for some $d\in [0,n)$, and $a>0$ is such that $F_0\subset [-a,a]^n$, then 
\begin{equation} \label{eq:destinf}
F := \sum_{z\in 4a\Z^n} \left(z+F_0\right)
\end{equation}
is a closed $d$-set but, for every $x\in F_0$, $x+4a\Z^n\subset F$ so that $\dimA(F)=\dim_{\rm AS}(F)\geq \dim_{\rm AS}(x+ 4a\Z^n)=n$.

We finish this appendix with a proposition linking the concept of uniform porosity of a closed set $F\subset \R^n$ (Definition \ref{def:up}) to the above notions of the Assouad dimension of $F$. We  use this result in the main text in Remark \ref{rem:multipliers} 
to relate our results on multipliers  of function spaces on rough domains $\Omega$, expressed in terms of $\dimA(\partial \Omega)$, to those in Triebel \cite{Triebel2003,Tri08}, expressed in terms of uniform porosity of $\partial \Omega$. 
\begin{prop} \label{prop:upA} Suppose that $F\subset \R^n$ is compact and uniformly porous, in which case, by Definition \ref{def:up}, $F$ is an $h$-set,  for some Radon measure $\mu$ and $h:[0,1]\to[0,1]$ that satisfy the conditions of Definition \ref{def:hset}, so that \eqref{eq:mucond} holds for some $c_1,c_2>0$, and, by Proposition \ref{prop:hsetp}, \eqref{eq:po2} holds for some $C>0$ and $\lambda\in (0,n)$.  Then $\dimA(F)\leq n-\lambda$.
\end{prop}

\vspace*{-4ex}

\begin{proof}
Suppose that $0<r\leq R\leq 1/3$ and $x\in F$. Let $\gamma_{\rho}(y):= B(y,\rho)\cap F$, for $y\in F$ and $\rho>0$. Since $\cup_{y\in \gamma_R(x)}B(y,r/12)$ is a cover for $\gamma_R(x)$, it follows by the Vitali covering theorem (e.g., \cite[Theorem 1.24]{EG:2015}) that there exists a countable subset $\mathcal{G}\subset \gamma_R(x)$ such that $B(y,r/12)\cap B(z,r/12)=\emptyset$, for $y,z\in \mathcal{G}$ with $y\neq z$, and such that $\gamma_R(x) \subset \cup_{y\in \mathcal{G}}B(y,5r/12)$. But this implies that  $N_r(\gamma_R(x))\leq \# \mathcal{G}$, the cardinality of $\mathcal{G}$ (which must be finite since $F$ is bounded). But also, by \eqref{eq:mucond}, and since $\gamma_{r/12}(y)\subset \gamma_{2R}(x)$, for $y\in \mathcal{G}$,
$$
c_2h(2R) \geq \mu(\gamma_{2R}(x)) \geq  \mu(\cup_{y\in \mathcal{G}}\gamma_{r/12}(y)) = \sum_{y\in \mathcal{G}} \mu(\gamma_{r/12}(y)) \geq  c_1h(r/12) \,\#\mathcal{G},
$$
so that, by \eqref{eq:po2}, 
$$
N_r(\gamma_R(x))\leq \#\mathcal{G} \leq \frac{c_2 h(2R)}{c_1h(r/12)} \leq   \frac{c_2 C}{c_1}\left(\frac{24 R}{r}\right)^{n-\lambda}.
$$
This implies in turn that $n-\lambda\in I'_{\rm AS}(\partial \Omega)$, as defined in  \eqref{eq:IASDef2}, so that the Assouad dimension of $F$, in the sense of \cite{Fraser:14}, is $\leq n-\lambda$. Since $F$ is bounded, it follows, as discussed above, that this Assouad dimension coincides with $\dim_{\rm AS}(F)$, defined by \eqref{eq:IASDef}, which coincides with $\dimA(F)$, defined by \eqref{eq:IADefr}, so that $\dimA(F)\leq n-\lambda$. 
\end{proof}


\vspace*{-2ex}

\bibliography{BEMbib_d=n} 
\bibliographystyle{siam}
\end{document}